\documentclass[oneside,english]{elsartdense}
\usepackage[T1]{fontenc}
\usepackage[latin9]{inputenc}
\usepackage{float}
\usepackage{fancybox}
\usepackage{calc}
\usepackage{amsmath}
\usepackage{graphicx}
\usepackage{amssymb}
\usepackage{esint}

\makeatletter

\newcommand{\lyxline}[1][1pt]{%
  \par\noindent%
  \rule[.5ex]{\linewidth}{#1}\par}
\providecommand{\tabularnewline}{\\}

\newenvironment{lyxlist}[1]
{\begin{list}{}
{\settowidth{\labelwidth}{#1}
 \setlength{\leftmargin}{\labelwidth}
 \addtolength{\leftmargin}{\labelsep}
 }}
{\end{list}}

\makeatother

\makeatother

\usepackage{babel}

\makeatother

\usepackage{babel}

\makeatother

\usepackage{babel}

\begin{document}
\begin{frontmatter}

\title{Parafermi Algebra and Interordinality}

\author{U. Merkel}

\address{D-70569 Stuttgart, Universitätsstr. 38}

\ead{merkel.u8@googlemail.com}
\begin{abstract}
Starting from the observation that for neighboring orders $p=2^{n}-1$,
$p'=2^{n+1}-1$ of the well-known Green's representations of parafermi
algebra there exists a specifiable interordinal relationship, matrices
with similar interordinality properties and intrinsic Catalan structure
are constructed which seem to have a bearing on Euclidean geometry
and applications via Nebe kissing numbers. \end{abstract}
\begin{keyword}
parafermi algebra \sep paraorder \sep Mersenne numbers \sep interordinality
\sep Catalan numbers\sep span parameter \sep stenoscopy \sep kissing
numbers \sep topological/geometrical operator \sep preon model \sep
cardioidic transformation \sep crossing ellipses \sep continued
fractions \sep sine-like and cosine-like structural relationships
\sep positional number systems \sep qphyletics
\end{keyword}
\end{frontmatter}

\section{Introduction -- parafermi operator and root-of-nilpotent sequences}

\setlength{\fboxsep}{1mm}

\noindent Parafermi structures have been studied both in modern quantum
field theory, following H.S. Green, and quantum information theory
\cite{Green98}. The term parafermion is specifically used for the
generalization of a spin-$1/2$ particle (fermion) to spin $p/2$.
Translated to operator language,\begin{equation}
b^{p+1}=\left(b^{+}\right)^{p+1}=0.\end{equation}
 In his original paper \cite{Green53}, Green supplied a $(p+1)\!\times\!(p+1)$
matrix representation for $b$, \begin{equation}
b_{\alpha\beta}=B_{\beta}\delta_{\alpha,\beta+1},\;\;(b^{+})_{\alpha\beta}=B_{\alpha}\delta_{\alpha+1,\beta},\qquad B_{\beta}=\sqrt{\beta(p-\beta+1)},\end{equation}
 which realizes the spin-$p/2$ representation\begin{equation}
\frac{1}{2}[b^{+},b]=\mathrm{diag}\textrm{\Large(}\frac{p}{2},\frac{p}{2}-1,\cdots,-\frac{p}{2}+1,-\frac{p}{2}\textrm{\Large)}\label{eq:b-spin}\end{equation}
 and the characteristic trilinear relations of parafermi algebra\begin{equation}
[[b^{+},b],b]=-2b,\;\:[[b^{+},b],b^{+}]=2b^{+}.\label{eq:trilinear}\end{equation}
 For the least paraorder, the parafermi operator coincides with the
fermi operator $f^{(1)}$ and satisfies the well-known algebra\begin{equation}
\{f^{(1)},(f^{(1)})^{+}\}=1,\:(f^{(1)})^{2}=0=((f^{(1)})^{+})^{2}.\end{equation}
 One fact that seems to have been neglected, or overlooked, is that
those representations, when of order $p=2^{n}-1$ and tensorially
extended by \textbf{1,} are related to those of order $p'=2p+1=2^{n+1}\!-1$
by an operator identity that could be named the Mersennian of parafermi
algebra, for the 17th-century scholar Marin Mersenne who studied the
properties of $2^{n}-1$: \begin{equation}
\frac{1}{2}\textrm{\large\{}b^{(p')},1\!\!\!\boldsymbol{1}^{\otimes n}\otimes b^{(1)}\textrm{\large\}}=b^{(p)}\otimes1\!\!\!\boldsymbol{1}.\label{eq:inter-b}\end{equation}
 For example,\[
\begin{array}{c}
\left\{ b^{(7)},1\!\!\!\boldsymbol{1}^{\otimes2}\otimes b^{(1)}\right\} =\\
\\\left(\begin{array}{cccccccc}
0 & 0 & 0 & 0 & 0 & 0 & 0 & 0\\
\sqrt{7} & 0 & 0 & 0 & 0 & 0 & 0 & 0\\
0 & \sqrt{12} & 0 & 0 & 0 & 0 & 0 & 0\\
0 & 0 & \sqrt{15} & 0 & 0 & 0 & 0 & 0\\
0 & 0 & 0 & 4 & 0 & 0 & 0 & 0\\
0 & 0 & 0 & 0 & \sqrt{15} & 0 & 0 & 0\\
0 & 0 & 0 & 0 & 0 & \sqrt{12} & 0 & 0\\
0 & 0 & 0 & 0 & 0 & 0 & \sqrt{7} & 0\end{array}\right)\left(\begin{array}{cccccccc}
0 & 0 & 0 & 0 & 0 & 0 & 0 & 0\\
1 & 0 & 0 & 0 & 0 & 0 & 0 & 0\\
0 & 0 & 0 & 0 & 0 & 0 & 0 & 0\\
0 & 0 & 1 & 0 & 0 & 0 & 0 & 0\\
0 & 0 & 0 & 0 & 0 & 0 & 0 & 0\\
0 & 0 & 0 & 0 & 1 & 0 & 0 & 0\\
0 & 0 & 0 & 0 & 0 & 0 & 0 & 0\\
0 & 0 & 0 & 0 & 0 & 0 & 1 & 0\end{array}\right)+\left(\begin{array}{cccccccc}
0 & 0 & 0 & 0 & 0 & 0 & 0 & 0\\
1 & 0 & 0 & 0 & 0 & 0 & 0 & 0\\
0 & 0 & 0 & 0 & 0 & 0 & 0 & 0\\
0 & 0 & 1 & 0 & 0 & 0 & 0 & 0\\
0 & 0 & 0 & 0 & 0 & 0 & 0 & 0\\
0 & 0 & 0 & 0 & 1 & 0 & 0 & 0\\
0 & 0 & 0 & 0 & 0 & 0 & 0 & 0\\
0 & 0 & 0 & 0 & 0 & 0 & 1 & 0\end{array}\right)\left(\begin{array}{cccccccc}
0 & 0 & 0 & 0 & 0 & 0 & 0 & 0\\
\sqrt{7} & 0 & 0 & 0 & 0 & 0 & 0 & 0\\
0 & \sqrt{12} & 0 & 0 & 0 & 0 & 0 & 0\\
0 & 0 & \sqrt{15} & 0 & 0 & 0 & 0 & 0\\
0 & 0 & 0 & 4 & 0 & 0 & 0 & 0\\
0 & 0 & 0 & 0 & \sqrt{15} & 0 & 0 & 0\\
0 & 0 & 0 & 0 & 0 & \sqrt{12} & 0 & 0\\
0 & 0 & 0 & 0 & 0 & 0 & \sqrt{7} & 0\end{array}\right)\\
=\left(\begin{array}{cccccccc}
0 & 0 & 0 & 0 & 0 & 0 & 0 & 0\\
0 & 0 & 0 & 0 & 0 & 0 & 0 & 0\\
\sqrt{12} & 0 & 0 & 0 & 0 & 0 & 0 & 0\\
0 & \sqrt{12} & 0 & 0 & 0 & 0 & 0 & 0\\
0 & 0 & 4 & 0 & 0 & 0 & 0 & 0\\
0 & 0 & 0 & 4 & 0 & 0 & 0 & 0\\
0 & 0 & 0 & 0 & \sqrt{12} & 0 & 0 & 0\\
0 & 0 & 0 & 0 & 0 & \sqrt{12} & 0 & 0\end{array}\right)\\
=2\left(\begin{array}{cccc}
0 & 0 & 0 & 0\\
\sqrt{3} & 0 & 0 & 0\\
0 & 2 & 0 & 0\\
0 & 0 & \sqrt{3} & 0\end{array}\right)\otimes\left(\begin{array}{cc}
1 & 0\\
0 & 1\end{array}\right)=2b^{(3)}\otimes1\boldsymbol{\!\!\!1}.\end{array}\]
\\
Since numbers of the form $p=2^{n}-1$ have the binary representation
1, 11, 111, $\dots$, we say that the above paraorders $p'$ and $p$
are in a carry-bit neighborhood to one another. While its physical
and information-theoretical meaning remain unclear, the operator identity
(\ref{eq:inter-b}) neatly carries over to nil\-potent operators
$f^{(p')}$ which are obtained by {}``extracting the square root''
of $f^{(p)}\otimes1\mathbf{\!\!\!1}$,%
\footnote{note that $f^{(p)}$ untensorized is not a proper exponentiation of
an operator%
} in a recursive process beginning with the fermi operator $f^{(1)}$.
To allow $f^{(p')}$ squared to act as a normalized- anticommutator
analog of Eq.$\,$(\ref{eq:inter-b}), \begin{equation}
(f^{(p')})^{2}=f^{(p)}\otimes1\mathbf{\!\!\!1},\label{eq:inter-f}\end{equation}
 the structure of $f^{(p')}$ has to be amalgamated with $1\mathbf{\!\!\!1}^{\otimes n}\otimes f^{(1)}$,
as we shall see. In matrix form, the structural parts are blockwise
composed of elements of the Clifford algebra $Cl(2,1)$ with basis\begin{equation}
\textrm{\Large\{}c_{1}=({\scriptstyle {1\atop 0}{0\atop -1}}),\: c_{2}=({\scriptstyle {0\atop 1}{1\atop 0}}),\: c_{3}=({\scriptstyle {0\atop -1}{1\atop 0}})\textrm{\Large\}.}\end{equation}
 The simplest representation of the initial operator in the recursive
process consists of a linear combination of one basis element per
signature, usually\begin{equation}
f^{(1)}=\frac{1}{2}(c_{2}-c_{3}).\end{equation}
 Equivalently, and closer to physics, one may start by Clifford algebra
$Cl(3)$ which has the set of Pauli matrices as basis, where $f^{(1)}$
is represented by combining one real basis element of grade 1 -- vector
$\sigma_{1}$ oder $\sigma_{3}$ -- with the only real basis element
of grade 2 -- bivector $\sigma_{31}$, the preferred choice being
\[
f^{(1)}=\frac{1}{2}(\sigma_{1}-\sigma_{31})=\frac{1}{2}(\sigma_{1}^{+}+\sigma_{31}^{+}).\]
 The simplification achieved is that the conjugations $^{^{\mathnormal{+}}}$
und $^{T}$ coincide.

\noindent Solving Eq. (\ref{eq:inter-f}) for $f^{(p')}$ is made
easy by requiring the main block dia\-gonal of the $f^{(p')}$ matrix
to coincide with $1\boldsymbol{\!\!\!1}^{\otimes n}\otimes f^{(1)}$,
and the triangular matrix\emph{ }below it (LTM) to mutually exclusively
consist of blocks \ $G_{\mu\nu}c_{3}$, $E_{\mu\nu}(f^{(1)})^{+}$
or $J_{\mu\nu}c_{2}$ $(\mu>\nu;\: G_{\mu\nu},E_{\mu\nu},J_{\mu\nu}\in\mathbb{Z})$.
As the action below the main diagonal shows,\vspace{-5mm}
 \[
\begin{array}{ccc}
\begin{array}{ccc}
\left(\begin{array}{rcccc}
0 & \,0 & \,\cdots\\
1 & \,0\\
0 & \, x_{1} & \!0 & 0\\
x_{1}y & \,0 & \!1 & 0\\
\vdots & \, & \, &  & \ddots\end{array}\right)\!\begin{array}[b]{l}
\\^{2}\\
\\\\\end{array} & = & \left(\begin{array}{rcccc}
0 & \,\,0 & \,\,\cdots\\
0 & \,\,0 & \,\,\\
1 & \,\,0 & \,\,\\
0 & \,\,1 & \,\,\,_{\ddots}\\
\vdots & \,\, & \,\,\end{array}\right)\end{array} &  & \leadsto x_{1}=1,\qquad(y\in\{-1,0,1\})\end{array}\]
 the whole task is thereby effectively reduced to a linear problem.
To distinguish the resulting matrix sequences from one another, we
call them root-$f$ sequence, root-$d$ sequence and root-$h$ sequence
according to their LTM content: \begin{equation}
f^{(p)}=1\!\!\!\boldsymbol{1}^{\otimes n-1}\otimes f^{(1)}+(G_{\mu\nu}^{(p)})\otimes c_{3},\label{eq:F-detached}\end{equation}
\begin{equation}
d^{(p)}=1\mathbf{\!\!\!1}^{\otimes n-1}\otimes f^{(1)}+(E_{\mu\nu}^{(p)})\otimes(f^{(1)})^{+},\label{eq:D-detached}\end{equation}
\begin{equation}
h^{(p)}=1\mathbf{\!\!\!1}^{\otimes n-1}\otimes f^{(1)}+(J_{\mu\nu}^{(p)})\otimes c_{2}.\label{eq:H-detached}\end{equation}
The solutions to the LPs start out behaving as expected: \[
\begin{array}{lrl}
f^{(1)}=\left(\begin{array}{cc}
0 & 0\\
1 & 0\end{array}\right),\end{array}\]
 \[
\sqrt{f^{(1)}\otimes1\mathbf{\!\!\!1}}=f^{(3)}=\left(\begin{array}{cc}
f^{(1)} & \mathbf{0}\\
c_{3} & f^{(1)}\end{array}\right),\]
 \begin{equation}
\sqrt{f^{(3)}\otimes1\mathbf{\!\!\!1}}=\: f^{(7)}\:=\:\left(\begin{array}{cccc}
f^{(1)} & \mathbf{0} & \mathbf{0} & \mathbf{0}\\
c_{3} & f^{(1)} & \mathbf{0} & \mathbf{0}\\
c_{3} & c_{3} & f^{(1)} & \mathbf{0}\\
c_{3} & c_{3} & c_{3} & f^{(1)}\end{array}\right);\label{eq:f-seq}\end{equation}
 \[
\]
thereafter, however, the bulk of the coefficients $G_{\mu\nu}$ $(\mu>\nu)$
begin to deviate from 1, as a snapshot of the LP $(f^{(15)})^{2}=f^{(7)}\otimes1\mathbf{\!\!\!1}$
when the upper left and lower right quadrants are already computed
--\begin{equation}
\begin{array}{rrrrrrrr}
f^{(1)} & 0 & 0 & 0\\
c_{3} & f^{(1)} & 0 & 0\\
c_{3} & c_{3} & f^{(1)} & 0\\
c_{3} & c_{3} & c_{3} & f^{(1)}\\
G_{51}c_{3} & G_{52}c_{3} & G_{53}c_{3} & G_{54}c_{3} & \; f^{(1)} & 0 & 0 & 0\\
G_{61}c_{3} & G_{62}c_{3} & G_{63}c_{3} & G_{64}c_{3} & \; c_{3} & \; f^{(1)} & 0 & 0\\
G_{71}c_{3} & G_{72}c_{3} & G_{73}c_{3} & G_{74}c_{3} & \; c_{3} & \; c_{3} & \; f^{(1)} & 0\\
G_{81}c_{3} & G_{82}c_{3} & G_{83}c_{3} & G_{84}c_{3} & \; c_{3} & \; c_{3} & \; c_{3} & \; f^{(1)}\\
\\\end{array}\label{eq:snapshot-f15}\end{equation}
-- shows:\medskip{}
 \begin{equation}
\begin{array}{ll}
\mathrm{row}\,5/\mathrm{col}\,4\,\mathrm{downto}\,1: & \begin{array}[t]{l}
G_{54}=1\\
G_{53}-G_{54}=0\;\leadsto G_{53}=1\\
G_{52}-G_{53}-G_{54}=1\;\leadsto G_{52}=3\\
G_{51}-G_{52}-G_{53}-G_{54}=0\;\leadsto G_{51}=5\end{array}\\
\mathrm{\mathrm{row}\,6/\mathrm{col}\,4\, downto\,1}: & \begin{array}[t]{l}
G_{64}-G_{54}=0\;\leadsto G_{64}=1\\
G_{63}-G_{64}-G_{53}=-1\;\leadsto G_{63}=1\\
G_{62}-G_{63}-G_{64}-G_{52}=0\;\leadsto G_{62}=5\\
G_{61}-G_{62}-G_{63}-G_{64}-G_{51}=-1\;\leadsto G_{61}=11\end{array}\\
\mathrm{\mathrm{\mathrm{row}\,7/}\mathrm{col}\,4\, downto\,1}: & \begin{array}[t]{l}
G_{74}-G_{54}-G_{64}=1\;\leadsto G_{74}=3\\
G_{73}-G_{74}-G_{53}-G_{63}=0\;\leadsto G_{73}=5\\
G_{72}-G_{73}-G_{74}-G_{52}-G_{62}=1\;\leadsto G_{72}=17\\
G_{71}-G_{72}-G_{73}-G_{74}-G_{51}-G_{61}=0\;\leadsto G_{71}=41\end{array}\\
\mathrm{\mathrm{row}\,8/\mathrm{\mathrm{col}\,4\, downto}\,1}: & \begin{array}[t]{l}
G_{84}-G_{54}-G_{64}-G_{74}=0\;\leadsto G_{84}=5\\
G_{83}-G_{84}-G_{53}-G_{63}-G_{73}=-1\;\leadsto G_{83}=11\\
G_{82}-G_{83}-G_{84}-G_{52}-G_{62}-G_{72}=0\;\leadsto G_{82}=41\\
G_{81}-G_{82}-G_{83}-G_{84}-G_{51}-G_{61}-G_{71}=-1\;\leadsto G_{81}=113\,.\end{array}\end{array}\label{eq:snapshot-f15-detail}\end{equation}
\[
\]

\noindent Thus,

\noindent \[
\begin{array}{c}
\begin{array}{r}
\sqrt{f^{(7)}\otimes1\mathbf{\!\!\!1}}=\: f^{(15)}\:=\:\\
\\\end{array}\end{array}\left(\begin{array}{cccccccc}
\!\! f^{(1)} & \mathbf{0} &  & \mathbf{} & \cdots &  &  & \!\!\mathbf{0}\\
\!\! c_{3} & \!\! f^{(1)}\\
\!\! c_{3} & c_{3} & \!\! f^{(1)}\\
\!\! c_{3} & c_{3} & c_{3} & \!\! f^{(1)}\\
\!\!5c_{3} & 3c_{3} & c_{3} & c_{3} & \!\! f^{(1)} & \ddots &  & \vdots\\
\!\!11c_{3} & 5c_{3} & c_{3} & c_{3} & c_{3} & \!\! f^{(1)}\\
\!\!41c_{3} & 17c_{3} & 5c_{3} & 3c_{3} & c_{3} & c_{3} & \!\! f^{(1)} & \mathbf{\!\!0}\\
\!\!113c_{3} & 41c_{3} & \!\!11c_{3} & 5c_{3} & c_{3} & c_{3} & \!\! c_{3} & \!\! f^{(1)}\end{array}\right),\cdots.\]
Mutatis mutandis for the root-$d$ and root-$h$ sequences: \[
\begin{array}{ccc}
d^{(1)}=f^{(1)},\end{array}\]
 \[
\]
 \[
\;\sqrt{d^{(1)}\otimes1\mathbf{\!\!\!1}}=d^{(3)}=\left(\begin{array}{cc}
f^{(1)} & \mathbf{0}\\
(f^{(1)})^{+} & f^{(1)}\end{array}\right),\]
 \[
\]
\begin{equation}
\sqrt{d^{(3)}\otimes1\mathbf{\!\!\!1}}=d^{(7)}=\left(\begin{array}{cccc}
f^{(1)} & \mathbf{0} & \mathbf{0} & \mathbf{0}\\
(f^{(1)})^{+} & f^{(1)} & \mathbf{0} & \mathbf{0}\\
\mathbf{\mathbf{0}} & (f^{(1)})^{+} & f^{(1)} & \mathbf{0}\\
\mathbf{\mathbf{0}} & \mathbf{\mathbf{0}} & (f^{(1)})^{+} & f^{(1)}\end{array}\right),\label{eq:d-seq}\end{equation}
 \[
\]
 \[
\begin{array}{c}
\begin{array}{r}
\sqrt{d^{(7)}\otimes1\mathbf{\!\!\!1}}=\: d^{(15)}\:=\:\\
\\\end{array}\:\end{array}\left(\begin{array}{cccccccc}
\!\! f^{(1)} & \mathbf{0} &  & \mathbf{} & \cdots &  &  & \!\!\!\mathbf{0}\\
\!\!(f^{(1)})^{+} & \!\! f^{(1)}\\
\mathbf{\!\!\!\!\!0} & \!\!(f^{(1)})^{+} & \!\! f^{(1)}\\
 &  & \!\!(f^{(1)})^{+} & \!\! f^{(1)}\\
\!\!\!\vdots &  & \!\!\!\ddots & \!\!(f^{(1)})^{+} & \!\! f^{(1)} & \ddots &  & \mathbf{\!\!\!\vdots}\\
 &  &  &  & \!\!(f^{(1)})^{+} & \!\!\! f^{(1)}\\
 &  &  &  &  & \!\!\!(f^{(1)})^{+} & \!\! f^{(1)} & \!\!\!\mathbf{0}\\
\!\!\!\mathbf{0} &  & \cdots &  &  & \!\!\!\!\!\mathbf{0} & \!\!\!(f^{(1)})^{+} & \!\!\! f^{(1)}\end{array}\!\!\right),\cdots;\]
 \[
\]
 \[
\]
 \[
\begin{array}{ccc}
h^{(1)}=f^{(1)},\end{array}\]
 \[
\]
 \[
\;\sqrt{h^{(1)}\otimes1\mathbf{\!\!\!1}}=h^{(3)}=\left(\begin{array}{cc}
f^{(1)} & \mathbf{0}\\
c_{2} & f^{(1)}\end{array}\right),\]
 \[
\]
\begin{equation}
\sqrt{h^{(3)}\otimes1\mathbf{\!\!\!1}}=h^{(7)}=\left(\begin{array}{cccc}
f^{(1)} & \mathbf{0} & \mathbf{0} & \mathbf{0}\\
c_{2} & f^{(1)} & \mathbf{0} & \mathbf{0}\\
\!\!\!-c_{2} & c_{2} & f^{(1)} & \mathbf{0}\\
3c_{2} & \!\!\!-c_{2} & c_{2} & f^{(1)}\end{array}\right),\label{eq:h-seq}\end{equation}
 \[
\]
 \[
\begin{array}{c}
\begin{array}{r}
\sqrt{h^{(7)}\otimes1\mathbf{\!\!\!1}}=\: h^{(15)}\:=\:\\
\\\end{array}\:\end{array}\left(\begin{array}{cccccccc}
f^{(1)} & \mathbf{0} &  & \mathbf{} & \cdots &  &  & \!\!\!\mathbf{0}\\
c_{2} & f^{(1)}\\
-c_{2} & c_{2} & f^{(1)}\\
3c_{2} & -c_{2} & c_{2} & f^{(1)}\\
-5c_{2} & c_{2} & -c_{2} & c_{2} & f^{(1)} & \ddots &  & \!\!\!\vdots\\
15c_{2} & -5c_{2} & 3c_{2} & -c_{2} & c_{2} & f^{(1)}\\
-43c_{2} & 15c_{2} & -5c_{2} & c_{2} & -c_{2} & c_{2} & f^{(1)} & \mathbf{0}\\
149c_{2} & -43c_{2} & 15c_{2} & -5c_{2} & 3c_{2} & -c_{2} & c_{2} & f^{(1)}\end{array}\!\!\right),\cdots.\]
 \[
\]
 While constant recurrence of $E_{\mu\nu}=\delta_{\mu,\nu+1}$ is
to be found throughout the root-$d$ sequence, the coefficients $G_{\mu\nu}$
and $J_{\mu\nu}$ are evolving away from 1, $G_{\mu\nu}$ at $f^{(15)}$,
and $J_{\mu\nu}$ at $h^{(7)}$. Even though not immediately apparent,
the way it happens is controlled by Catalan-type bookkeeping identities:
Where $C_{k}$ denotes the $k$th Catalan number, in Eqs. (\ref{eq:snapshot-f15})-(\ref{eq:snapshot-f15-detail})
there hold the identities\[
\begin{array}{c}
G_{51}=G_{62}=G_{73}=G_{84}=C_{3},\\
G_{52}+G_{61}=G_{74}+G_{83}=C_{4},\\
G_{53}+G_{71}=G_{64}+G_{82}=C_{5},\\
\sum_{i=0}^{3}G_{5+i,4-i}=C_{6},\end{array}\]
 and similar identities hold for the $J_{\mu\nu}$ of the lower left
quadrant of $h^{(7)}$:\[
\begin{array}{c}
-1=-C_{1},\\
1-3=-C_{2},\end{array}\]
as well as for those of the lower left quadrant of $h^{(15)}$:\[
\begin{array}{c}
-5=-C_{3},\\
1-15=-C_{4},\\
-1+43=C_{5},\\
1-3-15+149=C_{6}.\end{array}\]
At the same time, the coefficients $G_{\mu\nu}$ (or $J_{\mu\nu}$)
are linked to an important characteristics of Euclidean $D$-space,
namely the kissing number $L_{D}$ for densest packing of (hyper)spheres
of equal size in that space. For $f^{(7)}$, we have

\noindent \[
L_{1}=G_{4,1}+G_{4,2}=1+1=2,\]

\noindent \noindent and for $f^{(15)}$: \bigskip{}
 \begin{equation}
\begin{array}{ccc}
L_{2} & = & G_{5,1}+G_{5,2}-G_{5,3}-G_{5,4}=5+3-1-1=6,\\
L_{3} & = & G_{7,2}-G_{7,3}=17-5=12,\\
L_{4} & = & G_{7,1}-G_{7,2}=41-17=24,\\
L_{5} & = & G_{7,1}-G_{7,2}+G_{8,3}+G_{8,4}=41-17+11+5=40,\\
L_{6} & = & G_{8,1}+G_{8,2}=113-41=72,\\
L_{7} & = & G_{7,1}-C_{7,2}+G_{8,1}-G_{8,3}=41-17+113-11+=126.\end{array}\label{eq:kiss-function}\end{equation}
\bigskip{}

\noindent \noindent Thus, the members of the root-$f$- and the root-$h$
sequences first of all are partitioners of Catalan numbers in unfamilar
environs, and $f^{(p)}$ and $h^{(p)}$, as they climb up their parental
sequences, traverse all of these. That evolutionary behavior, and
the lack thereof in the root-$d$ sequence, calls for an examination
of whether and how it relates to the very different given of parafermi
algebra. Not surprisingly, the closest resemblance to the Green ansatz
is found among the members of root-$d$ sequence. Not only is the
nilpotence property $\left(d^{(p)}\right)^{p+1}\!\!=\left((d^{(p)})^{+}\right)^{p+1}\!=0$
satisfied, the structure as well is analogous: $d_{\alpha,\beta}=E_{\beta}\,\delta_{\alpha,\beta+1}$
with $E_{\beta}=1$. However, the spin-$(p/2)$ representation, as
we have learned, demands $b_{\alpha,\beta}=\sqrt{\beta(p-\beta+1)})\,\delta_{\alpha,\beta+1}$,
a condition that counters our premise that the main diagonal of $d^{(p')}$
$(p'=2^{n+1}-1)$ carries solely $f^{(1)}$ blocks. So the scope of
reference to parafermi algebra is quickly cut down to members of the
root-$f$- and root-$h$-sequences, and it's foremost the former that
we will pick up for exemplarily scrutinizing a possible relationship
with the Green ansatz. 

\noindent The rest of the paper is organized as follows. In Sect\emph{.
\ref{sec:-f-parafermi-algebra} }we outline the tenets of $f$-parafermi
algebra for order $p\in\left\{ 3,7\right\} $. In Sect. \emph{\ref{sec:A-variant-of},}
an alternative ansatz called heterotic $f$-parafermi algebra, again
of order $p\in\left\{ 3,7\right\} $, is discussed. In both sections,
parenthesized paraorder superscripts are used only where needed. The
structure of the members of the root-$f$ sequence is considered from
various perspectives in Sect\emph{. \ref{sec:Structure-of-the}}.
Specifically in Sect. \emph{\ref{sub:The-interordinal-aspect},} the
interordinal aspect between carry-bit neighbors $f^{(p)}$ and $f^{(p')}$
is elucidated and contrasted with the intraordinal aspect, and the
partial sequences of $G_{\mu\nu}$ representatives, $(G^{(3)})=(1)$,
$(G^{(7)})=(1)$, $(G_{\rho}^{(15)})=(3,5,11,17,41,113),\ldots$,
which already showed up more or less as a curiosity in \cite{Merk89},
are examined for their intrinsic Catalan-number related properties
in Sect\emph{. \ref{sub:Trace-structure-vs.}} and kissing-number
related properties in Sect. \emph{\ref{sub:Row-(column)-structure}};
in Sect\emph{. \ref{sub:The-factorization-aspect}}, the members of
$(G_{\rho}^{(p)})$ are further examined under the factorization aspect,
with an aside to the factorization of the Catalan number $C_{(p-3)/4}$.
As shown in Sect. \emph{\ref{sub:The-modulo-eight-aspect}}, the symmetries
underlying the $(\frac{p+1}{2})\times(\frac{p+1}{2})$ LTM $(G_{\mu\nu}^{(p)})$
allow a generalization of the results of Sects. \emph{\ref{sec:-f-parafermi-algebra}}
and \emph{\ref{sec:A-variant-of}} to orders $p\in\left\{ 15,31,\ldots\right\} $;
advantage is thereby taken of the persistence of symmetry properties
in residues left by the coefficients $G_{\mu\nu}$ after division
by eight. In Sect. \emph{\ref{sec:Structural-comparison-with},} the
analysis is extended to differences derivable from the members of
$(B_{\beta}^{(p)})$ and $(G_{\rho}^{(p)})$; the concept of interordinal
differences is developed in Sect. \emph{\ref{sub:Interordinal-differences}}
and applied to the kissing number problem in Sect\emph{. \ref{sub:Kissing-numbers---}}.
Even though cursorily, in Sect. \emph{\ref{sec:Synopsis-of-root-}}
a glance is taken at the root-$h$ sequence and its partial sequences
of representatives of $J_{\mu\nu}$: $(J^{(3)})=(1)$, $(J_{\omega}^{(7)})=(-1,3)$,
$(J_{\omega}^{(15)})=(1,-5,15,-43,149),\dots$; so-called synoptic
differences derivable from these sequence members and those of $(G_{\rho}^{(p)})$
shed light on the periodicity aspect of kissing number representation
(Sect. \emph{\ref{sec:Kissing-number-representation}}). After a brief
outline of an interordinal preon model (Sect. \emph{\ref{sec:An-interordinal-preon}}),
Sect. \emph{\ref{sec:A-proposal-for}} adds a proposal for a planar
geometric model that nicely fits with the root-$f$- and root-$h$
sequences, and an outline of the connection between Catalan structure
and the model's characteristic feature, the cardioidic arclength,
as well as continued fraction representations -- some of these in
extension of one discussed earlier in Sect\emph{. \ref{sub:Trace-structure-vs.},}
some adressing the kissing number problem from a qphyletic perspective
-- close this article. \bigskip{}

\section{\label{sec:-f-parafermi-algebra}{\em f}$\,$-parafermi algebra}

\noindent To begin with, for $p=2^{n}-1\:(n>1)$, with the exception
of the nilpotence property\vspace{1mm}
 \begin{equation}
f^{p+1}=(f^{+})^{p+1}=0,\end{equation}
 no relation from the Green ansatz is satisfied after substituting
$f$ for $b$. This necessitates an adaptation in form of an orthogonal
decomposition

\begin{equation}
f=\sum_{\upsilon}f_{\upsilon}\end{equation}
 such that\begin{equation}
\qquad\qquad\qquad\qquad\qquad\qquad f_{\upsilon}^{+}f_{\upsilon}=\begin{cases}
\begin{array}{ll}
\mathrm{diag}(\{0,1\}), & \upsilon=0,\\
\mathrm{diag}(\{0\}\cup\{\mathnormal{G_{\mu,s_{\upsilon}(\mu)}^{\,{\displaystyle {\scriptstyle 2}}}\}),} & \upsilon=1,\dots,(p-1)/2,\;\;\mu>s_{\upsilon}(\mu).\end{array} & \qquad\qquad\qquad\end{cases}\label{eq:G-squared}\end{equation}

\hspace*{1cm}\\
 How a $2^{n}\!\times2^{n}$ matrix (here with a granularity of
$2^{2(n-1)}$ blocks $A_{\mu,\nu}$, rather than $2^{2n}$ matrix
elements $m_{\alpha,\beta}$) is orthogonally decomposed into (here
$2^{n-1}={\scriptstyle (}p+{\scriptstyle 1)/}{\scriptstyle 2}$) basis
elements, whilst delineated in literature, is discovered each time
anew. Key part of the decomposition procedure is the index permutations
$s_{\upsilon}(\mu)$ $\cong\mathrm{Z_{2}^{\,\mathit{n}-1}}$. In Table
\ref{tab:Orthonal-decomposition}, the permutations $s_{\upsilon}(\mu)$
$\cong$ $\mathrm{Z}_{2}^{\,3}$ are shown, which are known under
various isomorphic maps from other fields of mathematics (octonions,
Fano plane). For the basis-element characterizations

\begin{equation}
\begin{array}{cc}
\left.f_{0}\!:\;\right. & \left.(a_{1,1})_{\kappa,\lambda}+(a_{2,2})_{\kappa,\lambda}+\dots=A_{\kappa,\lambda}(\delta_{1,\kappa}\delta_{\lambda,1}+\delta_{2,\kappa}\delta_{\lambda,2}+\dots)\right.\\
\left.\right. & \left.\Rightarrow\; f_{0}=a_{1,1}+a_{2,2}+\dots\,=\left(\begin{array}{ccc}
{{\atop A_{1,1}}\atop {\scriptstyle {{\scriptstyle 0}\atop \vdots}}} & {{\atop {\scriptstyle 0}}\atop {A_{1,1}\atop }} & {{\atop {\scriptstyle {\textstyle \cdots}}}\atop {\atop \ddots}}\end{array}\right)\right.\end{array}\end{equation}

\hspace*{1cm}\\
 etc. we use the shorthand $0\,$: $11\!+\!22+\dots\;$ etc. %
\begin{table}[h]
\caption{\label{tab:Orthonal-decomposition}}

\lyxline{\normalsize}

$\qquad$$\upsilon\qquad$$\qquad\qquad\qquad\qquad\qquad\qquad\qquad\qquad{\displaystyle {\textstyle \sum_{\mu}}}a_{\mu,s_{\upsilon}(\mu)}\qquad(\{s_{\upsilon}\}\cong\mathrm{Z}_{2}\!\times\mathrm{Z}_{2}\!\times\mathrm{Z}_{2})\qquad\qquad\qquad\qquad\qquad\qquad$\lyxline{\normalsize}

$\qquad0\qquad$$\qquad\qquad\qquad\qquad\qquad\qquad\qquad\qquad11+22+33+44+55+66+77+88\qquad\qquad\qquad\qquad\qquad\qquad$

$\qquad1\qquad$$\qquad\qquad\qquad\qquad\qquad\qquad\qquad\qquad\underline{12+21}+34+43+56+65+78+87\qquad\qquad\qquad\qquad\qquad\qquad$

$\qquad2\qquad$$\qquad\qquad\qquad\qquad\qquad\qquad\qquad\qquad13+24+31+42+57+68+75+86\qquad\qquad\qquad\qquad\qquad\qquad$

$\qquad3\qquad$$\qquad\qquad\qquad\qquad\qquad\qquad\qquad\qquad\underline{14+23+32+41}+58+67+76+85\qquad\qquad\qquad\qquad\qquad\qquad$

$\qquad4\qquad$$\qquad\qquad\qquad\qquad\qquad\qquad\qquad\qquad15+26+37+48+51+62+73+84\qquad\qquad\qquad\qquad\qquad\qquad$

$\qquad5\qquad$$\qquad\qquad\qquad\qquad\qquad\qquad\qquad\qquad16+25+38+47+52+61+74+83\qquad\qquad\qquad\qquad\qquad\qquad$

$\qquad6\qquad$$\qquad\qquad\qquad\qquad\qquad\qquad\qquad\qquad17+28+35+46+53+64+71+82\qquad\qquad\qquad\qquad\qquad\qquad$

$\qquad7\qquad$$\qquad\qquad\qquad\qquad\qquad\qquad\qquad\qquad\underline{18+27+36+45+54+63+72+81}\qquad\qquad\qquad\qquad\qquad\qquad$\smallskip{}
 \lyxline{\normalsize} 
\end{table}

\noindent Under the delineated proviso we get an $f$-parafermi algebra\\
 \begin{equation}
\frac{1}{2}[f_{0}^{+},f_{0}]+\sum_{\upsilon=1}^{(p-1)/2}[f_{\upsilon}^{+},f_{\upsilon}]=\mathrm{diag}\textrm{\Large(}\frac{p}{2},\frac{p}{2}-1,\cdots,-\frac{p}{2}+1,-\frac{p}{2}\textrm{\Large)},\label{eq:f-par1}\end{equation}
 \begin{equation}
\sum_{\upsilon=0}^{(p-1)/2}[[f_{\upsilon}^{+},f_{\upsilon}],f_{\upsilon}]=-2f,\;\:\sum_{\upsilon=0}^{(p-1)/2}[[f_{\upsilon}^{+},f_{\upsilon}],f_{\upsilon}^{+}]=2f^{+}.\label{eq:f-ar2}\end{equation}

\noindent Note that the above equations hold for $p\in\{3,7\}$; how
they can be generalized for $p\in\{15,31,\dots\}$ will be shown in
Section \emph{\ref{sub:The-modulo-eight-aspect}}.

\hspace*{1cm}\\
 For $p=3$, the orthogonal decomposition reads $f=f_{0}+f_{1}$
where, according to the shorthand prescription $0\!:11\!+\!22,$ $1\!:12\!+\!21$,
we have\[
(f_{0})_{1,1}=(f_{0})_{2,2}=(c_{2}-c_{3})/2\;,\;\;\;(f_{1})_{1,2}=\mathbf{0},\;\;(f_{1})_{2,1}=G_{2,1}c_{3}=c_{3}.{}\]
 \hspace*{1cm}\\
 Mutatis mutandis for case $p=7$.\\
\noindent The spin arithmetics differ in one respect: in Green's
algebra (Eq. \ref{eq:b-spin}), spin values emerge as a half times
differences of squares $B_{\beta^{*}}^{2}-B_{\beta^{*}-1}^{2}\;(\beta^{*}\in\{1,\ldots,p+1\},B_{0},B_{p+1}\equiv0)$,\[
{\scriptstyle \left.{\scriptstyle \begin{array}{cccccc}
\frac{1}{2} & {\textstyle (7-0)} & {\textstyle (12-7)} & {\textstyle (15-12)} & {\textstyle (16-15)\cdots} & {\textstyle (0-7)}\\
 & {\textstyle \frac{7}{2}} & {\textstyle \frac{5}{2}} & {\textstyle \frac{3}{2}} & {\textstyle \qquad\frac{1}{2}}\quad{\textstyle \cdots} & -\frac{7}{2}{\textstyle ,}\end{array}}\right.}\]
 \\
 in $f$-parafermi algebra (Eq. (\ref{eq:f-par1})), as sums of
linear terms,\[
{\scriptstyle \left.\begin{array}{ccccc}
\quad\quad{\textstyle \frac{1}{2}\quad\quad} & {\textstyle -\frac{1}{2}\quad\quad} & {\textstyle \frac{1}{2}}\quad\quad & {\textstyle -\frac{1}{2}}\quad\quad & {\textstyle \cdots}\quad-\frac{1}{2}\\
\quad\quad{\textstyle 3\quad\quad} & {\textstyle 3}\quad\quad & {\textstyle 1\quad\quad} & {\textstyle 1\quad\quad} & {\textstyle \cdots}{\textstyle \quad-3}\\
\quad\quad{\textstyle \frac{7}{2}}\quad\quad & {\textstyle \frac{5}{2}}\quad\quad & {\textstyle \frac{3}{2}}\quad\quad & {\textstyle \frac{1}{2}}\quad\quad & {\textstyle \cdots\quad-\frac{7}{2}.}\end{array}\right.}\]

\section{\label{sec:A-variant-of}A variant of {\em f}$\,$-parafermi algebra}

\noindent Even though it seems unlikely it can reveal parafermionic
aspects of the structure of $(G_{\mu\nu}^{(p)})$, a second version
of $f$-parafermi algebra is worth reviewing. One always finds a $g$
(a matrix with free parameters in general), for which\begin{equation}
[[f^{+}f],g]=-2f,\;\:[[f^{+},f],g^{+}]=2f^{+}.\label{eq:f-var1}\end{equation}
 As the system of linear equations embraced by $g$ is underdetermined,
one has to constrain the block structure of $g$ to -- compared with
$f$'s -- slightly relaxed linear combinations $H_{\mu\nu}c_{2}+K_{\mu\nu}c_{3}$$\;(\mu,\nu=1,\dots,{\scriptstyle (}p+{\scriptstyle 1}{\scriptstyle )}{\scriptstyle /{\textstyle {\scriptstyle 2}}})$
to get the solutions unique, or their range narrowed by further constraints,
and $g$ thus constructed. The spin-$p/2$ representation is recovered
by imposing the requirement $g=\sum_{\upsilon}g_{\upsilon},\;(g_{\upsilon})_{\mu,s_{\upsilon}{\scriptscriptstyle (\mu)}}=H_{\mu,s_{\upsilon}({\scriptscriptstyle \mu})}c_{2}\!+\! K_{\mu,s_{\upsilon}({\scriptscriptstyle \mu})}c_{3}\;\;(\{s_{\upsilon}\}\cong\mathrm{Z_{2}^{\,\mathnormal{n}-1}})$
and choosing the ansatz\vspace{1mm}
 \begin{equation}
\begin{array}{c}
\sum_{\upsilon=0}^{(p-1)/2}\textrm{\Large(}\chi[f_{\upsilon}^{+},f_{\upsilon}]+\sigma\left([f_{\upsilon}^{+},g_{\upsilon}]+[g_{\upsilon}^{+},f_{\upsilon}]\right)+\tau\left([f_{\upsilon},g_{\upsilon}]+[g_{\upsilon}^{+},f_{\upsilon}^{+}]\right)+\gamma[g_{\upsilon},g_{\upsilon}^{+}]\textrm{\Large)}\\
=\:\mathrm{diag}\textrm{\Large(}\frac{p}{2},\frac{p}{2}-1,\cdots,-\frac{p}{2}+1,-\frac{p}{2}\textrm{\Large)}\end{array}\label{eq:f-var2}\end{equation}
 so that Eqs. (\ref{eq:f-var1})--(\ref{eq:f-var2}) may be slated
as a heterotic version of $f$-parafermi algebra. Whatever relationship
there might exist between $B_{\beta}$ and $G_{\mu\nu}\;(\mu>\nu)$,
by the additional quantities $H_{\mu\nu}$,$K_{\mu\nu}$ $(\mu,\nu=1,\dots,{\scriptstyle (}p+{\scriptstyle 1)/}{\scriptstyle 2})$
and $\chi,\sigma,\tau,\gamma$ it is rather concealed than revealed.
The steps of computation to be taken shall nevertheless briefly be
expounded for paraorders 3 and 7, whilst postponing the question of
how a generalization for $p\in\{15,31,\dots\}$ might look like. The
LSE for $g^{(3)}$ has a unique solution which reads\[
g^{(3)}=\left(\begin{array}{cccc}
0 & 1 & 0 & \frac{1}{2}\\
0 & 0 & -\frac{3}{2} & 0\\
0 & \frac{3}{2} & 0 & 1\\
-\frac{1}{2} & 0 & 0 & 0\end{array}\right).\]
 \noindent In a similar way as $f^{(3)}$ was treated, $g^{(3)}$
is orthogonally decomposed by following the prescription $0\!:11\!+\!22,$
$1\!:12\!+\!21$, which yields $g=g_{0}+g_{1}$, with block structures\[
\begin{array}{cc}
(g_{0})_{1,1}=(g_{0})_{2,2}=\frac{1}{2}c_{2}+\frac{1}{2}c_{3}, & (g_{1})_{1,2}=c_{3}-\frac{1}{2}c_{2},\;(g_{1})_{2,1}=c_{3}+\frac{1}{2}c_{2}.\end{array}\]
 By the LSE of the spin-3/2 representation (Eq. (\ref{eq:f-var2}))
we then obtain a parametrized set of solutions for the normalizing
factors,\[
\begin{array}{l}
\chi^{(3)}=(4r_{2}+2r_{1}+2)/3,\\
\sigma^{(3)}=(-10r_{2}-2r_{1}+1)/2,\\
\tau^{(3)}=r_{2},\\
\gamma^{(3)}=r_{1}.\end{array}\;(r_{i}\,\mathrm{free}\,\mathrm{parameters})\]

\noindent Solving the LSE for $g^{(7)}$ raises a matrix with no less
than four free parameters! Of which we may free us -- not arbitrarily,
but by imposing on $g^{(7)}$ the very same symmetries that govern
$g^{(3)}$. Three types of these can be read from the above representation
of $g^{(3)}$ ($A^{T}$ transposed matrix, $A^{S}$ matrix reflected
in secondary diagonal): 1) $({A\atop -B^{T}}{B\atop A})$; 2) $({A\atop -B^{T}}{B\atop A^{S}})$;
3) $({A_{0}\atop C}{B\atop A_{0}^{S}})$, where the subscript in $A_{0}$
is indicative of a {}``zero block'' in the lower left part of the
secondary diagonal: $A_{0}=({U\atop \mathbf{0}}{V\atop W})$. In fact,
each of the symmetries 1) to 3) effects the complete eli\-mination
of degrees of freedom from $g^{(7)}$, leading to the LSE solutions
\[
\begin{array}{l}
\;\;1)\;\,\end{array}\; g^{(7)}=\left(\!\!\begin{array}{cccccccc}
0 & \frac{5}{8} & 0 & \frac{2}{5} & 0 & \frac{9}{40} & 0 & \frac{1}{10}\\
\frac{3}{8} & 0 & \frac{-1}{4} & 0 & \frac{-1}{8} & 0 & \frac{1}{4} & 0\\
0 & \frac{1}{4} & 0 & \frac{5}{8} & 0 & \frac{9}{20} & 0 & \frac{9}{40}\\
\frac{-2}{5} & 0 & \frac{3}{8} & 0 & \frac{-1}{5} & 0 & \frac{-1}{8} & 0\\
0 & \frac{1}{8} & 0 & \frac{1}{5} & 0 & \frac{5}{8} & 0 & \frac{2}{5}\\
\frac{-9}{40} & 0 & \frac{-9}{20} & 0 & \frac{3}{8} & 0 & \frac{-1}{4} & 0\\
0 & \frac{-1}{4} & 0 & \frac{1}{8} & 0 & \frac{1}{4} & 0 & \frac{5}{8}\\
\frac{-1}{10} & 0 & \frac{-9}{40} & 0 & \frac{-2}{5} & 0 & \frac{3}{8} & 0\end{array}\right),\;\;\left\{ \begin{array}{ll}
\chi^{(7)}=1, & \;\sigma^{(7)}=\frac{-(r_{1}+2)}{8},\\
\\\tau^{(7)}=\frac{r_{1}+2}{8}, & \;\gamma^{(7)}=r_{1}\end{array}\right\} ;\]

\[
\begin{array}{l}
2)\;\end{array}\; g^{(7)}=\left(\!\!\begin{array}{cccccccc}
0 & \frac{1}{24} & 0 & \frac{-1}{60} & 0 & \frac{-1}{40} & 0 & \frac{1}{60}\\
\frac{23}{24} & 0 & 0 & 0 & \frac{-5}{24} & 0 & \frac{-1}{6} & 0\\
0 & 0 & 0 & \frac{3}{8} & 0 & \frac{1}{5} & 0 & \frac{-1}{40}\\
\frac{1}{60} & 0 & \frac{5}{8} & 0 & \frac{-7}{60} & 0 & \frac{-5}{24} & 0\\
0 & \frac{5}{24} & 0 & \frac{7}{60} & 0 & \frac{3}{8} & 0 & \frac{-1}{60}\\
\frac{1}{40} & 0 & \frac{-1}{5} & 0 & \frac{5}{8} & 0 & 0 & 0\\
0 & \frac{1}{6} & 0 & \frac{5}{24} & 0 & 0 & 0 & \frac{1}{24}\\
\frac{-1}{60} & 0 & \frac{1}{40} & 0 & \frac{1}{60} & 0 & \frac{23}{24} & 0\end{array}\right),\;\;\left\{ \begin{array}{ll}
\chi^{(7)}=1, & \;\sigma^{(7)}=-1/4,\\
\\\tau^{(7)}=1/4, & \;\gamma^{(7)}=0\end{array}\right\} ;\]

\[
\begin{array}{l}
\quad\quad\;\;\;3)\;\end{array}\; g^{(7)}=\left(\!\!\begin{array}{cccccccc}
0 & \!\frac{1}{24} & \!0 & \!\frac{-1}{60} & \!0 & \!\frac{-1}{40} & \!0 & \!\frac{1}{60}\\
\frac{187}{200} & \!0 & \!\frac{-1}{100} & \!0 & \!\frac{-41}{200} & 0 & \!\frac{-3}{20} & \!0\\
0 & \!0 & 0 & \!\frac{3}{8} & \!0 & \!\frac{1}{5} & \!0 & \!\frac{-1}{40}\\
0 & \!0 & \!\frac{123}{200} & \!0 & \!\frac{-3}{25} & 0 & \!\frac{-41}{200} & \!0\\
0 & \!\frac{5}{24} & \!0 & \!\frac{7}{60} & \!0 & \!\frac{3}{8} & \!0 & \!\frac{-1}{60}\\
\frac{3}{200} & \!0 & \!\frac{-21}{200} & \!0 & \!\frac{123}{200} & \!0 & \!\frac{-1}{100} & \!0\\
0 & \!\frac{1}{6} & \!0 & \!\frac{5}{24} & \!0 & \!0 & \!0 & \!\frac{1}{24}\\
\frac{-1}{50} & \!0 & \!\frac{3}{200} & \!0 & \!0 & \!0 & \!\frac{187}{200} & \!0\end{array}\right),\;\;\left\{ \begin{array}{ll}
\chi^{(7)}=\frac{14183539}{14137018}, & \sigma^{(7)}=\frac{-1737725}{7068509},\\
\\\tau^{(7)}=\frac{1738225}{7068509}, & \gamma^{(7)}=\frac{147500}{7068509}\end{array}\right\} .\]
 \[
\]
 Conspicuously, variant 2) seems to bring out a {}``standard set''
of normalizing factors $\{\chi^{(7)}\!=\!1,$ $\sigma^{(7)}=\frac{-1}{4},\tau^{(7)}=\frac{1}{4},\gamma^{(7)}=0\}$,
which those belonging to 1) can be made to conform to by the choice
$r_{1}=0$ and which those belonging to 3) differ from by no more
than $\approx$ 2\%. Viewed in this light, $\{\chi^{(3)}\!\!\!=\!1,\sigma^{(3)}\!\!\!=\!\frac{-3}{4},\tau^{(3)}\!\!\!=\!\frac{1}{4},\gamma^{(3)}\!\!\!=\!0\}$
can be considered the standard normalizing factor set for $p=3$.
It cannot be excluded that other types of symmetries expand the range
of viable solutions; lack of symmetry however -- by simply setting
all four parameters in the general matrix of $g^{(7)}$ equal to zero
-- only results in $\{\}$ for the set of normalizing factors.

\noindent After this aside, we turn to interordinality as a way of
exploring the putative parafermionic nature of the coefficients $G_{\mu\nu}$
by relating them directly to the Green coefficients $B_{\beta}$.\newpage{}

\section{\label{sec:Structure-of-the}Structure of the members of root-{\em
f} sequence}

\subsection{\label{sub:The-interordinal-aspect}Interordinal aspect vs. intraordinal
aspect}

\noindent We first encountered the structural interordinal aspect
in links between the paraorder $p$ and its upper carry-bit neighbor
$p'$ (Eqs. (\ref{eq:inter-b}) and (\ref{eq:inter-f}) which pair
$b^{(p)}$, $b^{(p')}$ and $f^{(p)}$, $f^{(p')}$ respectively)
and with its lower carry-bit neighbor (the $(p+1)\!\times\!(p+1)$
array of $f^{(p)}$ and $h^{(p)}$ vs. the $(\frac{p-1}{2}+1)\!\times\!(\frac{p-1}{2}+1)$
structure of $(G_{\mu\nu}^{(p)})$ and $(J_{\mu\nu}^{(p)})$. We briefly
touched a further link existing between $p$ and its next but one
lower carry-bit neighbor (the Catalan ``accounting identities'' controlling
the $(\frac{p-3}{4}+1)\!\times\!(\frac{p-3}{4}+1)$ structure of the
lower left quadrant of $(G_{\mu\nu}^{(p)})$ and $(J_{\mu\nu}^{(p)})$),
and we shall soon find it necessary to enlarge the picture by one
or several more (higher or lower) carry-bit neighbors, so it seems
suitable to adopt a shorthand for them. In analogy to the denotation
of Mersenne numbers $M_{n}=2^{n}-1$, we occasionally find it convenient
to write $p_{n}$ for $p\equiv2^{n}+1$, $p_{n+1}$ for $p'=2p+1$
and so on, also $q_{n}$ for $q=(p-3)/4$ , $q_{n-1}$ for $(p-7)/8$
or $q_{n+1}$ for $q'=(p-1)/2$.\\
Let us begin with the link between $b^{(p)}$ and $b^{(p')}$ where
the structural interordinal aspect enters in a basic way: every coefficient
that falls within order $p$ is echoed by every second coefficient
that falls within order $p'$ via the relation $B_{2\beta}^{(p')}=2B_{\beta}^{(p)}$;
for instance\[
\begin{array}{cccccccccc}
p'=15: & \qquad\sqrt{15} & \sqrt{28} & \sqrt{39} & \sqrt{48} & \sqrt{55} & \sqrt{60} & \sqrt{63} & \sqrt{64} & \cdots,\\
p=\;7: &  & \sqrt{7} &  & \sqrt{12} &  & \sqrt{15} &  & \sqrt{16} & \cdots.\end{array}\]

\noindent A variant of this doubling effect then likely is to be
found when ascending from $f^{(p)}$ to $f^{(p')}$. Before elaborating,
let us address the exponential nature of the objects we deal with.
$f^{(31)}=\sqrt{\sqrt{\sqrt{\sqrt{f^{(1)}\otimes1\mathbf{\!\!\!1}}\otimes1\mathbf{\!\!\!1}}\otimes1\mathbf{\!\!\!1}}\otimes1\!\!\!\boldsymbol{1}}$
is the first of root-$f$ sequence members that is too wide to fit
standard paper size -- as we need navigation of some form for them,
however, we introduce a minimum of new notation, speaking of upper/lower
left, or upper/lower right, parts to retain some depictability. Specific
quadrants are determined by one-place navigation $\upharpoonleft\!\downharpoonright\leftrightharpoons\textrm{arg}$,
subquadrants by $\mathrm{\upharpoonleft\!\downharpoonright}\!\leftrightharpoons\!(\upharpoonleft\!\downharpoonright\leftrightharpoons\textrm{arg}),$
and so on. The first observation worth a mention is that all subquadrants
(and quadrants, as well as $f^{(p)}$ itself) show invariance under
reflection in the secondary diagonal -- sometimes called secondary
symmetry \cite{ALee76}: \begin{equation}
\mathrm{\upharpoonleft\!\downharpoonright}\!\leftrightharpoons(\upharpoonleft\!\downharpoonright\leftrightharpoons f^{(p)})=(\mathrm{\upharpoonleft\!\downharpoonright}\!\leftrightharpoons(\upharpoonleft\!\downharpoonright\leftrightharpoons f^{(p)}))^{S},\label{eq:secondary-sym}\end{equation}

\noindent where the sequence of symbols $\mathrm{\upharpoonleft\!\downharpoonright}\!\leftrightharpoons$
is the same on both sides of the equation. One further is that identical
subquadrant content appears at different places, namely at UL(LL$f$)
and LR(LL$f$). Also, at LL(UL$f$) and LL(LR$f$) and, flanked by
these, at UR(LL($f$)):%
\begin{figure}[H]

\caption{\label{fig:subq-content-1}Identical subquadrants}
\[
f^{(7)}=\left(\begin{array}{cccc}
f^{(1)} & \boldsymbol{0} & \boldsymbol{0} & \boldsymbol{0}\\
\fbox{{\ensuremath{\mathit{c}_{3}}}} & f^{(1)} & \boldsymbol{0} & \mathbf{0}\\
\fbox{{\fbox{\ensuremath{c_{3}}}}} & \fbox{{\ensuremath{\mathit{c}_{3}}}} & f^{(1)} & \boldsymbol{0}\\
c_{3} & \fbox{\fbox{{\ensuremath{\mathit{c}_{3}}}}} & \fbox{{\ensuremath{\mathit{c}_{3}}}} & f^{(1)}\end{array}\right),\qquad f^{(15)}=\left(\begin{array}{cccc}
\begin{array}{cc}
f^{(1)} & \boldsymbol{0}\\
c_{3} & f^{(1)}\end{array} & \begin{array}{cc}
\cdots\\
\ddots\end{array} & \begin{array}{cc}
\\\\\end{array} & \begin{array}{cc}
\boldsymbol{\cdots} & \boldsymbol{0}\\
 & \boldsymbol{\vdots}\end{array}\\
\fbox{\mbox{{\ensuremath{\begin{array}{cc}
 \mathit{c}_{3}  &  \mathit{c}_{3}\\
\mathit{c}_{3}  &  \mathit{c}_{3}\end{array}}}}} & \begin{array}{cc}
f^{(1)} & \boldsymbol{0}\\
c_{3} & f^{(1)}\end{array} &  & \begin{array}{cc}
\\\\\end{array}\\
\fbox{\fbox{\mbox{{\ensuremath{\begin{array}{cc}
 \mathit{\mathrm{5}c}_{3}  &  3\mathit{c}_{3}\\
11\mathit{c}_{3}  &  5\mathit{c}_{3}\end{array}}}}}} & \fbox{{\mbox{{\ensuremath{\begin{array}{cc}
 \mathit{c}_{3}  &  \mathit{c}_{3}\\
\mathit{c}_{3}  &  \mathit{c}_{3}\end{array}}}}}} & \begin{array}{cc}
f^{(1)} & \boldsymbol{0}\\
c_{3} & f^{(1)}\end{array} & \begin{array}{cc}
\\\ddots & \vdots\end{array}\\
\begin{array}{cc}
41c_{3} & 17c_{3}\\
113c_{3} & 41c_{3}\end{array} & \fbox{\fbox{\mbox{{\ensuremath{\begin{array}{cc}
 5\mathit{c}_{3}  &  3\mathit{c}_{3}\\
11\mathit{c}_{3}  &  5\mathit{c}_{3}\end{array}}}}}} & \fbox{{\mbox{{\ensuremath{\begin{array}{cc}
 \mathit{c}_{3}  &  \mathit{c}_{3}\\
\mathit{c}_{3}  &  \mathit{c}_{3}\end{array}}}}}} & \begin{array}{cc}
f^{(1)} & \boldsymbol{0}\\
c_{3} & f^{(1)}\end{array}\end{array}\right),\:\ldots\]

\end{figure}
\noindent  Now the content of the single framed blocks in $\mathrm{\mathit{f}^{(\mathit{p}')}}$
is identical to the quadrant $\mathrm{LL\mathit{f}^{(\mathit{p})}}$.
Thus, if the double framed blocks are suggestive of the notion of
subquadrantal in\emph{tra}ordinality, the single framed ones (shown
bracketed in what follows) may be attributed to what we here call
subquadrantal in{\em ter}ordina\-li\-ty (not to be confused with
the term used in statistics). 

\noindent  This property of two types of ordinality governing the
root-$f$ structure is further refined on the subsubquandral level,
emerging first with $p=7$, $\: p'=15$:%
\begin{figure}[H]
\caption{\label{fig:interord-subq-1}Interordinally related subsubquadrants}
\[
\mathrm{\hbox{\small{interordinal\, identity:}}\qquad\qquad\qquad\qquad\qquad\qquad\qquad\qquad\qquad}\qquad\qquad\qquad\]
\[
\]
\begin{equation}
\mathrm{UR(UL(LL{\mathit{f}}^{(\mathit{p'})}))}=\textrm{LL}(\textrm{LL}\mathit{f}^{(\mathit{p})})+2\,\mathrm{UR(LL{\mathit{f}}^{(\mathit{p})})}\label{eq:inter}\end{equation}
\[
\]
\[
\]
\[
\mathrm{LL}\, f^{(15)}=\left(\begin{array}{cc}
\begin{array}{cc}
5c_{3} & \fbox{\mbox{{3\ensuremath{\mathit{c}_{3}}}}}\\
11c_{3} & 5c_{3}\end{array} & \left[\begin{array}{cc}
c_{3} & \fbox{\mbox{{\ensuremath{\mathit{c}_{3}}}}}\\
\fbox{\mbox{{\ensuremath{\mathit{c}_{3}}}}} & c_{3}\end{array}\right]\\
\begin{array}{cc}
41c_{3} & 17c_{3}\\
113c_{3} & 41c_{3}\end{array} & \begin{array}{cc}
5c_{3} & 3c_{3}\\
11c_{3} & 5c_{3}\end{array}\end{array}\right),\]
\[
\]
\[
\]
\[
\mathrm{LL}\, f^{(31)}=\left(\begin{array}{cc}
\begin{array}{cc}
\begin{array}{cc}
429c_{3} & 155c_{3}\\
1275c_{3} & 429c_{3}\end{array} & \fbox{\mbox{{\ensuremath{\begin{array}{cc}
 43\mathit{c}_{3}  &  19\mathit{c}_{3}\\
115\mathit{c}_{3}  &  43\mathit{c}_{3}\end{array}}}}}\\
\begin{array}{cc}
4819c_{3} & 1595c_{3}\\
15067c_{3} & 4819c_{3}\end{array} & \begin{array}{cc}
429c_{3} & 155c_{3}\\
1275c_{3} & 429c_{3}\end{array}\end{array} & \left[\begin{array}{cc}
\begin{array}{cc}
5c_{3} & 3c_{3}\\
11c_{3} & 5c_{3}\end{array} & \fbox{\mbox{{\ensuremath{\begin{array}{cc}
 \mathit{c}_{3}  &  \;\mathit{c}_{3}\\
\mathit{c}_{3}  &  \;\mathit{c}_{3}\end{array}}}}}\\
\fbox{\mbox{{\ensuremath{\begin{array}{cc}
 41\mathit{c}_{3}  &  17\mathit{c}_{3}\\
113\mathit{c}_{3}  &  41\mathit{c}_{3}\end{array}}}}} & \begin{array}{cc}
5c_{3} & 3c_{3}\\
11c_{3} & 5c_{3}\end{array}\end{array}\right]\\
\begin{array}{cc}
\begin{array}{cc}
58781c_{3} & 18627c_{3}\\
189371c_{3} & 58781c_{3}\end{array} & \begin{array}{cc}
4905c_{3} & 1633c_{3}\\
15297c_{3} & 4905c_{3}\end{array}\\
\begin{array}{cc}
737953c_{3} & 227089c_{3}\\
2430289c_{3} & 737953c_{3}\end{array} & \begin{array}{cc}
58781c_{3} & 18627c_{3}\\
189371c_{3} & 58781c_{3}\end{array}\end{array} & \begin{array}{cc}
\begin{array}{cc}
429c_{3} & 155c_{3}\\
1275c_{3} & 429c_{3}\end{array} & \begin{array}{cc}
43c_{3} & 19c_{3}\\
115c_{3} & 43c_{3}\end{array}\\
\begin{array}{cc}
4819c_{3} & 1595c_{3}\\
15067c_{3} & 4819c_{3}\end{array} & \begin{array}{cc}
429c_{3} & 155c_{3}\\
1275c_{3} & 429c_{3}\end{array}\end{array}\end{array}\right),\:\ldots;\]

\end{figure}

\[
\]

\begin{figure}[H]
\caption{\label{fig:intraord-subq-1}Intraordinally related subsubquadrants}

\end{figure}
\[
\mathrm{\hbox{\small{intraordinal\, identity:}}\qquad\qquad\qquad\qquad\qquad\qquad\qquad\qquad\qquad}\qquad\qquad\qquad\]
\[
\]
\begin{equation}
\mathrm{\mathrm{UR(LL(LL{\mathit{f}}^{(\mathit{p'})}))}=LL(UL(LL\mathit{f}^{(\mathit{p'})}))+2\, UR(UL(LL{\mathit{f}}^{(\mathit{p'})}))}\label{eq:intra}\end{equation}
\[
\]
\[
\]
\[
\mathrm{LL}\, f^{(15)}=\left(\begin{array}{cc}
\begin{array}{cc}
5c_{3} & \fbox{\fbox{\mbox{{3\ensuremath{\mathit{c}_{3}}}}}}\\
\fbox{\fbox{\mbox{{11\ensuremath{\mathit{c}_{3}}}}}} & 5c_{3}\end{array} & \left[\begin{array}{cc}
{{\textstyle c_{3}}\atop } & {{\textstyle \; c_{3}}\atop }\\
{\atop {\textstyle c_{3}}} & {\atop {\textstyle \; c_{3}}}\end{array}\right]\\
\begin{array}{cc}
41c_{3} & \fbox{\fbox{\mbox{{17\ensuremath{\mathit{c}_{3}}}}}}\\
113c_{3} & 41c_{3}\end{array} & \begin{array}{cc}
{{\textstyle 5c_{3}}\atop } & {{\textstyle 3c_{3}}\atop }\\
{\atop {\textstyle 11c_{3}}} & {\atop {\textstyle 5c_{3}}}\end{array}\end{array}\right),\]
\newpage{}\hfill{}{\small Fig.} \ref{fig:intraord-subq-1} {\small continued:}\hfill{}

\[
\]
\[
\mathrm{LL}\, f^{(31)}=\left(\begin{array}{cc}
\begin{array}{cc}
\begin{array}{cc}
429c_{3} & 155c_{3}\\
1275c_{3} & 429c_{3}\end{array} & \fbox{\fbox{\mbox{{\ensuremath{\begin{array}{cc}
 43\mathit{c}_{3}  &  19\mathit{c}_{3}\\
115\mathit{c}_{3}  &  43\mathit{c}_{3}\end{array}}}}}}\\
\fbox{\fbox{\mbox{{\ensuremath{\begin{array}{cc}
 4819\mathit{c}_{3}  &  1595\mathit{c}_{3}\\
15067\mathit{c}_{3}  &  4819\mathit{c}_{3}\end{array}}}}}} & \begin{array}{cc}
429c_{3} & 155c_{3}\\
1275c_{3} & 429c_{3}\end{array}\end{array} & \left[\begin{array}{cc}
\begin{array}{cc}
5c_{3} & 3c_{3}\\
11c_{3} & 5c_{3}\\
\\\end{array} & \begin{array}{cccccc}
 &  &  & c_{3} &  & \: c_{3}\\
 &  &  & c_{3} &  & \: c_{3}\\
\\\end{array}\\
\begin{array}{cc}
41c_{3} & 17c_{3}\\
113c_{3} & 41c_{3}\end{array} & \begin{array}{cccccc}
 &  &  & 5c_{3} &  & 3c_{3}\\
 &  &  & 11c_{3} &  & 5c_{3}\end{array}\end{array}\right]\\
\begin{array}{cc}
\begin{array}{cc}
58781c_{3} & 18627c_{3}\\
189371c_{3} & 58781c_{3}\end{array} & \fbox{\fbox{\mbox{{\ensuremath{\begin{array}{cc}
 4905\mathit{c}_{3}  &  1633\mathit{c}_{3}\\
15297\mathit{c}_{3}  &  4905\mathit{c}_{3}\end{array}}}}}}\\
\begin{array}{cc}
737953c_{3} & 227089c_{3}\\
2430289c_{3} & 737953c_{3}\end{array} & \begin{array}{cc}
58781c_{3} & 18627c_{3}\\
189371c_{3} & 58781c_{3}\end{array}\end{array} & \begin{array}{cc}
\begin{array}{cc}
429c_{3} & 155c_{3}\\
1275c_{3} & 429c_{3}\\
\\\end{array} & \begin{array}{cc}
43c_{3} & 19c_{3}\\
115c_{3} & 43c_{3}\\
\\\end{array}\\
\begin{array}{cc}
4819c_{3} & 1595c_{3}\\
15067c_{3} & 4819c_{3}\end{array} & \begin{array}{cc}
429c_{3} & 155c_{3}\\
1275c_{3} & 429c_{3}\end{array}\end{array}\end{array}\right),\:\ldots;\]
 \medskip{}

\noindent \negthinspace{}\normalsize identities that, due to the
structural symmetries noted above, find an equivalent in\[
\]
\begin{equation}
\mathrm{\mathrm{UR(LR(LL{\mathit{f}}^{(\mathit{p'})}))}=LL(LL\mathit{f}^{(\mathit{p})})}+2\,\mathrm{UR(LL{\mathit{f}}^{(\mathit{p})})}\label{eq:sine-like}\end{equation}
 and\begin{equation}
\mathrm{UR(LL(LL{\mathit{f}}^{(\mathit{p'})}))}=\mathrm{LL(LR(LL\mathit{f}^{(\mathit{\mathit{p}'})}))}+2\,\mathrm{UR(LR(LL}{f}^{(\mathit{p'})}))\label{eq:cosine-like}\end{equation}
\[
\]
respectively. As will be explained in Sect. \emph{\ref{sub:Continued-fraction-representation-}},
the interordinal identities (\ref{eq:inter},\ref{eq:sine-like})
can be classified as sine-like, and the intraordinal identities (\ref{eq:intra},\ref{eq:cosine-like})
as cosine-like.

\noindent \noindent The logical consequence of the structural interordinal
aspect is that it restricts the domain from which to select specific
coefficients as representatives falling within order $p$. We define
the representatives $G_{\rho}^{(p)}$ as to be taken from those $G_{\mu\nu}^{(p)}$
that spring from the nonbracketed part of $\mathrm{LL}\mathit{f}^{(\mathit{p})}$,
denoted $\neg\mathrm{UR(LL\mathit{f}^{(\mathit{p})}})$. To find out
how many such representatives $G_{\rho}^{(p)}$ exist, we have to
address Catalan structure next.

\subsection{\label{sub:Trace-structure-vs.}Secondary trace structure vs. stub
structure}

\noindent As mentioned in the introduction, the bookkeeping on the
coefficients of $\mathrm{LL}(G_{\mu\nu}^{(p)})$ is done by way of
$\textrm{str}()$, the symbol signifying traces over the secondary
and adjacent diagonals, in that quadrant. Taking these as gross traces,\[
\textrm{\textrm{gstr}}\, G_{q+1+\xi-\zeta,\xi+\zeta}^{(p)}\quad(\xi,\zeta+1\in\{1,\dots,q+1\}),\]
 yields sums of $C_{q-1+\xi}$ in general. For $\mathrm{LL}(G_{\mu\nu}^{(31)})$
for instance:\[
\begin{array}{lclccc}
\qquad G_{9,1} & = & C_{7}, & \Sigma_{\zeta=0}^{1}G_{10-\zeta,1+\zeta} & = & C_{8},\\
\Sigma_{\zeta=0}^{2}G_{11-\zeta,1+\zeta} & = & C_{9}+C_{7}, & \Sigma_{\zeta=0}^{3}G_{12-\zeta,1+\zeta} & = & C_{10},\\
\Sigma_{\zeta=0}^{4}G_{13-\zeta,1+\zeta} & = & C_{11}+C_{9}+C_{7}, & \Sigma_{\zeta=0}^{5}G_{14-\zeta,1+\zeta} & = & C_{12},\\
\Sigma_{\zeta=0}^{6}G_{15-\zeta,1+\zeta} & = & C_{13}+C_{11}+C_{7}, & \Sigma_{\zeta=0}^{7}G_{16-\zeta,1+\zeta} & = & C_{14}.\\
\\\end{array}\]
An alternative way to keep records is with net secondary traces, $\textrm{nstr}()$:\bigskip{}

\noindent \[
\textrm{\textrm{nstr}}\, G_{q_{n}+1+\xi-\zeta,\xi+\zeta}^{(p)}=C_{q_{n}-1+\xi}\]
where\[
\xi,\zeta+1\in\{1,\dots,q_{n}+1\};\; G_{q_{n}+1+\xi-\zeta,\xi+\zeta}^{(p)}\neq G_{\mu\nu}^{*(p)}.\]
\bigskip{}

\begin{figure}[H]
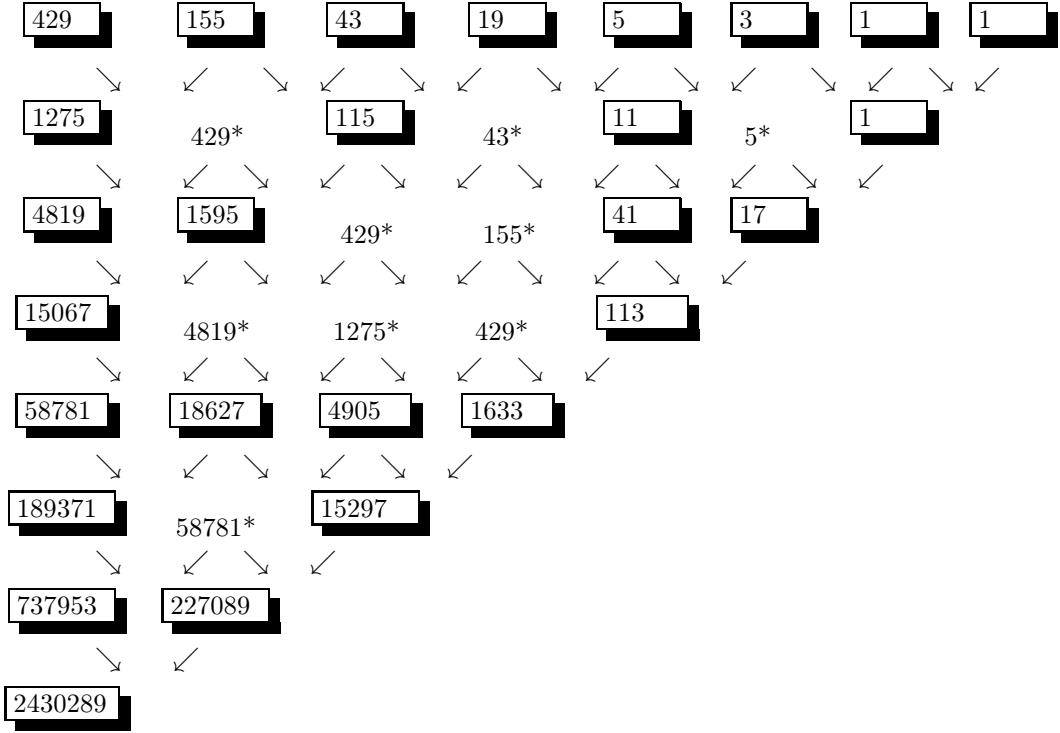

\caption{\label{fig:near-traces}Secondary trace structure of $\mathrm{LL}(G_{\mu\nu}^{(31)})$}
\bigskip{}

\qquad{}\qquad{}\qquad{}\begin{tabular}{cccccccc}
\noalign{\vskip5mm}
\shadowbox{\begin{minipage}[t]{8mm}%
429%
\end{minipage}} & %
\shadowbox{\begin{minipage}[t]{8mm}%
155%
\end{minipage}} & %
\shadowbox{\begin{minipage}[t]{8mm}%
43%
\end{minipage}} & %
\shadowbox{\begin{minipage}[t]{8mm}%
19%
\end{minipage}} & %
\shadowbox{\begin{minipage}[t]{8mm}%
5%
\end{minipage}} & %
\shadowbox{\begin{minipage}[t]{8mm}%
3%
\end{minipage}} & %
\shadowbox{\begin{minipage}[t]{8mm}%
1%
\end{minipage}} & %
\shadowbox{\begin{minipage}[t]{8mm}%
1%
\end{minipage}}\tabularnewline
\noalign{\vskip5mm}
\shadowbox{\begin{minipage}[t]{8mm}%
1275%
\end{minipage}} & %
\begin{minipage}[t]{8mm}%
429{*}$\;$%
\end{minipage} & %
\shadowbox{\begin{minipage}[t]{8mm}%
115%
\end{minipage}} & %
\begin{minipage}[t]{8mm}%
43{*}$\;$%
\end{minipage} & %
\shadowbox{\begin{minipage}[t]{8mm}%
11%
\end{minipage}} & %
\begin{minipage}[t]{8mm}%
5{*}$\;$%
\end{minipage} & %
\shadowbox{\begin{minipage}[t]{8mm}%
1%
\end{minipage}} & \tabularnewline
\noalign{\vskip5mm}
\shadowbox{\begin{minipage}[t]{8mm}%
4819%
\end{minipage}} & %
\shadowbox{\begin{minipage}[t]{8mm}%
1595%
\end{minipage}} & %
\begin{minipage}[t]{8mm}%
429{*}$\;$%
\end{minipage} & %
\begin{minipage}[t]{8mm}%
155{*}$\;$%
\end{minipage} & %
\shadowbox{\begin{minipage}[t]{8mm}%
41%
\end{minipage}} & %
\shadowbox{\begin{minipage}[t]{8mm}%
17%
\end{minipage}} &  & \tabularnewline
\noalign{\vskip5mm}
\shadowbox{\begin{minipage}[t]{10mm}%
15067%
\end{minipage}} & %
\begin{minipage}[t]{10mm}%
4819{*}$\;$%
\end{minipage} & %
\begin{minipage}[t]{10mm}%
1275{*}$\;$%
\end{minipage} & %
\begin{minipage}[t]{10mm}%
429{*}$\;$%
\end{minipage} & %
\shadowbox{\begin{minipage}[t]{10mm}%
113%
\end{minipage}} &  &  & \tabularnewline
\noalign{\vskip5mm}
\shadowbox{\begin{minipage}[t]{10mm}%
58781%
\end{minipage}} & %
\shadowbox{\begin{minipage}[t]{10mm}%
18627%
\end{minipage}} & %
\shadowbox{\begin{minipage}[t]{10mm}%
4905%
\end{minipage}} & %
\shadowbox{\begin{minipage}[t]{10mm}%
1633%
\end{minipage}} &  &  &  & \tabularnewline
\noalign{\vskip5mm}
\shadowbox{\begin{minipage}[t]{12mm}%
189371%
\end{minipage}} & %
\begin{minipage}[t]{12mm}%
58781{*}$\;$%
\end{minipage} & %
\shadowbox{\begin{minipage}[t]{12mm}%
15297%
\end{minipage}} &  &  &  &  & \tabularnewline
\noalign{\vskip5mm}
\shadowbox{\begin{minipage}[t]{12mm}%
737953%
\end{minipage}} & %
\shadowbox{\begin{minipage}[t]{12mm}%
227089%
\end{minipage}} &  &  &  &  &  & \tabularnewline
\noalign{\vskip5mm}
\shadowbox{\begin{minipage}[t]{13mm}%
2430289%
\end{minipage}} &  &  &  &  &  &  & \tabularnewline
\end{tabular}

\vspace{-9.8cm}
$\qquad\qquad\qquad\qquad$\begin{tabular}{cccccccc}
\noalign{\vskip5mm}
\begin{minipage}[t]{12mm}%
$\quad\;\searrow$%
\end{minipage} & %
\begin{minipage}[t]{14mm}%
$\!\swarrow\qquad\searrow$%
\end{minipage} & %
\begin{minipage}[t]{14mm}%
$\!\swarrow\qquad\searrow$%
\end{minipage} & %
\begin{minipage}[t]{14mm}%
$\!\swarrow\qquad\searrow$%
\end{minipage} & %
\begin{minipage}[t]{14mm}%
$\!\swarrow\qquad\searrow$%
\end{minipage} & %
\begin{minipage}[t]{14mm}%
$\!\swarrow\qquad\searrow$%
\end{minipage} & %
\begin{minipage}[t]{14mm}%
$\!\swarrow\quad\;\searrow$%
\end{minipage} & %
\begin{minipage}[t]{13mm}%
$\negthinspace\!\!\!\negthinspace\!\negthinspace\!\swarrow$%
\end{minipage}\tabularnewline[1mm]
\noalign{\vskip5mm}
\begin{minipage}[t]{12mm}%
$\quad\;\searrow$%
\end{minipage} & %
\begin{minipage}[t]{14mm}%
$\!\swarrow\quad\;\searrow$%
\end{minipage} & %
\begin{minipage}[t]{14mm}%
$\!\swarrow\quad\;\searrow$%
\end{minipage} & %
\begin{minipage}[t]{14mm}%
$\!\swarrow\quad\;\searrow$%
\end{minipage} & %
\begin{minipage}[t]{14mm}%
$\!\swarrow\quad\;\searrow$%
\end{minipage} & %
\begin{minipage}[t]{14mm}%
$\!\swarrow\quad\;\searrow$%
\end{minipage} & %
\begin{minipage}[t]{13mm}%
$\!\!\!\!\swarrow$%
\end{minipage} & \tabularnewline[1mm]
\noalign{\vskip5mm}
\begin{minipage}[t]{12mm}%
$\quad\;\searrow$%
\end{minipage} & %
\begin{minipage}[t]{14mm}%
$\!\swarrow\quad\;\searrow$%
\end{minipage} & %
\begin{minipage}[t]{14mm}%
$\!\swarrow\quad\;\searrow$%
\end{minipage} & %
\begin{minipage}[t]{14mm}%
$\!\swarrow\quad\;\searrow$%
\end{minipage} & %
\begin{minipage}[t]{14mm}%
$\!\swarrow\quad\;\searrow$%
\end{minipage} & %
\begin{minipage}[t]{13mm}%
$\!\!\!\!\swarrow$%
\end{minipage} &  & \tabularnewline[1mm]
\noalign{\vskip5mm}
\begin{minipage}[t]{12mm}%
$\quad\;\searrow$%
\end{minipage} & %
\begin{minipage}[t]{14mm}%
$\!\swarrow\quad\;\searrow$%
\end{minipage} & %
\begin{minipage}[t]{14mm}%
$\!\swarrow\quad\;\searrow$%
\end{minipage} & %
\begin{minipage}[t]{14mm}%
$\!\swarrow\quad\;\searrow$%
\end{minipage} & %
\begin{minipage}[t]{13mm}%
$\!\!\!\!\swarrow$%
\end{minipage} &  &  & \tabularnewline[1mm]
\noalign{\vskip5mm}
\begin{minipage}[t]{12mm}%
$\quad\;\searrow$%
\end{minipage} & %
\begin{minipage}[t]{14mm}%
$\!\swarrow\quad\;\searrow$%
\end{minipage} & %
\begin{minipage}[t]{14mm}%
$\!\swarrow\quad\;\searrow$%
\end{minipage} & %
\begin{minipage}[t]{13mm}%
$\!\!\!\!\swarrow$%
\end{minipage} &  &  &  & \tabularnewline[1mm]
\noalign{\vskip5mm}
\begin{minipage}[t]{12mm}%
$\quad\;\searrow$%
\end{minipage} & %
\begin{minipage}[t]{14mm}%
$\!\swarrow\quad\;\searrow$%
\end{minipage} & %
\begin{minipage}[t]{13mm}%
$\!\!\!\!\swarrow$%
\end{minipage} &  &  &  &  & \tabularnewline[1mm]
\noalign{\vskip5mm}
\begin{minipage}[t]{12mm}%
$\quad\;\searrow$%
\end{minipage} & %
\begin{minipage}[t]{13mm}%
$\!\!\!\!\swarrow$%
\end{minipage} &  &  &  &  &  & \tabularnewline[1mm]
\noalign{\vskip5mm}
\begin{minipage}[t]{13mm}%
$\qquad$%
\end{minipage} &  &  &  &  &  &  & \tabularnewline[1mm]
\end{tabular}
\end{figure}

\noindent \noindent This one requires a preprocessing where each
main-type diagonal, the main diagonal itself and its adjacents, are
traversed from the upper leftmost entry to the lower rightmost and
checked for duplicates which are marked when identified. The $\textrm{nstr}()$
then simply skips marked entries (see Fig. \ref{fig:near-traces}
where duplicates are marked by an asterisk).

\noindent %

\noindent \noindent Now, any summand lying in a (net secondary) trace
$r$ positions away from that of $C_{q}$ obeys an upper bound $4^{r}C_{q}$
because ${\displaystyle \lim_{n\rightarrow\infty}}C_{n+1}/C_{n}=4$.
Thus $G_{\textrm{max}}^{(p)}\equiv G_{2q+2,1}^{(p)}$, though the
summand largest in its trace, does stay well below this bound: $113\approx2.82^{3}\cdot5$;
and $2430289\approx3.43^{7}\cdot429$. From $G_{\textrm{max}}^{(p)}=(\Phi^{(p)})^{q}\, C_{q}$
we may define a span parameter $\Phi^{(p)}$ with continued fraction
representation (CFR)\begin{equation}
\Phi^{(p)}=\left(\frac{G_{\textrm{max}}^{(p)}}{C_{q}}\right)^{1/q}=\phi_{0}^{(p)}+\frac{1}{\phi_{1}^{(p)}+{\displaystyle \frac{1}{\;\phi_{2}^{(p)}{\displaystyle +\frac{1}{\phi_{3}^{(p)}+{\displaystyle \ddots}}}}}}\;\equiv\;[\phi_{0}^{(p)};\phi_{1}^{(p)},\phi_{2}^{(p)},\phi_{3}^{(p)},\ldots]\,.\label{eq:contfrac}\end{equation}
\[
\]
\noindent The $\mathrm{LL}(G_{\mu\nu}^{(p)})$ coefficients $<G_{\textrm{max}}^{(p)}$
then become coefficients $\phi_{\alpha+j(2q+2)-\left\lfloor (j+1)/2\right\rfloor }^{(p)}$
of the span parameter for some start value $\alpha$. Those of $\mathrm{LL}(G_{\mu\nu}^{(15)})$
take the form $\phi_{5+8j-\left\lfloor (j+1)/2\right\rfloor }^{(15)}$
for $j=0,1,\ldots,5$, while the special form $G_{\textrm{max}}^{(15)}\equiv G_{8,1}^{(15)}=\phi_{50}^{(15)}-\Sigma_{i=1}^{4}p_{i}$
is assumed for $j=6$:\emph{}%
\footnote{Some conspicuous near matches are springing up incidentally: Five
of kissing numbers with deviations $\Delta_{L}<G_{\textrm{max}}^{(15)}$,
\emph{viz}. $\phi_{50}^{(15)}=139=L_{7}+13$, $\phi_{173}^{(15)}=230=L_{8}-10$,
$\phi_{403}^{(15)}=431=L_{11}-7$, $\phi_{308}^{(15)}=1432=L_{14}+10$,
$\phi_{128}^{(15)}=10558=L_{19}-110$; and two of Catalan numbers
with $\Delta_{C}<C_{3}$, \emph{viz.} $\phi_{308}^{(15)}=1432=C_{8}+2$,
$\phi_{403}^{(15)}=431=C_{7}+2$. We shall come back to the relationship
between kissing numbers and Catalan numbers before long.%
}\[
\begin{array}{cccccccccccccccccccc}
\Phi^{(15)}= & \left[2;\right. & \ldots & 1, & \ldots & 3, & \ldots, & \;5, & \ldots & \,11, & \ldots & \;17, & \ldots & \;41, & \ldots & 139, & \left.\ldots\right]\\
{\scriptscriptstyle } & {\scriptscriptstyle } & {\scriptscriptstyle } & {\scriptscriptstyle } & {\scriptscriptstyle } & {\scriptscriptstyle } & {\scriptscriptstyle } & {\scriptscriptstyle } & {\scriptscriptstyle } & {\scriptscriptstyle } & {\scriptscriptstyle } & {\scriptscriptstyle } & {\scriptscriptstyle } & {\scriptscriptstyle } & {\scriptscriptstyle } & {\scriptscriptstyle } & {\scriptscriptstyle } & {\scriptscriptstyle } & {\scriptscriptstyle } & {\scriptscriptstyle }\\
{\atop } & {0\atop } & {\atop } & {5\atop } & {\atop } & {12\atop } & {\atop } & {20\atop } & {\atop } & {27\atop } & {\atop } & {35\atop } & {\atop } & {42\atop } & {\atop } & {50\atop } & {\atop } & {\atop } & {\atop } & {\atop }\end{array}\]
\[
\]

\noindent A corollary to its secondary symmetry is what we call stub
structure of $\mathrm{LL\,}(G_{\mu\nu}^{(p)})$, defined by the coefficient-wise
homogeneity of its main diagonal and the main diagonals of its subquadrants,
subsubquadrants etc. to either side. It may also be recognized that
the subquadrants UL(LL$(G_{\mu\nu}^{(p)}))$ and LL(LL$(G_{\mu\nu}^{(p)}))$
suffice as sources for $G_{\rho}^{(p)}$, and what is more, intraordinal
relation (\ref{eq:intra}) guarantees them to be independent. The
number of different $G_{\rho}^{(p)}$, denoted $\mathcal{T}_{p}$
here, is readily computed by inspecting these subquadrants using their
secondary symmetry and stub structure. As square matrices of dimension
$\frac{q+1}{2}\times\frac{q+1}{2}$ with secondary symmetry, they
have at most $\frac{(q+1)(q+3)}{8}$ distinct elements, the entries
to a secondary, skewed triangular matrix each:%
\begin{figure}[H]
\caption{\label{fig:stubs in ULLL and LLLL}Stub structure}

\bigskip{}
\bigskip{}

\hfill{}$\mathrm{LL}\,(G_{\mu\nu}^{(15)}):\qquad\left(\begin{array}{c}
\\\\\\\\\end{array}\right.$\begin{tabular}{cc}
$\left\{ \!\begin{array}{c}
\\\\\end{array}\right.$\begin{tabular}{cc}
\shadowbox{\begin{minipage}[t]{8mm}%
5%
\end{minipage}} & %
\shadowbox{\begin{minipage}[t]{8mm}%
3%
\end{minipage}}\tabularnewline
\shadowbox{\begin{minipage}[t]{8mm}%
11%
\end{minipage}} & \tabularnewline
\end{tabular}$\left.\begin{array}{c}
\\\\\end{array}\!\right\} $ & $\left[\begin{array}{c}
\\\\\end{array}\right.$\begin{tabular}{cc}
\begin{minipage}[t]{8mm}%
1%
\end{minipage} & %
\begin{minipage}[t]{8mm}%
1%
\end{minipage}\tabularnewline
\begin{minipage}[t]{8mm}%
1%
\end{minipage} & \tabularnewline
\end{tabular}\hspace{-3mm}$\left]\begin{array}{c}
\\\\\end{array}\right.$\tabularnewline
$\left\{ \!\begin{array}{c}
\\\\\end{array}\right.$\begin{tabular}{cc}
\shadowbox{\begin{minipage}[t]{8mm}%
41%
\end{minipage}} & %
\shadowbox{\begin{minipage}[t]{8mm}%
17%
\end{minipage}}\tabularnewline
\shadowbox{\begin{minipage}[t]{8mm}%
113%
\end{minipage}} & \tabularnewline
\end{tabular}$\left.\begin{array}{c}
\\\\\end{array}\!\right\} $ & \begin{tabular}{cc}
 & \tabularnewline
 & \tabularnewline
\end{tabular}\tabularnewline
\end{tabular}$\negthickspace\left)\begin{array}{c}
\\\\\\\\\end{array}\right.$,\hfill{}\bigskip{}
\bigskip{}
$\mathrm{LL}\,(G_{\mu\nu}^{(31)}):\qquad$$\left(\begin{array}{c}
\\\\\\\\\\\\\\\\\end{array}\right.$\begin{tabular}{cc}
$\left\{ \begin{array}{c}
\\\\\\\\\end{array}\right.$\begin{tabular}{cccc}
\shadowbox{\begin{minipage}[t]{12mm}%
429%
\end{minipage}} & %
\shadowbox{\begin{minipage}[t]{12mm}%
155%
\end{minipage}} & %
\shadowbox{\begin{minipage}[t]{12mm}%
43%
\end{minipage}} & %
\shadowbox{\begin{minipage}[t]{12mm}%
19%
\end{minipage}}\tabularnewline
\shadowbox{\begin{minipage}[t]{12mm}%
1275%
\end{minipage}} & %
\begin{minipage}[t]{8mm}%
429{*}%
\end{minipage} & %
\shadowbox{\begin{minipage}[t]{12mm}%
115%
\end{minipage}} & \tabularnewline
\shadowbox{\begin{minipage}[t]{12mm}%
4819%
\end{minipage}} & %
\shadowbox{\begin{minipage}[t]{12mm}%
1595%
\end{minipage}} &  & \tabularnewline
\shadowbox{\begin{minipage}[t]{12mm}%
15067%
\end{minipage}} &  &  & \tabularnewline
\end{tabular}$\left.\begin{array}{c}
\\\\\\\\\end{array}\right\} $ & $\left[\begin{array}{c}
\\\\\\\\\end{array}\right.$\begin{tabular}{cccc}
\begin{minipage}[t]{8mm}%
5%
\end{minipage} & %
\begin{minipage}[t]{8mm}%
3%
\end{minipage} & %
\begin{minipage}[t]{8mm}%
1%
\end{minipage} & %
\begin{minipage}[t]{8mm}%
1%
\end{minipage}\tabularnewline
\begin{minipage}[t]{8mm}%
11%
\end{minipage} & %
\begin{minipage}[t]{8mm}%
\end{minipage} & %
\begin{minipage}[t]{8mm}%
1%
\end{minipage} & \tabularnewline
\begin{minipage}[t]{8mm}%
41%
\end{minipage} & %
\begin{minipage}[t]{8mm}%
17%
\end{minipage} &  & \tabularnewline
\begin{minipage}[t]{8mm}%
113%
\end{minipage} &  &  & \tabularnewline
\end{tabular}\hspace{-4mm}$\left.\begin{array}{c}
\\\\\\\\\end{array}\right]$\tabularnewline
$\left\{ \begin{array}{c}
\\\\\\\\\end{array}\right.$\begin{tabular}{cccc}
\shadowbox{\begin{minipage}[t]{12mm}%
58781%
\end{minipage}} & %
\shadowbox{\begin{minipage}[t]{12mm}%
18627%
\end{minipage}} & %
\shadowbox{\begin{minipage}[t]{12mm}%
4905%
\end{minipage}} & %
\shadowbox{\begin{minipage}[t]{12mm}%
1633%
\end{minipage}}\tabularnewline
\shadowbox{\begin{minipage}[t]{12mm}%
189371%
\end{minipage}} & %
\begin{minipage}[t]{12mm}%
58781{*}%
\end{minipage} & %
\shadowbox{\begin{minipage}[t]{12mm}%
15297%
\end{minipage}} & \tabularnewline
\shadowbox{\begin{minipage}[t]{12mm}%
737953%
\end{minipage}} & %
\shadowbox{\begin{minipage}[t]{12mm}%
227089%
\end{minipage}} &  & \tabularnewline
\shadowbox{\begin{minipage}[t]{12mm}%
2430289%
\end{minipage}} &  &  & \tabularnewline
\end{tabular}$\left.\begin{array}{c}
\\\\\\\\\end{array}\right\} $ & \begin{tabular}{cccc}
 &  &  & \tabularnewline
 &  &  & \tabularnewline
 &  &  & \tabularnewline
 &  &  & \tabularnewline
\end{tabular}\tabularnewline
\end{tabular}$\left)\begin{array}{c}
\\\\\\\\\\\\\\\\\end{array},\:\cdots\right.$
\end{figure}

\noindent Redundant copies of entries on the subquadrantal main-diagonal,
and subsubquadrantal main-diagonal etc. stubs to either side (see
Fig. \ref{fig:stubs in ULLL and LLLL}) have to be subtracted yet.
For $p=15,$ $\frac{(q+1)(q+3)}{8}=3$, there are no subtractions,
hence $\mathcal{T}_{15}=6$. For $p=31,$ $\frac{(q+1)(q+3)}{8}=10$,
there's one subtraction for either subquadrant, hence $\mathcal{T}_{31}=18$.
For $p=63,$ $\frac{(q+1)(q+3)}{8}=36$, tedious but straightforward
calculations show that nine upper-subquadrantal entries and nine lower-subquadrantal
entries are redundant, hence $\mathcal{T}_{63}=54$. Thus, while the
subtractions seem nontrivial, the result boils down to a simple formula
for the number of distinct representatives, \begin{equation}
\mathcal{T}_{p}=2\cdot3^{\log_{2}(p+1)-3}\qquad(p=15,31,63,\ldots).\label{eq:num-rep}\end{equation}

\noindent See Table \ref{tab:2-partioning1-3-5} for a summary:\pagebreak{}

\noindent %
\begin{table}[h]
\caption{Number of representatives, $\mathcal{T}_{p}\equiv$ \#$G_{\rho}^{(p)}$
\smallskip{}
 \label{tab:2-partioning1-3-5}}

\lyxline{\normalsize}

$\qquad\quad\neg\textrm{UR}(\textrm{LL}\mathrm{\mathit{f}^{(15)})}$
$\mathrm{\mathrm{\qquad\qquad\neg UR(LL}\mathrm{\mathit{f}^{(31)})}}$$\mathrm{\qquad\qquad\qquad\neg UR(LL}\mathrm{\mathit{f}^{(63)}})$$\mathrm{\mathrm{\qquad\qquad\qquad\neg UR(LL}\mathrm{\mathit{f}^{(127)}}})$$\mathrm{\mathrm{\mathrm{\qquad\qquad\qquad\neg UR(LL}\mathrm{\mathit{f}^{(255)}}})}$

\lyxline{\normalsize}

\smallskip{}
$\mathcal{T}_{p}\qquad\qquad6$$\mathrm{\mathrm{\qquad\qquad\qquad\qquad\qquad18}}$$\mathrm{\mathrm{\qquad\qquad\qquad\qquad\qquad\qquad54}}$$\mathrm{\mathrm{\qquad\qquad\qquad\qquad\qquad\qquad162\mathrm{\mathrm{\qquad\qquad\qquad\qquad\qquad\quad486}}}}$\\

\noindent $\qquad(3-0)+(3-0)$$\qquad\qquad$$\quad(10-1)+(10-1)$$\mathrm{\mathrm{\qquad\qquad}}$$\quad(36-9)+(36-9)\mathrm{\mathrm{\qquad\qquad}}$$\mathrm{\mathrm{\quad(\textrm{136-55})+(\textrm{136-55})\qquad\qquad}}$$\quad(\textrm{528-285})+(\textrm{528-285})$

\lyxline{\normalsize} \medskip{}

\end{table}

\noindent \noindent  The stub structure-based ansatz takes account
of subtractions (row 2) and reads\begin{equation}
2\left(\dfrac{(q+1)(q+3)}{8}-s\right).\label{eq:subtractions}\end{equation}
The subtractions $s_{n},\; n=\log_{2}(p_{n}+1)$, build from two types
of atoms, $m_{\mu}\equiv p_{\mu}+p_{\mu-1}$, and $o_{\nu}\equiv(p_{\nu}+p_{\nu-2})/2$.
For $p_{5}=31$, $s_{5}=m_{1}=1$. For $p_{6}=63$, $s_{6}=o_{4}=9$.
For $p_{7}=127$, the first more complex case, a mix shows up, $s_{7}=m_{5}+o_{4}=55$.
This mix begins to branch in increasingly complex ways, $s_{8}=(m_{7}+m_{5}+m_{1})+(o_{6}+o_{4})=285$,
$s_{9}=\ldots$, but should ultimaltely lead to a proof of the following
conjecture: Let $\mathcal{L}_{s}$ be the set of numbers $\left\lfloor \log_{2}(C_{q*}C_{2q*+1})\right\rfloor >s$,
where $q*\in\{1,3,7,17,\ldots\}$. Then the least element $l_{\textrm{min}}\in\mathcal{L}_{s}$
satisfies $l_{\textrm{min}}-s\equiv u\,\textrm{(mod 16)}$ where $u\in\{1,5,9,13\}$.
For $\frac{(q+1)(q+3)}{8}=10$, $s=1$, we find $\left\lfloor \log_{2}(C_{1}C_{3})\right\rfloor -s\equiv1\,\textrm{\textrm{(mod 16)}};$
also for $\frac{(q+1)(q+3)}{8}=36$, $s=9$, $\left\lfloor \log_{2}(C_{3}C_{7})\right\rfloor -s\equiv1\,\textrm{(mod 16)};$
the next instances are $\frac{(q+1)(q+3)}{8}=136$, $s=55$ with $\left\lfloor \log_{2}(C_{15}C_{31})\right\rfloor -s\equiv5\,\textrm{\textrm{(mod 16)}}$,
$\frac{(q+1)(q+3)}{8}=528$, $s=285$ with $\left\lfloor \log_{2}(C_{63}C_{127})\right\rfloor -s\equiv9$
(mod 16), and so on.

\subsection{\label{sub:Row-(column)-structure}Row (column) structure by way
of kissing-number stenoscopy}

\noindent Much in the same way as the (secondary-diagonal) trace
structure and (diagonal) stub structure of $\mathrm{LL}(G_{\mu\nu}^{(p)})$
are governed by Catalan numbers, rows or columns, by way of additive
partitions of their elements which in the simplest case are of length
$q\pm1$, harbor kissing numbers that are sandwiched between Catalan
numbers in a stenoscopic way as shown in Table \ref{tab:Stenoscopy-of-kissing}.
\\
For $p=15$, the partitions in question have already been presented
in Eqs.(\ref{eq:kiss-function}); of these, those which are of length
4 are\bigskip{}

\[
\begin{array}{lll}
L_{2}=6: & \qquad5+3-1-1 & \qquad{G_{5,1}\atop \centerdot}\:{G_{5,2}\atop \centerdot}\:{G_{5,3}\atop \centerdot}\:{G_{5,4}\atop \centerdot}\qquad\qquad\left(\mathrm{UL(LL}(G_{\mu\nu}^{(p)})),\mathrm{UR(LL}(G_{\mu\nu}^{(p)}))\right)\\
\\L_{5}=40: & \qquad41-17+11+5 & \qquad{G_{7,1}\atop \centerdot}\:{G_{7,2}\atop \centerdot}\:{\centerdot\atop G_{8,3}}\:{\centerdot\atop G_{8,4}}\qquad\qquad\left(\mathrm{LL(LL}(G_{\mu\nu}^{(p)})),\mathrm{LR(LL}(G_{\mu\nu}^{(p)}))\right)\\
\\L_{7}=126: & \qquad41-17+113-11 & \qquad{G_{7,1}\atop G_{8,1}}\:{G_{7,2}\atop \centerdot}\:{\centerdot\atop G_{8,3}}\;\enskip{\centerdot\atop \centerdot}\qquad\qquad\quad\left(\mathrm{LL(LL}(G_{\mu\nu}^{(p)})),\mathrm{LR(LL}(G_{\mu\nu}^{(p)}))\right),\end{array}\]

\bigskip{}
and those of length 2 \bigskip{}
\[
\begin{array}{lll}
L_{3}=12: & \qquad\enskip17-5 & \enskip\quad\qquad\qquad\quad{\centerdot\atop \centerdot}\,\enskip{G_{7,2}\atop \centerdot}\;{G_{7,3}\atop \centerdot}\:\quad{\centerdot\atop \centerdot}\qquad\qquad\quad\left(\mathrm{LL(LL}(G_{\mu\nu}^{(p)})),\mathrm{LR(LL}(G_{\mu\nu}^{(p)}))\right)\\
\\L_{4}=24: & \qquad\enskip41-17 & \enskip\quad\qquad\qquad\;{G_{7,1}\atop \centerdot}\:{G_{7,2}\atop \centerdot}\,\enskip{\centerdot\atop \centerdot}\qquad{\centerdot\atop \centerdot}\qquad\qquad\quad\left(\mathrm{LL(LL}(G_{\mu\nu}^{(p)})),\mathrm{LR(LL}(G_{\mu\nu}^{(p)}))\right)\\
\\L_{6}=72: & \qquad\enskip113-41 & \enskip\quad\qquad\qquad\;{\centerdot\atop G_{8,1}}\:{\centerdot\atop G_{8,2}}\,\enskip{\centerdot\atop \centerdot}\qquad{\centerdot\atop \centerdot}\qquad\qquad\quad\left(\mathrm{LL(LL}(G_{\mu\nu}^{(p)})),\mathrm{LR(LL}(G_{\mu\nu}^{(p)}))\right).\end{array}\]
\pagebreak{}

\noindent For $p'=31$, there is an interordinal corridor in which
the remaining kissing numbers $\Sigma_{i=1}^{q+1}G_{2q+2,i}^{(p)}>L_{x}<C_{q'}(\!\,=C_{7})$
reside; they are given by additive partitions of length $\rho\in\{q\pm1,q'\pm1\}$:%
\footnote{entries $G_{\mu\nu}^{(p')}$ that owe their existence to the subsubquadrantal
identity (\ref{eq:inter}) are set in parentheses, while those springing
from $\textrm{UR}(\textrm{LL}\mathrm{(\mathit{G}_{\mu\nu}^{(\mathit{p'})}))}$
are set in brackets%
}\bigskip{}
\[
\begin{array}{lll}
L_{8}=240: & \;\textrm{429-155-(43-19)-[5+3+1+1]} & {{{\scriptstyle \quad G_{9,1}\enskip G_{9,2}\enskip G_{9,3}\enskip G_{9,4}\enskip G_{9,5}\enskip G_{9,6}\enskip G_{9,7}\enskip G_{9,8}}\atop {\scriptstyle \quad\centerdot\qquad\centerdot\qquad\centerdot\qquad\centerdot\qquad\,\centerdot\qquad\centerdot\qquad\centerdot\qquad\centerdot}}\atop {{\scriptstyle {\scriptstyle \quad\centerdot\qquad\centerdot\qquad\centerdot\qquad\centerdot\qquad\,\centerdot\qquad\centerdot\qquad\centerdot\qquad\centerdot}}\atop {\scriptstyle {\scriptstyle \quad\centerdot\qquad\centerdot\qquad\centerdot\qquad\centerdot\qquad\,\centerdot\qquad\centerdot\qquad\centerdot\qquad\centerdot}}}}\:\quad\left(\mathrm{UL(LL}(G_{\mu\nu}^{(p')})),\mathrm{UR(LL}(G_{\mu\nu}^{(p')}))\right)\\
\\L_{9}=272: & \;\textrm{429-115-[41+17]+[11+5]} & {{{\scriptstyle {\scriptstyle \centerdot\qquad\centerdot\qquad\centerdot\qquad\centerdot\qquad\centerdot\qquad\,\centerdot\qquad\centerdot\qquad\centerdot}}\atop {\scriptstyle {\scriptstyle \centerdot\quad G_{10,2}\quad G_{10,3}\,\enskip\centerdot\qquad\centerdot\qquad\,\centerdot\qquad\centerdot\qquad\centerdot}}}\atop {{\scriptstyle \centerdot\qquad\centerdot\qquad\centerdot\qquad\centerdot\;\quad G_{11,5}\enskip G_{11,6}\quad\centerdot\qquad\centerdot}\atop {\scriptstyle {\scriptstyle {\scriptstyle \quad\!\centerdot\qquad\centerdot\qquad\centerdot\qquad\centerdot\qquad\centerdot\qquad\centerdot\quad G_{12,7}\enskip G_{12,8}}}}}}\quad\left(\mathrm{UL(LL}(G_{\mu\nu}^{(p')})),\mathrm{UR(LL}(G_{\mu\nu}^{(p')}))\right)\\
\\L_{10}=336: & \;\textrm{429-155+(43+19)} & \qquad\qquad\qquad\qquad\,{{{\scriptstyle {\scriptstyle \: G_{13,5}\enskip G_{13,6}\enskip G_{13,7}\enskip G_{13,8}}}\atop {\scriptstyle \quad\centerdot\qquad\centerdot\qquad\centerdot\qquad\centerdot}}\atop {{\scriptstyle \quad\centerdot\qquad\centerdot\qquad\centerdot\qquad\centerdot}\atop {\scriptstyle \quad\centerdot\qquad\centerdot\qquad\centerdot\qquad\centerdot}}}\quad\quad\mathrm{LR(LL}(G_{\mu\nu}^{(p')})).\end{array}\]

\bigskip{}
\noindent We adopt the name corridor $G$-set for the collection
of $G_{\mu\nu}^{(p)}$ that potentially partake in additive partitions
realizing kissing numbers that reside in the interordinal corridor
$]\Sigma_{i=1}^{(q+1)/2}G_{q+1,i}^{((p-1)/2)},C_{q}[$ . For $p=15,31$,
the corridor $G$-set is $G_{\textrm{cor}}^{(15)}=\{[1],(3)\}$, $G_{\textrm{cor}}^{(31)}=\{[1,3,5,11,17,41],(19,43),115,155,429\}$.
\\
Little is known for certain about higher kissing numbers, but Table
\ref{tab:Stenoscopy-of-kissing} gives valuable hints on not directly
accessible details,%
\footnote{the table favors the latest results by Nebe (summarized in \cite{Nebe11})
over older ones regarding Mordell-Weil lattice, Barnes-Wall lattice%
} especially about how many of them would belong to $\mathrm{LL}(G_{\mu\nu}^{(31)})$
and how many to $\mathrm{LL}(G_{\mu\nu}^{(63)})$ \emph{etc}.%
\begin{table}

\caption{\label{tab:Stenoscopy-of-kissing}Stenoscopy of kissing numbers of
Euclidean $D$-space relative to Catalan numbers and least-row-of-LL
sums $\Sigma_{i=1}^{q+1}G_{2q+2,i}^{(p)}$}

\lyxline{\normalsize}

$\; D$$\qquad\qquad\qquad\qquad\qquad\qquad\qquad\qquad\qquad\qquad$$L_{D}\qquad\qquad\qquad\qquad\quad$
bounding Catalan numbers and least-row-of-LL sums

\lyxline{\normalsize}\medskip{}
$\begin{array}{crc}
\,1 & \qquad\qquad\qquad\qquad\qquad\qquad\qquad\qquad\qquad\qquad\:2 & \qquad\qquad\qquad\qquad\qquad\qquad\leq C_{2}\boldsymbol{<C_{q}=C_{3}=5}\end{array}$

\vspace{-1mm}
- - - - - - - - - - - - - - - - - - - - - - - - - - - - - - - - -
- - - - - - - - - - - - - - - - - - - - - - - - - - - - - - - - -
- \vspace{-1mm}

$\begin{array}{crc}
\,2 & \qquad\qquad\qquad\qquad\qquad\qquad\qquad\qquad\qquad\qquad\:6 & \qquad\qquad\qquad\qquad\qquad\qquad\qquad\boldsymbol{>C_{3}=5}\\
\,3 & 12 & \qquad\qquad\qquad\qquad\qquad\qquad\qquad<C_{4}=14\\
\,4 & 24 & \qquad\qquad\qquad\qquad\qquad\qquad\qquad\\
\,5 & 40 & \qquad\qquad\qquad\qquad\qquad\qquad\qquad<C_{5}=42\\
\,6 & 72 & \qquad\qquad\qquad\qquad\qquad\qquad\qquad\\
\,7 & 126 & \qquad\qquad\qquad\qquad\leq C_{6}(\!\,=132)\leq\Sigma_{i=1}^{4}G_{8,i}^{(15)}(\!\,=170)\end{array}$

\vspace{-1mm}

- - - - - - - - - - - - - - - - - - - - - - - - - - - - - - - - -
- - - - - - - - - - - - - - - - - - - - - - - - - - - - - - - - -
- \vspace{-1mm}

$\begin{array}{crc}
\,8 & \qquad\qquad\qquad\qquad\qquad\qquad\qquad\qquad\qquad\enskip240 & \;\qquad\qquad\qquad\qquad\qquad\qquad\qquad>\Sigma_{i=1}^{4}G_{8,i}^{(15)}\\
9 & \qquad\qquad\qquad\qquad\qquad\qquad\enskip\quad(\textrm{\textrm{\ensuremath{{\scriptstyle \textrm{Leech lattice}}}}})\:272\\
10 & 336 & \qquad\qquad\qquad\qquad\qquad\qquad\quad\boldsymbol{<C_{q'}=C_{7}=}\boldsymbol{429}\end{array}$

. . . . . . . . . . . . . . . . . . . . . . . . . . . . . . . . .
. . . . . . . . . . . . . . . . . . . . . . . . . . . . . . . . .
. . . . . . . 

$\begin{array}{crc}
11 & \qquad\qquad\qquad\qquad\qquad\qquad\qquad\qquad\qquad\quad438 & \qquad\qquad\qquad\qquad\qquad\qquad\qquad\boldsymbol{>C_{7}=}\boldsymbol{429}\end{array}$

$\begin{array}{crc}
12 & 648 & \qquad\qquad\qquad\qquad\qquad\qquad\qquad\\
13 & 906 & \quad\qquad\qquad\qquad\qquad\qquad\qquad\\
14 & 1422 & \qquad\qquad\qquad\qquad\qquad\qquad\qquad<C_{8}=1430\\
15 & 2340 & \quad\qquad\qquad\qquad\qquad\qquad\qquad\\
16 & \qquad\qquad\qquad\qquad\qquad\qquad\qquad\:({\scriptstyle \textrm{certified}})\:4320 & \qquad\qquad\qquad\qquad\qquad\qquad\qquad<C_{9}=4862\\
? & ? & \qquad\qquad\qquad\qquad\qquad\qquad\qquad<C_{10}=16796\\
? & ? & \qquad\qquad\qquad\qquad\qquad\qquad\qquad<C_{11}=58786\\
24 & ({\scriptstyle \textrm{certified}})\:196560 & \qquad\qquad\qquad\qquad\qquad\qquad\qquad<C_{12}=208012\\
? & ? & \qquad\qquad\qquad\qquad\qquad\qquad\qquad<C_{13}=742900\\
28 & ? & \qquad\qquad\qquad\qquad\qquad\qquad\qquad\leq C_{14}=2674440\end{array}$

\vspace{-1mm}
. . . . . . . . . . . . . . . . . . . . . . . . . . . . . . . . .
. . . . . . . . . . . . . . . . . . . . . . . . . . . . . . . . .
. . . . . . .

\vspace{-1mm}
$\begin{array}{crc}
29 & \:\;\quad\qquad\qquad\qquad\qquad\qquad\qquad\qquad\qquad\qquad? & \qquad\qquad\qquad\qquad\qquad\qquad\qquad>C_{14}=2674440\\
30 & \qquad\qquad\qquad\qquad\qquad\qquad\enskip\quad?\\
31 & ? & \qquad\qquad\qquad\qquad\qquad\quad\leq\Sigma_{i=1}^{8}G_{16,i}^{(31)}(\!\,=3437984)\end{array}$

\vspace{-1mm}
- - - - - - - - - - - - - - - - - - - - - - - - - - - - - - - - -
- - - - - - - - - - - - - - - - - - - - - - - - - - - - - - - - -
-

\vspace{-1mm}
$\begin{array}{crc}
32 & \:\;\quad\qquad\qquad\qquad\qquad\qquad\qquad\qquad\qquad\qquad? & \qquad\qquad\qquad\qquad\qquad\qquad\qquad>\Sigma_{i=1}^{8}G_{16,i}^{(31)}\\
? & \qquad\qquad\qquad\qquad\qquad\qquad\enskip\quad?\\
40 & ? & \qquad\qquad\qquad\qquad\qquad\boldsymbol{<C_{q''}=C_{15}=}\boldsymbol{9694845}\end{array}$

\vspace{-1mm}
. . . . . . . . . . . . . . . . . . . . . . . . . . . . . . . . .
. . . . . . . . . . . . . . . . . . . . . . . . . . . . . . . . .
. . . . . . .

\vspace{-1mm}
$\begin{array}{crc}
41 & \:\quad\qquad\qquad\qquad\qquad\qquad\qquad\qquad\qquad\qquad? & \qquad\qquad\qquad\qquad\qquad\boldsymbol{>C_{q''}=C_{15}=}\boldsymbol{9694845}\\
? & \qquad\qquad\qquad\qquad\qquad\qquad\quad?\\
48 & ({\scriptstyle \textrm{Nebe}})\:52416000 & \qquad\qquad\qquad\qquad\qquad\qquad<C_{17}=129644790\\
? & ?\\
72 & ({\scriptstyle \textrm{Nebe}})\:6218175600 & \qquad\qquad\qquad\qquad\qquad\qquad<C_{20}=6564120420\\
? & ?\\
\!102 & ? & \qquad\qquad\qquad\qquad\qquad\leq C_{30}=3814986502092304\end{array}$

\vspace{-1mm}
. . . . . . . . . . . . . . . . . . . . . . . . . . . . . . . . .
. . . . . . . . . . . . . . . . . . . . . . . . . . . . . . . . .
. . . . . . .

\vspace{-1mm}
$\begin{array}{crc}
\!103 & \qquad\qquad\qquad\qquad\qquad\qquad\qquad\qquad\qquad\quad\:? & \qquad\qquad\qquad\qquad\qquad\qquad\qquad\qquad>C_{30}\\
? & ?\end{array}$

\medskip{}
 \lyxline{\normalsize}
\end{table}

\begin{conj}
\label{con:In--(or}Let $R_{a}^{(n)}$ $(n=\log_{2}(p+1))$ be the
\# of kissing numbers in $]\Sigma_{i=1}^{(q+1)/2}G_{q+1,i}^{((p-1)/2)},C_{q}[$
, $R_{b}^{(n)}$ the \# of those in $]C_{q},C_{2q}]$ , and $R_{c}^{(n)}$
the \# of those in $]C_{2q},\Sigma_{i=1}^{q+1}G_{2q+2,i}^{(p)}]$
, all representable by suitably chosen additive partitions from rows
(columns) in $\mathit{\mathrm{LL}}(G_{\mu\nu}^{(p)})$ (or $\mathit{\mathrm{LL}}(J_{\mu\nu}^{(p)})\,$).
Then\begin{equation}
\begin{array}{c}
R_{a}^{(n)}=R_{c}^{(n)}=\frac{1}{2}\mathcal{T}_{(p-1)/2}\\
R_{b}^{(n)}=\mathcal{T}_{p}\end{array}\qquad(n>4),\label{eq:add-part}\end{equation}
\textup{\noindent where} $R_{a}^{(n)}$ (marked by dashed-line, dotted-line
delimiters in Table \ref{tab:Stenoscopy-of-kissing}) determines the
corridor $G$- ($J$)-set, $G_{\textrm{cor}}^{(p)}$ ($J_{\textrm{cor}}^{(p)}$).
\end{conj}
\noindent Case $p=7$ is degenerate, with interval $]1,2]$ harboring
one kissing number, $L_{1}=2$. We set $R_{a}^{(3)}=R_{c}^{(3)}=0,R_{b}^{(3)}=1.$

\noindent Case $p=15$: No kissing number lives in $]2,5]$ : $R_{a}^{(4)}=0$
-- which is equivalent to saying the interordinal corridor $G$-set
$G_{\textrm{cor}}^{(15)}$ has no unbracketed, unparenthesized entries
--, and in $]132,170]$ live none either: $R_{c}^{(4)}=0$, thus the
case is degenerate, too. Since interval $]5,132]$ harbors six kissing
numbers -- 6,12,24,40,72,126 --, we set $R_{b}^{(4)}=6$.

\noindent Case $p=31$: In $]170,429[$ we find $R_{a}^{(5)}=3$
kissing numbers, indicated by the presence of unbracketed, unparenthesized
entries in the corridor $G$-set $G_{\textrm{cor}}^{(31)}$. In $]429,2674440]$
, $R_{b}^{(5)}=18$, and in $]2674440,3437984]$ , $R_{c}^{(5)}=3$
(most of them uncertified).

\noindent Case $p=63$: In $]3437984,9694845[$ live $R_{a}^{(6)}=9$
kissing numbers, in $]9694845,C_{30}]$ $R_{b}^{(6)}=54$ and in $]C_{30},\Sigma_{i=1}^{16}G_{32,i}^{(63)}]$
$R_{c}^{(6)}=9$. 

\noindent And so on for the cases $p=127,255,\ldots$ 

\noindent All partitions realizing kissing numbers may equivalently
be defined over columns, and the reader is invited to identify the
corresponding rows and columns in $\mathrm{LL}(J_{\mu\nu}^{(p)})$,
too. 

\pagebreak{}

\subsection{\label{sub:The-factorization-aspect}The factorization aspect}

\medskip{}
\noindent The question of whether the partial sequences $(G^{(3)})=(1),\;(G^{(7)})=(1),\;(G_{\rho}^{(15)})=(\,3,\!5,\!11,\!17,\!41,\!113\,),\,...\,$
always consist of prime-numbered representatives at paraorder fifteen
or higher may have seemed intriguing at first: after all, the expressions
$B_{p}^{2}-2$ themselves started out with primes for $p=7,15,31,63$.
Notable as these facts may be, we have afterwards seen that the bulk
of $G_{\rho}^{(p)}$ do \emph{not} stay prime as we go along with
the computation of the next higher members of the root-$f$ sequence.
In fact, the prime-numbered $G_{\rho}^{(p)}$ decrease in number --
six in $\mathrm{\neg UR\mathrm{(LL}(\mathit{G_{\mu\nu}^{\mathrm{(15)}}}))}$,
four in $\mathrm{\mathrm{\neg UR(LL}(\mathit{G_{\mu\nu}^{\mathrm{(31)}}})})$
--, circumstantial evidence that might indicate a trend. For a first,
identity (\ref{eq:inter}) provides an opportunity for a quick prime
number test. For $p\!=\!15$, $p'\!=\!31$, Eq. (\ref{eq:inter})
furnishes 19 und 43, the paraorder-thirty-one twins of the paraorder-fifteen
primes 17 und 41. And for para\-orders thirty-one and sixty-three,
this equation yields \begin{equation}
\mathrm{\mathrm{UR(UL(LL}(\mathit{G_{\mu'\nu'}^{(\mathrm{63)}}})))=\left(\mathit{G}_{16+\xi,4+\zeta}^{(63)}\right)}=\left(\begin{array}{rrrr}
58791 & 18633 & 4907 & 1635\\
189393 & 58791 & 15299 & 4907\\
738035 & 227123 & 58791 & 18633\\
2430515 & 738035 & 189393 & 58791\end{array}\right),\label{eq:search-prime-63}\end{equation}
\\
 adding one additional prime, $\pi_{1787}=15299$. But prime-numbered
$G_{\rho'}^{(63)}$ can spring from any other quad\-rant of $\mathrm{\neg UR(LL(\mathit{G}_{\mu\nu}^{(63)}))}$.
However small the increase in knowledge to expect from the endeavor,
we thought the question intriguing enough to undertake a complete
{\tt int64} computation of all 54 representatives of $\mathrm{\neg UR(LL(\mathit{G}_{\mu\nu}^{(63)}))}$
$-$ and found three more primes, $\pi_{2\,364\,489}=38\,792\,251$,
$\pi_{?}=69\,531\,783\,535\,237$ and $\pi_{?}=283\,858\,869\,110\,417$.
If primes refuse to go, maybe their number settles on 4 from Mersennian
paraorder 31 on.

\noindent More generally, the quantities $G_{\rho}^{(p)}>1$ can
be classified by their factorization. We distinguish pure prime numbers
$\pi_{r}$, factorization into two or three prime factors, $\pi_{r}\!\cdot\!\pi_{s}$
and $\pi_{r}\cdot\pi_{s}\cdot\pi_{t}$, as well as factorization into
one or more prime factors exponentiated $\pi_{r}^{z_{r}}\!(\cdot\pi_{s}^{z_{s}}\cdot\dots),$
$z_{r}\!>\!1$ $(\vee\: z_{s}\!>\!1\cdots)$. Only the conditions
up to $p=63$ have been analyzed; so far, there's nothing that contradicts
the assumption that the \# of factorization types, with more complex
ones pooled as $\pi_{r}\cdot\pi_{s}\cdot\pi_{t}\cdot\ldots$ and $\pi_{r}^{z_{r}}\cdot\ldots$,
respectively, stays as even-numbered as it turns out to be for $\mathrm{\mathrm{\neg UR(LL}(\mathit{G_{\mu\nu}^{(p)}})),\;\mathit{p}=15,31,63}$:%
\begin{table}[H]
\caption{$G_{\rho}$ according to factorization type \smallskip{}
 \label{tab:Factorization-types-1}}

\lyxline{\normalsize}\# $G_{\rho}$ fact'd as $\mathrm{\qquad\qquad\qquad\qquad\qquad\neg UR(LL}(\mathit{G_{\mu\nu}^{\mathrm{(15)}}}))$$\mathrm{\qquad\qquad\qquad\neg UR(LL}(\mathit{G_{\mu\nu}^{\mathrm{(31)}}}))$$\mathrm{\qquad\qquad\qquad\neg UR(LL}(\mathit{G_{\mu\nu}^{\mathrm{(63)}}}))$
\lyxline{\normalsize}

$\begin{array}{cccc}
\pi_{r} & \qquad\qquad\qquad\:\qquad\qquad\qquad6 & \qquad\qquad\qquad\qquad\qquad\qquad4 & \qquad\qquad\qquad\qquad\qquad\qquad4\\
\pi_{r}\cdot\pi_{s} & \qquad\qquad\qquad\:\qquad\qquad\qquad\hbox{-} & \qquad\qquad\qquad\qquad\qquad\qquad6 & \qquad\qquad\qquad\qquad\qquad\qquad16\\
\pi_{r}\cdot\pi_{s}\cdot\pi_{t} & \qquad\qquad\qquad\:\qquad\qquad\qquad\hbox{-} & \qquad\qquad\qquad\qquad\qquad\qquad6 & \qquad\qquad\qquad\qquad\qquad\qquad14\\
\pi_{r}\cdot\pi_{s}\cdot\pi_{t}\cdot\ldots & \qquad\qquad\qquad\:\qquad\qquad\qquad\hbox{-} & \qquad\qquad\qquad\qquad\qquad\qquad\hbox{-} & \qquad\qquad\qquad\qquad\qquad\qquad8\\
\pi_{r}^{z_{r}}\cdot\ldots & \qquad\qquad\qquad\:\qquad\qquad\qquad\hbox{-} & \qquad\qquad\qquad\qquad\qquad\qquad2 & \qquad\qquad\qquad\qquad\qquad\qquad12\\
\sum & \qquad\qquad\qquad\:\qquad\qquad\qquad6 & \qquad\qquad\qquad\qquad\qquad\qquad18 & \qquad\qquad\qquad\qquad\qquad\qquad54\end{array}$

\smallskip{}
 \lyxline{\normalsize} \smallskip{}

\end{table}

\noindent One further observation is that composite $G_{\rho}^{(31)}$
-- indeed the bulk of them -- are missing pure prime numbers with
minimal spacings that lie in the same range as the \# of $G_{\rho}^{(31)}$
congruent with $(7-2k)(\mathrm{mod}\:8)$ $(k=1,2,3)$ -- see Table
\ref{tab:Extended-partioning-1-3-5-7} for a summary. Conversely,
just as a $G_{\mu\nu}$ congruent with 7 modulo 8 is absent from $\neg\mathrm{UR}\mathrm{(LL\mathit{f}}^{(31)})$,
so is the minimal spacing 14 involving the factor 7:%
\begin{figure}[H]

\caption{\label{fig:-interpolating-prime}$G_{\rho}^{(31)}$ interpolating
adjoining prime numbers}

\[
\begin{array}{rlllll}
\quad19 & = & \pi_{8} &  &  & \qquad\mathrm{UL}(\mathrm{LL}\,(G_{\mu\nu}^{(31)}))\\
\quad43 & = & \pi_{14} &  &  & \qquad\quad\downarrow\;\;\\
\quad115 & = & \pi_{30}+2 & \quad\qquad(= & \pi_{31}-12)\\
\quad155 & = & \pi_{36}+4 & \quad\qquad(= & \pi_{37}-2)\\
\quad429 & = & \pi_{82}+8 & \quad\qquad(= & \pi_{83}-2)\\
\quad1275 & = & \pi_{205}+16 & \quad\qquad(= & \pi_{206}-2)\\
\quad1595 & = & \pi_{250}+12 & \quad\qquad(= & \pi_{251}-2)\\
\quad4819 & = & \pi_{649}+2 & \quad\qquad(= & \pi_{650}-12)\\
\quad15067 & = & \pi_{1759}+6 & \quad\qquad(= & \pi_{1760}-6)\end{array}\]

\[
\begin{array}{rlllll}
1633 & = & \pi_{258}+6 & \:\qquad(= & \pi_{259}-4) & \mathrm{\quad\; LL(LL\mathit{\mathrm{(}G_{\mu\nu}^{\mathrm{(31)}}}))}\\
4905 & = & \pi_{655}+2 & \:\qquad(= & \pi_{656}-4) & \quad\;\quad\downarrow\;\;\\
15297 & = & \pi_{1786}+8 & \:\qquad(= & \pi_{1787}-2)\\
18627 & = & \pi_{2129}+10 & \:\qquad(= & \pi_{2130}-10)\\
58781 & = & \pi_{5946}+10 & \:\qquad(= & \pi_{5947}-6)\\
189371 & = & \pi_{17110}+10 & \:\qquad(= & \pi_{17111}-6)\\
227089 & = & \pi_{20185}\\
737953 & = & \pi_{59377}+24 & \:\qquad(= & \pi_{59378}-16)\\
2430289 & = & \pi_{178344}\end{array}\]

\end{figure}

\noindent At the interval in question, prime numbers are relatively
close to one another, so instead of surmising some lawfulness behind
this phenomenon, suffice it to say in this section that the \# of
$G_{\rho}\equiv(7-2k)\:(\mathrm{mod}\:8)$$(k=0,1,2,3)$ (as shown
in Table \ref{tab:Extended-partioning-1-3-5-7}) and the prime-number
interpolations seem to follow a common structural pattern.

\noindent \noindent Let us now address a phenomenon that we considered
important enough to assign its characteristic order a separate letter,
$q$, consistently meaning $(p-3)/4$. As we have seen, there are
places where the coefficients $G_{\rho}^{(p)}$ directly intersect
with the root-$h$ associated coefficients $J_{\omega}^{(p)}$ times
minus one -- demonstrated in Fig. \ref{fig:catal-rep} as framed entries
for $\mathrm{LL}\,(G_{\mu\nu}^{(31)})$ and $\mathrm{LL}\,(J_{\mu\nu}^{(31)})$
respectively; the appendant entries for $\mathrm{LL}\,(G_{\mu\nu}^{(p)})$
and $\mathrm{LL}\,(J_{\mu\nu}^{(p)})$ for $p=15$, $p=7$ are set
off in the bracketed parts, upper right:%
\begin{figure}[H]
\caption{\label{fig:catal-rep}Catalan representative $G_{q+2,1}^{(p)}$ constituting
SCPF primes}
\[
\mathrm{LL}\,(G_{\mu\nu}^{(31)})=\left(\begin{array}{cc}
\begin{array}{rrrr}
\fbox{429} & 155 & 43 & 19\\
1275 & \fbox{429} & 115 & 43\\
4819 & 1595 & \fbox{429} & 155\\
15067 & 4819 & 1275 & \fbox{429}\end{array} & \left[\begin{array}{cc}
\begin{array}{cc}
\fbox{5} & 3\\
11 & \fbox{5}\end{array} & \left[\begin{array}{cc}
\fbox{1} & 1\\
1 & \fbox{1}\end{array}\right]\\
\begin{array}{cc}
\!\!\!\!41 & \quad17\\
\!\!\!113 & \quad41\end{array} & \begin{array}{cc}
\fbox{5} & 3\\
11 & \fbox{5}\end{array}\end{array}\right]\\
\begin{array}{rrrr}
58781 & 18627 & 4905 & 1633\\
189371 & 58781 & 15297 & 4905\\
737953 & 227089 & 58781 & 18627\\
2430289 & 737953 & 189371 & 58781\end{array} & \begin{array}{cccc}
\fbox{429} & 155 & 43 & 19\\
1275 & \fbox{429} & 115 & 43\\
4819 & 1595 & \fbox{429} & 155\\
15067 & 4819 & 1275 & \fbox{429}\end{array}\end{array}\right),\]
\medskip{}
 \[
\mathrm{LL}\,(J_{\mu\nu}^{(31)})=\left(\begin{array}{cc}
\begin{array}{rrrr}
\fbox{-429} & 117 & -41 & 13\\
1547 & \fbox{-429} & 143 & -41\\
-4903 & 1343 & \fbox{-429} & 117\\
18269 & -4903 & 1547 & \fbox{-429}\end{array} & \left[\begin{array}{cc}
\begin{array}{cc}
\fbox{-5} & 1\\
15 & \fbox{-5}\end{array} & \left[\begin{array}{cc}
\fbox{-1} & 1\\
3 & \fbox{-1}\end{array}\right]\\
\begin{array}{cc}
\!\!\!\!-43 & \quad15\\
\!\!\!149 & \quad-43\end{array} & \begin{array}{cc}
\fbox{-5} & 1\\
15 & \fbox{-5}\end{array}\end{array}\right]\\
\begin{array}{rrrr}
-58791 & 15547 & -4823 & 1319\\
223573 & -58791 & 17989 & -4823\\
-747765 & 194993 & -58791 & 15547\\
2886235 & -747765 & 223573 & -58791\end{array} & \begin{array}{cccc}
\fbox{-429} & 117 & -41 & 13\\
1547 & \fbox{-429} & 143 & -41\\
-4903 & 1343 & \fbox{-429} & 117\\
18269 & -4903 & 1547 & \fbox{-429}\end{array}\end{array}\right).\]

\end{figure}

\noindent  The associated Catalan number $C_{q}$ takes a special
role here, \emph{viz}. \begin{equation}
G_{q+1+\xi,\xi}^{(p)}+J_{q+1+\zeta,\zeta}^{(p)}=C_{q}-C_{q}=0\qquad\quad(\xi,\zeta\in\{1,\dots,q+1\}).\end{equation}
This relation is part of a larger underlying symmetry: while the quadrant
sum $\mathrm{LL}(G_{\mu\nu}^{(p)})+\mathrm{LL}(J_{\mu\nu}^{(p)})$
is secondary symmetric and $(n-2)$-fold traceless, one trace vanishing
main-, the next to the right submain-, and so on, the Lie bracket
of the summands is secondary antisymmetric,\[
[\mathrm{LL}(G_{\mu\nu}^{(p)}),\mathrm{LL}(J_{\mu\nu}^{(p)})]=-[\mathrm{LL}(G_{\mu\nu}^{(p)}),\mathrm{LL}(J_{\mu\nu}^{(p)})]^{S},\]

\noindent and $(p_{n-1}-1)$-fold traceless, with $p_{n-1}-2$ traces
vanishing main- and adjacent- on either side due to secondary antisymmetry,
and one secondary- due to secondary-diagonal zeroing. We may single
out the upper left coefficient $G_{q+2,1}^{(p)}$, say, to get what
may be called an overarching Catalan representative of $(G_{\rho}^{(p)})$
and $(J_{\omega}^{(p)})$. What makes the latter unique is that, from
$p=15$ on, it displays a peculiar type of factorization that involves
a $\mathbf{\underline{s}}$uffix of $\mathbf{\underline{c}}$onsecutive
$\mathbf{\underline{p}}$rime $\mathbf{\underline{f}}$actors (SCPF)
lying in the interval $]q+1,2q[$. To wit, 

\noindent  $G_{5,1}^{(15)}=C_{3}$ corresponds to the suffix (underlined)
\[
\underline{5}\,,\]
while to $G_{9,1}^{(31)}=C_{7}$ there belongs 

\[
429=3\cdot\underline{11\cdot13}\,,\]
 and to $G_{17,1}^{(63)}=C_{15}$\[
9694845=3^{2}\cdot5\cdot\underline{17\cdot19\cdot23\cdot29}\,,\]
 \\
 followed by%
\footnote{special thanks go to wolframalpha.com through whose good offices larger
Catalan numbers have now become widely accessible%
} $\mathit{G}_{33,1}^{(127)}=\mathit{C}_{31}=$ \[
\mathrm{7\cdot11\cdot17\cdot19\cdot\underline{37\cdot41\cdot43\cdot47\cdot53\cdot59\cdot61}}\,,\]
\bigskip{}
\negthinspace{}\noindent further followed by $\mathit{G}_{65,1}^{(255)}=\mathit{C}_{63}=$

\[
\mathrm{3\cdot5^{3}\cdot11^{2}\cdot\ldots\cdot41\cdot\underline{67\cdot71\cdot73\cdot79\cdot83\cdot89\cdot97\cdot101\cdot103\cdot107\cdot109\cdot113}}\,,\]
\[
\]

\noindent and so on, which can be summarized in the multiplicative
Euler-product partition

\[
2\cdot\prod_{{\scriptstyle {\textstyle _{p\,\textrm{Mersenne prime}}}}}\!\negthinspace\!\!\negthinspace\!\! p\:\cdot\prod_{\; p=15,31,\ldots}\mathrm{\: SCPF}(G_{q+2,1}^{(p)})\;=\;\prod\!\pi_{r}\,.\]
\[
\]

\noindent With the denotation $\mathrm{SCPF}_{p}:\textrm{the set of prime factors contained in}\mathrm{\: SCPF}(G_{q+2,1}^{(p)})\;(p=15,31,\ldots),$
the set of all prime numbers becomes the disjoint union of the singleton
$\{2\}$, the Mersenne prime numbers and the $\mathrm{SCPF}_{p}$'s.
According to the prime-number density theorem, the number of factors
contained in $\mathrm{\: SCPF}_{p}$, denoted $S_{p}$ here, is of
order $\frac{2q}{\log(2q)}-\frac{q+1}{\log(q+1)}$. In Table \ref{tab:Order-()-numbers},
the values of $S_{p}$ for $p=15,31,\ldots$ are listed together with
two other order-($\frac{2q}{\log(2q)}-\frac{q+1}{\log(q+1)}$) numbers.
Where $n_{q}=\log_{2}(q+1)$: the Catalan numbers of half-integer
index, $C_{1+n_{q}/2}$, defined by $\frac{2^{2(1+n_{q}/2)}\Gamma(1+(1+n_{q})/2)}{\sqrt{3.14\cdots}\,\Gamma(3+n_{q}/2)}$,
and the kissing numbers $L_{n_{q}-3}$.%
\footnote{$\,$Strictly speaking, only lower and upper bounds are known for
them in some places -- e.g. (40,44), (72,78) and (126,134); but, according
to prevailing knowledge, they do stop being order-($\frac{2q}{\log(2q)}-\frac{q+1}{\log(q+1)}$)
numbers after dimension eight. %
}

\medskip{}
\begin{table}[H]
\caption{{\small \label{tab:Order-()-numbers}Orde}r-($\frac{2q}{\log(2q)}-\frac{q+1}{\log(q+1)}$)
{\small numbers} ${\scriptstyle {\textstyle \qquad\qquad\qquad\qquad\qquad\qquad\qquad\qquad(n_{q}=\log_{2}(q+1),\: p=4q+3)}}$}

\medskip{}
\lyxline{\normalsize}

\noindent $\quad q$$\quad\qquad\qquad\qquad$$\frac{2q}{\log(2q)}-\frac{q+1}{\log(q+1)}$$\qquad\enskip\qquad\qquad S_{p}$$\qquad\qquad\quad\qquad\qquad C_{1+n_{q}/2}$$\qquad\qquad\qquad\qquad L_{n_{q}-3}$

\lyxline{\normalsize}\vspace{-5mm}
$\begin{array}{cccccc}
 & \qquad\qquad & \qquad\qquad & \qquad\qquad\qquad\qquad & \qquad\qquad\qquad\qquad & \quad\qquad\qquad\qquad\\
3 & \qquad & \qquad\qquad0.46 & \qquad\qquad\qquad\qquad1 & \qquad\qquad\qquad\qquad2 & \quad\qquad\qquad\qquad-\\
7 & \qquad & \qquad\qquad1.45 & \qquad\qquad\qquad\qquad2 & \qquad\qquad\qquad\qquad3.10 & \quad\qquad\qquad\qquad-\\
15 & \qquad & \qquad\qquad3.04 & \qquad\qquad\qquad\qquad4 & \qquad\qquad\qquad\qquad5 & \quad\qquad\qquad\qquad2\\
31 & \qquad & \qquad\qquad5.78 & \qquad\qquad\qquad\qquad7 & \qquad\qquad\qquad\qquad8.27 & \quad\qquad\qquad\qquad6\\
63 & \qquad & \qquad\qquad10.66 & \qquad\qquad\qquad\qquad12 & \qquad\qquad\qquad\qquad14 & \quad\qquad\qquad\qquad12\\
127 & \qquad & \qquad\qquad19.48 & \qquad\qquad\qquad\qquad23 & \qquad\qquad\qquad\qquad24.08 & \quad\qquad\qquad\qquad24\\
255 & \qquad & \qquad\qquad35.63 & \qquad\qquad\qquad\qquad43 & \qquad\qquad\qquad\qquad42 & \quad\qquad\qquad\qquad(40,44)\\
511 & \qquad & \qquad\qquad65.41 & \qquad\qquad\qquad\qquad75 & \qquad\qquad\qquad\qquad74.09 & \quad\qquad\qquad\qquad(72,78)\\
1023 & \qquad & \qquad\qquad120.64 & \qquad\qquad\qquad\qquad137 & \qquad\qquad\qquad\qquad132 & \quad\qquad\qquad\qquad(126,134)\\
2047 & \qquad & \qquad\qquad223.62 & \qquad\qquad\qquad\qquad255 & \qquad\qquad\qquad\qquad237.11 & \quad\qquad\qquad\qquad240\\
4095 & \qquad & \qquad\qquad419.48 & \qquad\qquad\qquad\qquad463 & \qquad\qquad\qquad\qquad429 & \quad\qquad\qquad\qquad\times\end{array}$ \lyxline{\normalsize}\bigskip{}

\end{table}
\medskip{}
\noindent With the definition \[
\textrm{SCPF primes : }\enskip\bigcup_{p=15,31,\ldots}\mathrm{\! SCPF}_{p}\]

\noindent we face a dilemma: they almost form the class of all prime
numbers that lie between two consecutive Mersenne numbers, just like,
where $C_{\textrm{-}1/2}=0$, $\mathcal{C}=(C_{\textrm{-}1/2};C_{\textrm{-}1/2};C_{\textrm{-}1/2};C_{2};C_{4},\ldots,C_{6};C_{8},\ldots,C_{14};\ldots)$
almost forms the class of all Catalan numbers of non-Mersenne-numbered
index that ensue from net traces over the secondary and adjacent diagonals
of $\mathrm{LL}(G_{\mu\nu}^{(p)})$. There clearly is one element
missing in either case. Regarding the $\textrm{SCP}\textrm{F}_{p}$'s:
the prime number 2 lying between the Mersenne numbers 1 and 3, and
regarding the sequence $\mathcal{C}$: the one net trace belonging
to $\mathrm{LL}(G_{\mu\nu}^{(3)})$. Incorporating the missing items,
we respectively get \begin{equation}
\textrm{SCPF}{}^{+}\textrm{ prime numbers : SCPF prime numbers }\cup\left\{ 2\right\} ,\label{eq:SCPFplus}\end{equation}
and\begin{equation}
\mathcal{C}^{+}=(C_{\textrm{-}1/2};C_{\textrm{-}1/2};C_{\textrm{-}1/2},C_{1};C_{2};C_{4},\ldots,C_{6};C_{8},\ldots,C_{14};\ldots),\label{eq:Upsilonplus}\end{equation}

\noindent with the effect that one even number is included among
otherwise odd numbers in the former case, and one odd number among
even numbers in the latter. 

\noindent The distinction created between SCPF and Mersenne primes
becomes vital when it comes to determining a possible set membership
of $G_{\mu\nu}^{(p)}$ in $\mathrm{SCPF}_{p_{N}}$, where $N$ denotes
the identity \begin{equation}
N={\scriptscriptstyle -}\left\lceil {\scriptscriptstyle {\textstyle (n+1)/{\scriptstyle 2}}}\right\rceil {\scriptscriptstyle +}\Sigma_{i=1}^{n-1}p_{i}=\left\lfloor \log_{2}C_{q'}\right\rfloor .\label{eq:log C identity}\end{equation}
\noindent Then Mersenne primes $>3$, it turns out, are not among
the factors of $p_{N}$, but the factors of $q_{N}=p_{N-2}$. Let
us, for example, check the membership of $G_{\textrm{max}}^{(p)}$
in $\mathrm{SCPF}_{p_{N}}$ for the few known cases: \[
G_{\textrm{max}}^{(15)}=\pi_{6\cdot5}=113\in\mathrm{SCPF}_{2^{8}-1}\rightsquigarrow\quad p_{N}=2^{8}-1=3\cdot5\cdot17,\quad q_{N}=2^{6}-1=3^{2}\cdot\underbrace{7},\]

\noindent and\[
G_{\textrm{max}}^{(31)}=\pi_{6\cdot29724}=2430289\in\mathrm{SCPF}_{2^{23}-1}\rightsquigarrow\quad p_{N}=2^{23}-1=47\cdot178\,481,\quad q_{N}=2^{21}-1=\underbrace{7^{2}}\cdot\underbrace{127}\cdot\:337,\]

\noindent where Mersenne primes $>3$ are marked with underbraces.
The next instance, $G_{\textrm{max}}^{(63)}$, is a composite. We
expect its three prime factors $3,\:613\;\textrm{and}\;1,\!910,\!047,\!210,\!943\;$
to be elements of $\mathrm{SCPF}_{2^{53}-1}$ because we find the
conjectured pattern\[
p_{N}=2^{53}-1=6361\cdot69431\cdot20\,394\,401,\quad q_{N}=2^{51}-1=\underbrace{7}\cdot103\cdot\:2143\cdot11\,119\cdot\underbrace{131\,071}.\]

\subsection{\label{sub:The-modulo-eight-aspect}The modulo-eight aspect}

\noindent Taking the modulo-eight aspect into account allows us briefly
to resume the subject of Sects. \emph{2-3} to show how $f$-parafermi
algebra (Eqs. (\ref{eq:f-par1})-(\ref{eq:f-ar2})) can be made to
hold for $p\in\{15,31,\dots\}$. Thus, as odd numbers $2n+1$ come
with the identity $(2n+1)^{2}=8\sum n+1$ and with all $G_{\mu\nu}$
odd-numbered, one has $G_{\mu\nu}^{2}\equiv1(\mathrm{mod}\:8)$. Hence,
by Eq. (\ref{eq:G-squared}),

\noindent \begin{equation}
\frac{1}{2}[f_{0}^{+},f_{0}]+\sum_{\upsilon=1}^{(p-1)/2}[f_{\upsilon}^{+},f_{\upsilon}]_{\mathrm{mod}8}=\mathrm{diag}\textrm{\Large(}\frac{p}{2},\frac{p}{2}-1,\cdots,-\frac{p}{2}+1,-\frac{p}{2}\textrm{\Large)},\label{eq:f-par1alter}\end{equation}
 \begin{equation}
\sum_{\upsilon=0}^{(p-1)/2}[[f_{\upsilon}^{+},f_{\upsilon}]_{\mathrm{mod}8}\:,f_{\upsilon}]=-2f,\;\:\sum_{\upsilon=0}^{(p-1)/2}[[f_{\upsilon}^{+},f_{\upsilon}]_{\mathrm{mod}8}\:,f_{\upsilon}^{+}]=2f^{+}.\label{eq:f-ar2alter}\end{equation}
 \noindent But it's worthwhile to have a look at the very arrangement
of residues left by $G_{\mu\nu}$ after division by eight,%
\footnote{This is not to say that larger moduli are less important. One can
e.g. notice the interesting fact that $G_{\mathrm{max}}^{(15)}=113\equiv7^{2}(\mathrm{mod}\:64)$
and $G_{\mathrm{max}}^{(31)}=2430289\equiv9^{2}(\mathrm{mod}\:128)$.
The modulo-8 approach is chosen here because it is in agreement with
the closure effect that can spring from the group $\mathrm{Z_{2}^{\,3}}$
through its various isomorphic maps. For\ octonions, it marks the
loss of associativity of the hyper\-complex number system; for $f$-parafermi
algebra bar the modulo-8 approach, the loss of consistency. Plus,
kissing numbers $L_{n_{q}-3}$ with $3<n_{q}=\log_{2}(q+1)$ lose
their order-($\frac{2q}{\log(2q)}-\frac{q+1}{\log(q+1)}$)-number
characteristics after dimension eight (see Table \ref{tab:Order-()-numbers}).%
} \[
\begin{array}{lll}
\left(G_{\mu\nu}^{(7)}\right)_{\mathrm{mod}\,8}=\left(G_{\mu\nu}^{(7)}\right)=\left(\begin{array}{cccc}
0\\
1 & \quad0\\
1 & \quad1 & \quad0\\
1 & \quad1 & \quad1 & \quad0\end{array}\right), &  & \left(G_{\mu\nu}^{(15)}\right)_{\mathrm{mod}\,8}=\left(\begin{array}{cccccccc}
0\\
1 & \quad0\\
1 & \quad1 & \quad0\\
1 & \quad1 & \quad1 & \quad0\\
5 & \quad3 & \quad1 & \quad1 & \quad0\\
3 & \quad5 & \quad1 & \quad1 & \quad1 & \quad0\\
1 & \quad1 & \quad5 & \quad3 & \quad1 & \quad1 & \quad0\\
1 & \quad1 & \quad3 & \quad5 & \quad1 & \quad1 & \quad1 & \quad0\end{array}\right)\end{array},\ldots,\]
 which shows that underneath the overt secondary symmetry of $f^{(7)}$
the original main symmetry of $\mathrm{LL(\mathit{G_{\mu\nu}^{\mathrm{(7)}}})}$
exerts its influence on paraorders beyond that mark. Its persistence
in modulo-8 form, quadrantwise in LL, subquadrantwise in LLUL, LLLR,
URLL etc., makes clear how the heterotic variant of $f$-parafermi
algebra (Eqs. \ref{eq:f-var1}-\ref{eq:f-var2}) may be reshaped in
order to have it hold for $p\in\{15,31,\dots\}$, namely:\\
 \begin{equation}
[[(f^{\circ})^{+},f^{\circ}],g]=-2f^{\circ},\;\:[[(f^{\circ})^{+},f^{\circ}],g^{+}]=2(f^{\circ})^{+},\label{eq:f-var1alter}\end{equation}
 \begin{equation}
\begin{array}{c}
\sum_{\upsilon=0}^{(p-1)/2}\textrm{\Large(}\chi[(f^{\circ})_{\upsilon}^{+},(f^{\circ})_{\upsilon}]+\sigma\left([(f^{\circ})_{\upsilon}^{+},g_{\upsilon}]+[g_{\upsilon}^{+},(f^{\circ})_{\upsilon}]\right)+\tau\left([(f^{\circ})_{\upsilon},g_{\upsilon}]+[g_{\upsilon}^{+},(f^{\circ})_{\upsilon}^{+}]\right)+\gamma[g_{\upsilon},g_{\upsilon}^{+}]\textrm{\Large)}\\
=\:\mathrm{diag}\textrm{\Large(}\frac{p}{2},\frac{p}{2}-1,\cdots,-\frac{p}{2}+1,-\frac{p}{2}\textrm{\Large)},\end{array}\label{eq:f-var2alter}\end{equation}
 \\
 where $f\equiv f^{\circ}(\mathrm{mod}\:8)$, or explicitly, $(f^{\circ})^{(p)}=1\!\!\!\boldsymbol{1}^{\otimes n-1}\otimes f^{(1)}+(G_{\mu\nu}^{(p)})_{\mathrm{mod\,8}}\otimes c_{3}$
. 

\noindent Apart from its consequences for \emph{f}-parafermi algebra,
the persistence even of main symmetry in modulo-8 form allows a very
compact way of describing the LL part: \[
\mathrm{LL(\mathsf{\mathit{G_{\mu\nu}^{\mathrm{(}p\mathrm{)}}}})_{mod\,8}}=\mathrm{\! sym(d_{m}}({\scriptstyle {.\,.\atop .\,.}}),{\scriptstyle \ldots}\mathrm{,d_{s}}({\scriptstyle {.\,.\atop .\,.}})).\]
 Applied to the case $p=31$, say, the expression\begin{equation}
\mathrm{\mathrm{LL(\mathit{G_{\mu\nu}^{\mathrm{(31)}}})_{mod\,8}}=sym(d_{m}}({\scriptstyle {5\:3\atop 3\:5}}),({\scriptstyle {3\:3\atop 3\:3}}),({\scriptstyle {5\:3\atop 3\:5}})\mathrm{,d_{s}}({\scriptstyle {1\:1\atop 1\:1}}))\label{eq:LL31-mod8}\end{equation}
 can be read as a shorthand for the evolution\[
\begin{array}{ccc}
\mathrm{\underrightarrow{\qquad\textrm{d}_{\mathrm{m}}\qquad}}\qquad\left(\begin{array}{cccccccc}
5 & 3 & 3 & 3 & 5 & 3 & 1 & 1\\
3 & 5 & 3 & 3 & 3 & 5 & 1 & 1\\
 &  & 5 & 3\\
 &  & 3 & 5\\
 &  &  &  & 5 & 3\\
 &  &  &  & 3 & 5\\
 &  &  &  &  &  & 5 & 3\\
 &  &  &  &  &  & 3 & 5\end{array}\right) &  & \mathrm{\underrightarrow{\qquad\quad\textrm{d}_{\mathrm{s}}\qquad}}\qquad\left(\begin{array}{cccccccc}
5 & 3 & 3 & 3 & 5 & 3 & 1 & 1\\
3 & 5 & 3 & 3 & 3 & 5 & 1 & 1\\
 &  & 5 & 3 & 1 & 1\\
 &  & 3 & 5 & 1 & 1\\
 &  & 1 & 1 & 5 & 3\\
 &  & 1 & 1 & 3 & 5\\
1 & 1 &  &  &  &  & 5 & 3\\
1 & 1 &  &  &  &  & 3 & 5\end{array}\right)\\
\\\end{array}\]
\[
\begin{array}{ccc}
\mathrm{\underrightarrow{\textrm{main\, symmetry}}}\;\left(\begin{array}{cccccccc}
5 & 3 & 3 & 3 & 5 & 3 & 1 & 1\\
3 & 5 & 3 & 3 & 3 & 5 & 1 & 1\\
3 & 3 & 5 & 3 & 1 & 1\\
3 & 3 & 3 & 5 & 1 & 1\\
5 & 3 & 1 & 1 & 5 & 3\\
3 & 5 & 1 & 1 & 3 & 5\\
1 & 1 &  &  &  &  & 5 & 3\\
1 & 1 &  &  &  &  & 3 & 5\end{array}\right) &  & \mathrm{\underrightarrow{\textrm{secondary\, symmetry}}}\;\left(\begin{array}{cccccccc}
5 & 3 & 3 & 3 & 5 & 3 & 1 & 1\\
3 & 5 & 3 & 3 & 3 & 5 & 1 & 1\\
3 & 3 & 5 & 3 & 1 & 1 & 5 & 3\\
3 & 3 & 3 & 5 & 1 & 1 & 3 & 5\\
5 & 3 & 1 & 1 & 5 & 3 & 3 & 3\\
3 & 5 & 1 & 1 & 3 & 5 & 3 & 3\\
1 & 1 & 5 & 3 & 3 & 3 & 5 & 3\\
1 & 1 & 3 & 5 & 3 & 3 & 3 & 5\end{array}\right),\\
\\\end{array}\]
where it is understood that the last two steps are recursively repeated
on subquadrants etc. in case of positions left blank -- such as would
be the case with paraorder 63, 127 etc. \\
We have argued that representatives $G_{\rho}^{(p)}$ are to be
sought among those $G_{\mu\nu}^{(p)}$ that originate from $\mathrm{\neg UR(LL}\mathrm{\mathit{f}^{(\mathit{p})})}$.
So far, this yielded $G^{(3)}=1,G^{(7)}=1,(G_{\rho}^{(15)})=(3,5,11,17,41,113)$.
$7(\mathrm{mod}\:8)$-congruence did not occur among them, nor does
it by the new arrivals from $\mathrm{\neg UR(LL}\mathit{f}^{(\mathit{\mathrm{31}})})$:
only $(7-2k)\,(\mathrm{mod}\:8)$ $(k=1,2,3)$-congruence is to be
found among these. So the question arises: Does $7(\mathrm{mod}\:8)$-congruence
finally show up in $(G_{\rho}^{(63)})$? We stop short of listing
the entire $64\times64$ matrix $f^{(63)}$ as a quick inspection
of the first row of the LL part of $(G_{\mu\nu}^{(63)})$ already
answers the question in the affirmative:

\begin{equation}
\begin{array}{ccccccccccccccccc}
 & G_{17,1} & G_{17,2} & G_{17,3} & G_{17,4} & G_{17,5} & G_{17,6} & G_{17,7} & G_{17,8} & G_{17,9} & G_{17,10} & G_{17,11} & G_{17,12} & G_{17,13} & G_{17,14} & G_{17,15} & G_{17,16}\\
{\atop \mathrm{mod}\,8} & {9694845\atop 5} & {2926323\atop 3} & {747891\atop 3} & {230395\atop 3} & {58791\atop 7} & {18633\atop 1} & {4907\atop 3} & {1635\atop 3} & {429\atop 5} & {155\atop 3} & {43\atop 3} & {19\atop 3} & {5\atop 5} & {3\atop 3} & {1\atop 1} & {1\atop 1}\end{array}\label{eq:row-G17}\end{equation}
 \\
The reality of $G_{\rho}^{(p)}$ $\equiv$ $(7-2k)\,(\mathrm{mod}\,8)\:(k=0,1,2,3)$
$(p=63,127,\ldots)$ simply is a consequence of interordinal identity
(\ref{eq:inter}) applied modulo eight:

\begin{equation}
\mathrm{UR(UL(LL(\mathit{G_{\mu\nu}^{\mathrm{(63)}}})_{mod\,8}))}=\mathrm{sym(}({\scriptstyle {5,3\atop 3,5}})\mathrm{,d_{s}}({\scriptstyle {1,1\atop 1,1}}))+2\,\mathrm{sym(}({\scriptstyle {5,3\atop 3,5}}),\mathrm{d_{s}}({\scriptstyle {1,1\atop 1,1}}))\equiv\mathrm{sym(}({\scriptstyle {7,1\atop 1,7}}),({\scriptstyle {3,3\atop 3,3}}))\,(\mathrm{mod}\,8).\label{eq:inter-mod8}\end{equation}
 \\
The numbers of $G_{\rho}^{(p)}$ partitioned by congruence with
$(7-2k)\,(\mathrm{mod}\:8)$ $(k=0,1,2,3)$ up to $p=63$ are listed
in Table \ref{tab:Extended-partioning-1-3-5-7}. %
\begin{table}[H]
\caption{$G_{\rho}$ according to congruence with $(7-2k)\,(\mathrm{mod}\:8)$
$(k=0,1,2,3)$ up to $p=63$\smallskip{}
 \label{tab:Extended-partioning-1-3-5-7}}

\lyxline{\normalsize}\# $G_{\rho}$ cong't w/$\mathrm{\quad\quad\qquad\qquad\neg UR(LL}\mathit{f}^{(\mathit{\mathrm{15}})})$$\mathrm{\enskip\qquad\qquad\qquad\neg UR(LL}\mathit{f}^{(31)})$$\mathrm{\quad\qquad\qquad\qquad\neg UR(LL}f^{(63)})$
\lyxline{\normalsize}

$\begin{array}{cccc}
1(\mathrm{mod}\:8) & \qquad\qquad\:\qquad\qquad\qquad3 & \qquad\qquad\qquad\qquad\qquad\qquad6 & \qquad\qquad\qquad\qquad\qquad\qquad16\\
3(\mathrm{mod}\:8) & \qquad\qquad\:\qquad\qquad\qquad2 & \qquad\qquad\qquad\qquad\qquad\qquad10 & \qquad\qquad\qquad\qquad\qquad\qquad36\\
5(\mathrm{mod}\:8) & \qquad\qquad\:\qquad\qquad\qquad1 & \qquad\qquad\qquad\qquad\qquad\qquad2 & \qquad\qquad\qquad\qquad\qquad\qquad6\\
7(\mathrm{mod}\:8) & \qquad\qquad\:\qquad\qquad\qquad\hbox{-} & \qquad\qquad\qquad\qquad\qquad\qquad\hbox{-} & \qquad\qquad\qquad\qquad\qquad\qquad4\\
\sum & \qquad\qquad\:\qquad\qquad\qquad6 & \qquad\qquad\qquad\qquad\qquad\qquad18 & \qquad\qquad\qquad\qquad\qquad\qquad62\end{array}$\smallskip{}
 \lyxline{\normalsize} \bigskip{}

\end{table}

\noindent The composite map $\Lambda=(\mathrm{mod}\:8)\circ(\times3)$\emph{
}ensuing from

\emph{\[
\mathrm{LL(LL({\mathit{G_{\mu\nu}^{\mathrm{(}p\mathrm{)}}}})_{\mathrm{mod}\,8})}+2\mathrm{\, UR(LL(\mathit{G_{\mu\nu}^{\mathrm{(}p\mathrm{)}}})_{mod\,8})=\mathrm{\left(LL(LL({\mathit{G_{\mu\nu}^{\mathrm{(}p\mathrm{)}}}}))^{\mathit{T}}\right){}_{\mathrm{mod}\,8}}+2\mathrm{\, UR(LL(\mathit{G_{\mu\nu}^{\mathrm{(}p\mathrm{)}}})_{mod\,8})=3\mathrm{\, UR(LL(\mathit{G_{\mu\nu}^{\mathrm{(}p\mathrm{)}}})_{mod\,8})},}}\]
} \\
 offers an illuminating side to it, as can be gleaned from Table
\ref{tab:Continuation-of-interordinal} where we list arguments and
outputs to emphasize two things:\medskip{}

\begin{table}[H]
\caption{$\mathrm{Structural\: interordinality\: under\:\Lambda\: up\: to\:}p=127,p'=255$\smallskip{}
 \label{tab:Continuation-of-interordinal} }

\lyxline{\normalsize}

${\textstyle {\scriptstyle {\textstyle \Lambda:\:\mathrm{UR(LL(\mathit{G_{\mu\nu}^{\mathrm{(}p\mathrm{)}}}){}_{mod\,8}}}{\textstyle )\;\mapsto\;}}\mathrm{{\scriptstyle {\textstyle UR(UL(LL(\mathit{G_{\mu'\nu'}^{\mathrm{(}p'\mathrm{)}}})_{mod\,8}))\qquad\qquad\qquad\quad}}}}$
${\mathrm{i})\mathrm{{\textstyle \;\;{\scriptstyle argument\, row}}}\atop \mathrm{ii})\mathrm{{\textstyle {\scriptstyle \quad output\, row}\;\;}}}$\lyxline{\normalsize}\medskip{}

$\begin{array}{ccccc}
{{\scriptstyle }^{{\textstyle {\scriptstyle p=7}}}\atop _{{\scriptstyle p'=15}}} &  &  & \qquad\qquad\qquad\qquad\qquad\qquad\qquad\qquad\qquad\qquad & {^{{\scriptstyle (1)}}\atop _{{\scriptstyle (3)}}}\\
{\scriptstyle {\scriptscriptstyle }} & {\scriptstyle {\scriptscriptstyle }} & {\scriptstyle {\scriptscriptstyle }} & {\scriptstyle {\scriptscriptstyle }} & {\scriptstyle {\scriptscriptstyle }}\\
{{\textstyle {\scriptstyle }}^{{\textstyle {\scriptstyle p=15}}}\atop _{{\scriptstyle p'=31}}} &  &  &  & {^{{\scriptstyle {\scriptstyle \mathrm{{\scriptstyle sym}}(}{\scriptstyle 1\;1)}}}\atop _{{\scriptstyle \mathrm{{\scriptstyle sym}}{\scriptstyle (3\;3}{\scriptstyle )}}}}\\
{\scriptscriptstyle } & {\scriptscriptstyle } & {\scriptscriptstyle } & {\scriptscriptstyle } & {\scriptscriptstyle }\\
{{\scriptstyle }^{{\textstyle {\scriptstyle p=31}}}\atop _{_{{\scriptstyle p'=63}}}} &  &  &  & {^{{\scriptstyle {\scriptstyle \mathrm{sym(}}({\scriptstyle {\scriptstyle {5\;3\atop 3\;5}}}{\scriptstyle ),\mathrm{d_{s}}(}{\scriptstyle {1\;1\atop 1\;1}}{\scriptstyle ){\scriptstyle )}}}}\atop _{{\scriptstyle {\scriptstyle \mathrm{sym}(}{\scriptstyle (}{\scriptstyle {7\;1\atop 1\;7}}{\scriptstyle ),(}{\scriptstyle {3\;3\atop 3\;3}}{\scriptstyle )}{\scriptstyle )}}}}\\
{\scriptscriptstyle } & {\scriptscriptstyle } & {\scriptscriptstyle } & {\scriptscriptstyle } & {\scriptscriptstyle }\\
{{\textstyle {\textstyle {\scriptstyle }}^{{\textstyle {\scriptstyle p=63}}}}\atop _{{\scriptstyle p'=127}}} &  &  &  & {{\scriptstyle ^{{\scriptstyle \mathrm{sym(}({\scriptstyle {\scriptstyle {5\;3\atop 3\;5}}}),({\scriptstyle {3\;3\atop 3\;3}}),\mathrm{({\scriptstyle {\scriptstyle {5\;3\atop 3\;5}}}),d_{s}}({\scriptstyle {1\;1\atop 1\;1}}))}}}\atop _{{\scriptstyle \mathrm{sym}(({\scriptstyle {7\;1\atop 1\;7}}),({\scriptstyle {1\;1\atop 1\;1}}),({\scriptstyle {7\;1\atop 1\;7}}),({\scriptstyle {3\;3\atop 3\;3}}))}}}\\
{\scriptscriptstyle } & {\scriptscriptstyle } & {\scriptscriptstyle } & {\scriptscriptstyle } & {\scriptscriptstyle }\\
{{\scriptstyle }^{{\scriptstyle p=127}}\atop _{{\scriptstyle p'=255}}} &  &  &  & {^{{\scriptstyle \mathrm{sym(}({\scriptstyle {\scriptstyle {5\;3\atop 3\;5}}}),({\scriptstyle {3\;3\atop 3\;3}}),({\scriptstyle {7\;1\atop 1\;7}}),({\scriptstyle {3\;3\atop 3\;3}}),({\scriptstyle {\scriptstyle {5\;3\atop 3\;5}}}),({\scriptstyle {3\;3\atop 3\;3}}),\mathrm{({\scriptstyle {\scriptstyle {5\;3\atop 3\;5}}}),d_{s}}({\scriptstyle {1\;1\atop 1\;1}}))}}\atop _{{\scriptstyle \mathrm{sym}(({\scriptstyle {7\;1\atop 1\;7}}),({\scriptstyle {1\;1\atop 1\;1}}),({\scriptstyle {\scriptstyle {5\;3\atop 3\;5}}}),({\scriptstyle {1\;1\atop 1\;1}}),({\scriptstyle {7\;1\atop 1\;7}}),({\scriptstyle {1\;1\atop 1\;1}}),({\scriptstyle {7\;1\atop 1\;7}}),({\scriptstyle {3\;3\atop 3\;3}}))}}}\end{array}$\bigskip{}
 \lyxline{\normalsize} \bigskip{}

\end{table}

\noindent i) Secondary diagonal patterns, among others, are left
intact upon crossing the paraorder boundary $p$ to $p'$, as can
be seen from the tail $\mathrm{d_{s}()}$ in the respective arguments;
\smallskip{}
 \\
\noindent ii) patterns subject to the map, at least with the values
we know of, do oscillate: $\left(\begin{array}{cc}
{5\;\atop 3\;} & {3\atop 5}\end{array}\right)\leftrightarrow\left(\begin{array}{cc}
{7\;\atop 1\;} & {1\atop 7}\end{array}\right),\;\left(\begin{array}{cc}
{3\;\atop 3\;} & {3\atop 3}\end{array}\right)\leftrightarrow\left(\begin{array}{cc}
{1\;\atop 1\;} & {1\atop 1}\end{array}\right).$\\
\medskip{}

\noindent The oscillatory appearance is corroborated by the observation
that the $\mathrm{LL\mathit{f}_{mod\,8}}$ determinant (rank) alternates
between 0 and a nonzero (deficient and a complete) value among neighboring
orders $p$ and $p'$:

\begin{equation}
\det(\mathrm{LL\mathit{f}_{mod\,8}^{(3)})}=1,\;\det(\mathrm{LL\mathit{f}_{mod\,8}^{(7)})}=0,\;\det(\mathrm{LL\mathit{f}_{mod\,8}^{(15)})}=240^{2},\;\det(\mathrm{LL\mathit{f}_{mod\,8}^{(31)})}=0,\,\dots\end{equation}
\bigskip{}

\section{\label{sec:Structural-comparison-with}Structural comparison with
relation to differences}

\noindent Differences have thus far arisen at two stages in our analysis:
in Green's model, differences of \emph{squares} $\beta(p-\beta+1)$
are responsible for the capture of spin values; and, in the context
of $f$- (or $h$-) parafermions, coefficient differences constitute
in part, by fitting certain kissing numbers in the simplest case of
length-2 additive partitions, the row (column) structure of $\textrm{LL}(G_{\mu\nu}^{(p)})$.
It is therefore natural to ask which types of differences else might
be structurally constitutive, the first time so with paraorder fifteen.%
\footnote{As we have seen, at paraorders three and seven the respective partial
sequences are monomial%
}

\subsection{\label{sub:Naive-differences}Naive differences}

\noindent Let the members of the partial sequence $\:\left(G_{\rho}^{(p)}\right)\:$
be arranged in ascending order and differences springing from member
and predecessor denoted by $\Delta G_{\rho^{*}}^{(p)}$. One runs
across a peculiarity then. For $p=15,$ one gets a monotonously nondecreasing
sequence of differences,\begin{equation}
\left(\Delta G_{\rho^{*}}^{(15)}\right)=\left(2,6,6,24,72\right),\label{eq:naive-G15}\end{equation}
 whereas the related sequence for $p=31$ misses monotonicity of nondecrease:
\\
 \begin{equation}
\left(\Delta G_{\rho'^{*}}^{(31)}\right)=\left([24,72,]\,40,274,846,320,38,3186,86,10162,230,3330,40154,130590,37718,510864,1692336\right).\label{eq:naive-G31}\end{equation}
\[
\]
 Part of the order clash goes to the account of overlap of sequence
members entangled in interordinality (bracketed terms), the remaining
warps are due to intraordinal effects.\medskip{}

\noindent One way out is to economize on the number of selectable
differences, as expounded in Sect\emph{. \ref{sub:Oblique-differences}}.
Another way out is to follow the opposite track, as our demonstration
in Sect. \emph{\ref{sub:Interordinal-differences}} aims to achieve.
Different from Eqs. (\ref{eq:naive-G15})-(\ref{eq:naive-G31}) though
they may look, the types of differences thereby earmarked are clearly
linked to one another. Of particular interest turn out to be differences
derivable from Green's squares along a succession of individual carry-bit
neighbors as they shed light on the themes of Sects.\emph{ \ref{sub:The-factorization-aspect}}
to \emph{\ref{sub:The-modulo-eight-aspect}}; bar none so those derivable
from Green's squares for an \emph{enlarged} neighborhood which lead
us to believe that $b$- and $f$-parafermions may with their inherent
interordinal maps blend into a topological operator (Section \emph{\ref{sub:Kissing-numbers---}}).
But even differences with no seeming coming from Green's squares but
combining $f$- and $h$-parafermion lineage would hold topological
information, as our comments in Sects. \emph{\ref{sec:Synopsis-of-root-}}
and \emph{\ref{sec:Kissing-number-representation} }try to make clear.

\subsection{\label{sub:Oblique-differences}Skewed differences}

\noindent Recalling the way Catalan numbers were partitioned (see
Fig. \ref{fig:near-traces}), namely in form of net traces over the
secondary diagonal and adjacent diagonals in $\mathrm{LL}(G_{\mu\nu}^{(p)})$,
one would expect more meaningful differences to spring from a skewed
pairing of coefficients. There indeed exists a reduced set of $\,\frac{(p_{n-3}+1)(p_{n-3}+2)}{2}\cdot4\,$
differences, $\,\partial G_{\kappa}^{(p)}\:(p=p_{n}=15,31,\ldots)$,
that increase monotonously when subtraction is performed subsubquadrantwise
along a tilted path from the upper right to the lower left:%
\footnote{we omit configurations which stay the same upon reflection in the
secondary diagonal%
}%
\begin{figure}[H]

\caption{\label{fig:Monotonously-increasing-skewed}Monotonously increasing
skewed differences}

\begin{equation}
\begin{array}{l}
\\\mathrm{\!\neg UR(LL(G_{\mu\nu}^{\mathrm{(15)}}\mathrm{))}}:\:\qquad\qquad\qquad\qquad\qquad\qquad\\
\!\downarrow\\
\\\\\\\\\\\\\\\end{array}\begin{array}{cc}
\left.\begin{array}{ccc}
5 & - & 3\\
 &  & {{{\scriptstyle |}\atop }\atop {\textstyle 5}}\\
\\\end{array}\right. & \left.\qquad\qquad\qquad\qquad5-3=2,\right.\\
\left.\begin{array}{ccc}
11 & - & 5\\
 &  & {{{\scriptstyle |}\atop }\atop {\textstyle 17}}\\
\\\end{array}\right. & \left.\qquad\qquad\qquad\qquad\begin{array}{c}
{{\textstyle 11-5=6,}\atop }\\
{{\textstyle 17-5=12,}\atop }\end{array}\right.\\
\left.\begin{array}{ccc}
41 & - & 17\\
 &  & {{{\scriptstyle |}\atop }\atop {\textstyle 41}}\\
\\\end{array}\right. & \left.\qquad\qquad\qquad\qquad41-17=24,\right.\end{array}\label{eq:oblique-G15}\end{equation}
 \begin{equation}
\left.\begin{array}{l}
\qquad\neg\mathrm{UR}(\mathrm{LL}(G_{\mu\nu}^{\mathrm{(31)}}\mathrm{)}):\qquad\qquad\qquad\qquad\\
\qquad\downarrow\\
\\\\\\\end{array}\right.\left.\begin{array}{cc}
\left.\begin{array}{ccc}
\left.\left(\begin{array}{cc}
429 & 155\\
1275 & 429\end{array}\right)\right. & -\!\!- & \left.\!\!\left(\begin{array}{cc}
43 & 19\\
115 & 43\end{array}\right)\right.\\
 &  & |\\
\left.\right. &  & \left.\!\!\left(\begin{array}{cc}
429 & 155\\
1275 & 429\end{array}\right)\right.\end{array}\right. & \left.\begin{array}{c}
155-19=136,\\
429-43=386,\\
1275-115=1160,\end{array}\right.\end{array}\right.\label{eq:oblique-G31-head}\end{equation}
 
\end{figure}
\bigskip{}
 \begin{equation}
\begin{array}{l}
\\\left.\begin{array}{l}
\qquad\qquad\mathrm{\neg\mathrm{UR}(\mathrm{LL}\mathit{\mathrm{(}G_{\mu\nu}^{\mathrm{(31)}}\mathrm{))}}\, contd}.:\quad\qquad\\
\qquad\qquad\downarrow\\
\\\\\\\end{array}\right.\\
\left.\begin{array}{l}
\\\\\\\\\\\end{array}\right.\end{array}\begin{array}{cc}
\\\left.\begin{array}{ccc}
\left.\left(\begin{array}{cc}
4819 & 1595\\
15067 & 4819\end{array}\right)\right. & -\!\!- & \left.\!\!\left(\begin{array}{cc}
429 & 155\\
1275 & 429\end{array}\right)\right.\\
 &  & |\\
\left.\right. &  & \left.\!\!\left(\begin{array}{cc}
4905 & 1633\\
15297 & 4905\end{array}\right)\right.\end{array}\right. & \left.\begin{array}{c}
1595-155=1440,\\
1633-155=1478,\\
4819-429=4390,\\
4905-429=4476,\\
15067-1275=13792,\\
15297-1275=14022,\end{array}\right.\\
\left.\begin{array}{ccc}
\left.\left(\begin{array}{cc}
58791 & 18627\\
189371 & 58791\end{array}\right)\right. & -\!\!- & \left.\!\!\left(\begin{array}{cc}
4905 & 1633\\
15297 & 4905\end{array}\right)\right.\\
 &  & |\\
\left.\right. &  & \left.\!\!\left(\begin{array}{cc}
58791 & 18627\\
189371 & 58791\end{array}\right)\right.\end{array}\right. & \left.\begin{array}{c}
18627-1633=16994,\\
58791-4905=53886,\\
189371-15297=174074.\end{array}\right.\end{array}\label{eq:oblique-G31}\end{equation}

\subsection{\label{sub:Interordinal-differences}Interordinal differences}

\noindent An alternative to economizing on differences is to dovetail
ones non\-decreasingly into an enlarged sequence $(\Delta G_{\zeta}^{(p,p')})$,
taking from $\{G_{\rho}^{(p)}\}\cup\{G_{\rho'}^{(p')}\}$; for the
paraorder window (15,31), this procedure yields a sequence $(2,6,22,40,70,274,\dots,$
$1692336)$. The conjecture thereby dawning on the scrutator is that
in order for structural consistency with $\Delta G_{\zeta}$ to be
achieved, not so much intra- as interordinal differences of Green
squares are of importance. Those we first define with a single index
$\alpha$ by\\
 \begin{equation}
\vartheta_{\alpha}^{(p,p')}=\alpha(p'-\alpha+1)-\alpha(p-\alpha+1)=\alpha(p'-p)\;\;(\alpha=1,\dots,p):\label{eq:inter-green}\end{equation}
\begin{figure}[h]

\caption{Interordinal differences of Green's squares}
\[
\]
 \[
\left.\begin{array}{cccccc}
p'=3:\quad & 3\quad\quad & p'=7:\quad & 7 & 12 & 15\\
p=1:\quad & 1\quad\quad & p=3:\quad & 3 & 4 & 3\\
\quad & -\quad\quad & \quad & -- & -- & --\\
\quad\vartheta\;\;:\quad & 2\quad\quad & \quad\vartheta\;\;:\quad & 4 & 8 & 12\end{array}\right.\]
 \[
\]
 \[
\begin{array}{cccccccc}
p'=15:\quad & 15 & 28 & 39 & 48 & 55 & 60 & 63\\
p=\;7:\quad & 7 & 12 & 15 & 16 & 15 & 12 & 7\\
\quad & -- & -- & -- & -- & -- & -- & --\\
\quad\vartheta\;\;:\quad & 8 & 16 & 24 & 32 & 40 & 48 & 56\end{array}\]
 \[
\]
 $\begin{array}{cccccccccccccccccccccccc}
 &  &  &  &  &  &  &  & p'=31:\quad & 31 & 60 & 87 & 112 & 135 & 156 & 175 & 192 & 207 & 220 & 231 & 240 & 247 & 252 & 255\\
 &  &  &  &  &  &  &  & p=15:\quad & 15 & 28 & 39 & 48 & 55 & 60 & 63 & 64 & 63 & 60 & 55 & 48 & 39 & 28 & 15\\
\quad & \quad & \quad & \quad & \quad & \quad & \; &  & \quad & -- & -- & -- & -- & -- & -- & -- & -- & -- & -- & -- & -- & -- & -- & --\\
 &  &  &  &  &  &  &  & \quad\vartheta\quad:\quad & 16 & 32 & 48 & 64 & 80 & 96 & 112 & 128 & 144 & 160 & 176 & 192 & 208 & 224 & 240\end{array}$

\end{figure}
\bigskip{}

\noindent Compatibility is now achieved in that each $\Delta G_{\zeta}$
allows for a para\-fermial representation $\sum_{\mathfrak{q}\in\mathfrak{Q}_{\zeta}}\!\vartheta_{\mathfrak{g}_{\mathfrak{q}}}^{(\mathfrak{q},\mathfrak{q}')}\;(\mathfrak{q}'=2\mathfrak{q}+1)$.
For instance, for $\mathfrak{Q}_{\zeta}$ and $\mathfrak{g}_{\mathfrak{q}}$,
the ansatz\[
\mathfrak{Q}_{\zeta}=\left\{ 2^{\iota}-1\:\:\vert\:\:\iota\in P_{\zeta}\subset\left\{ 1,\dots,p\right\} \right\} ,\]
 \medskip{}
 \[
\mathfrak{g}_{\mathfrak{q}}\in\mathfrak{G}_{\zeta}\subset\left(G_{\rho}\right)\:|_{\mathrm{thru\,\, to\,\, order}\,\, p},\]
 \\
 can be probed to achieve representations like that of $\Delta G_{\mathrm{max}}^{(15,31)}$,\\
 \[
1692336=41\cdot2^{14}+113\cdot2^{13}+17\cdot2^{12}+11\cdot2^{11}+5\cdot2^{9}+1\cdot2^{6}+3\cdot2^{5}+1\cdot2^{4}.\]
 \\
 Whereas both for the \# of factorization types (Table \ref{tab:Factorization-types-1})
and for the spacings in prime interpolations (Fig. \ref{fig:-interpolating-prime})
linear para\-fermial expressions \begin{equation}
\mathrm{p.e.}=\sum_{i\in I_{\mathfrak{q}}}n_{i}\vartheta_{i}^{(\mathfrak{q},\mathfrak{q}')}+\!\sum_{j\in I_{\mathfrak{s}}}n_{j}\vartheta_{j}^{(\mathfrak{s},\mathfrak{s}')}\dots\qquad(\: I_{\mathfrak{q}}\subset\{1,\dots,\mathfrak{q}\},\: I_{\mathfrak{s}}\subset\{1,\dots,\mathfrak{s}\},\,\dots;\mathfrak{q}'=2\mathfrak{q}+1,\mathfrak{s}'=2\mathfrak{s}+1)\label{eq:parafermial-expression}\end{equation}
 should be as significant as for the partitions  of the \# of $G_{\rho}^{(p)}$
according to their congruence with $(7-2k)(\mathrm{mod}\:8)\:(k=0,1,2,3)$
(Table \ref{tab:Extended-partioning-1-3-5-7}).

\subsection{\label{sub:Kissing-numbers---}Kissing numbers -- the parafermion
as a topological operator}

\noindent In the introduction and in Sect. \emph{\ref{sub:Row-(column)-structure}},
it was noted that kissing numbers are tied to the row (column) structure
of $\mathrm{LL}(G_{\mu\nu}^{(p)})$. This connection is borne out
by the front members of the naive partial sequences $(\Delta G_{\rho^{*}}^{(15)})\!=\!(2,6,6,24,72)$,
$(\Delta G_{\rho'^{*}}^{(31)})\!=\!(24,72,40,{\scriptstyle \ldots})$
and the skewed partial sequence $\!(\partial G_{\kappa}^{(15)})=\!(2,6,12,24)$,
which yield the first six of them: $L_{1}=2,L_{2}\!=\!6,L_{3}\!=\!12,L_{4}\!=\!24$,
$L_{5}\!=\!40,$ $L_{6}=72$. Yet, there also exists a spin-based
connection which ties the spin-defining Green square differences --
now generalized with two indices $\alpha$ and $\beta$ to $\vartheta_{\alpha\beta}^{(p_{l},p_{u})}\equiv\beta(p_{u}-\beta+1)-\alpha(p_{l}-\alpha+1)$
where $p_{l}$ and $p_{u}$ are Mersenne numbers defined by\emph{
$p_{l}\in\left\{ 1,3,7,\ldots,2^{l}-1\,|\,2^{l}\leq2D\right\} $},
$p_{u}\in\left\{ 3,7,15,\ldots2^{u}-1\,|\,2^{u}\leq32D\right\} $
$(l,u,D\in\mathbb{N}$; $p_{l}<p_{u})$ and $\alpha$ and $\beta$
are respectively running from 1 to $p_{l}$ and $p_{u}$ -- to kissing
numbers:%
\footnote{As the bulk of kissing numbers shown in Table \ref{tab:The-first-sixteen}
are not certified, we are on risky ground here.%
}%
\begin{table}
\caption{{\small \label{tab:The-first-sixteen}The first sixteen kissing numbers
of Euclidean $D$-space as represented by $\vartheta_{\alpha\beta}^{(p_{l},p_{u})}$,
including spinless quantities $\vartheta_{\alpha\beta}^{(p_{l},p_{u})}=0$ }}
\medskip{}

\noindent \setlength{\fboxsep}{1mm}\lyxline{\normalsize}

\begin{tabular}{c|ccccc}
$\:\quad p_{l}\:\quad$ & $p_{u}$ &  & $\vartheta_{\alpha\alpha}^{(p_{l},p_{u})}\enskip(\alpha\;\textrm{odd)}\equiv\vartheta_{\lambda}^{(p_{l},p_{u})}\:^{*}$ & $\qquad\vartheta_{\alpha\alpha}^{(p_{l},p_{u})}\enskip(\alpha\;\textrm{even)}$ & $\qquad\qquad\qquad\qquad\vartheta_{\alpha\beta}^{(p_{l},p_{u})}\quad(\alpha\neq\beta)$\tabularnewline
\end{tabular}

\lyxline{\normalsize}\begin{tabular}{c|cccccc}
$\:1$ & $\:3$ & $\quad$ & ${\fbox{\ensuremath{{\textstyle \vartheta_{1}^{(1,3)}=\mathit{L}_{1}=2}}}\atop }$ & $\quad$ & $\qquad\quad\textrm{-}\quad\qquad$ & -\tabularnewline
\hline 
$\:1$ & $\:7$ &  & ${\atop {\fbox{\ensuremath{{\textstyle \vartheta_{1}^{(1,7)}=\mathit{L}_{2}=6}}}\atop }}$ &  & - & -\tabularnewline
$\:3$ & $\:7$ &  & ${\fbox{\ensuremath{{\textstyle \vartheta_{3}^{(3,7)}=\mathit{L}_{3}=12}}}\atop }$ &  & - & -\tabularnewline
\hline 
$\:1$ & $15$ &  & - &  & - & -\tabularnewline
$\:3$ & $15$ &  & - &  & $L_{4}=24$ & $L_{3}=12$\tabularnewline
$\:7$ & $15$ &  & ${\fbox{\ensuremath{{\textstyle \vartheta_{3}^{(7,15)}=\mathit{L}_{4}=24,}}}\atop {{\textstyle \fbox{\ensuremath{{\textstyle \vartheta_{5}^{(7,15)}=\mathit{L}_{5}=40}}}}\atop }}$ &  & - & $0,L_{3}=12$\tabularnewline
\hline 
$1,3$ & $31$ &  & - &  & - & -\tabularnewline
$\:7$ & $31$ &  & ${{\textstyle \vartheta_{1,1}^{(7,31)}=L_{4}=24,}\atop {{\textstyle \fbox{\ensuremath{{\textstyle \vartheta_{3}^{(7,31)}=\mathit{L}_{6}=72}}}}\atop }}$ &  & - & ${\textstyle \fbox{\ensuremath{{\textstyle \mathit{L}_{8}=240}}}}$\tabularnewline
$15$ & $31$ &  & - &  & - & $0,L_{3}=12,L_{4}=24,L_{6}=72,L_{8}=240$\tabularnewline
\hline 
$1,3$ & $63$ &  & - &  & - & -\tabularnewline
$\:7$ & $63$ &  & - &  & $L_{10}=336$ & ${\textstyle \fbox{\ensuremath{{\textstyle \mathit{L}_{12}=756}}}}$\tabularnewline
$15$ & $63$ &  & ${{\textstyle \vartheta_{5,5}^{(15,63)}=L_{8}=240,}\atop {\textstyle {\textstyle \fbox{\ensuremath{{\textstyle \vartheta_{7}^{(15,63)}=\mathit{L}_{10}=336}}}}}}$ &  & - & $0,L_{4}=24$\tabularnewline
$31$ & $63$ &  & - &  & - & $0,L_{3}=12,L_{5}=40,L_{8}=240,L_{10}=336$\tabularnewline
\hline 
$\:1$ & $127$ &  & ${{\textstyle }\atop {\textstyle {\textstyle \fbox{\ensuremath{{\textstyle \vartheta_{1}^{(1,127)}=\mathit{L}_{7}=126}}}}}}$ &  & - & -\tabularnewline
$\:3$ & $127$ &  & - &  & - & -\tabularnewline
$\:7$ & $127$ &  & - &  & $L_{8}=240$ & -\tabularnewline
$15$ & $127$ &  & $\vartheta_{3,3}^{(15,127)}=L_{10}=336$ &  & - & $L_{6}=72$\tabularnewline
$31$ & $127$ &  & - &  & - & ${{\textstyle 0,L_{3}=12,L_{5}=40,L_{8}=240},{\textstyle \fbox{\ensuremath{{\textstyle \mathit{L}_{15}=2340}}}}\atop }$\tabularnewline
$63$ & $127$ &  & - &  & - & ${{\textstyle 0,L_{3}=12,L_{6}=72},{\textstyle \fbox{\ensuremath{{\textstyle \mathit{L}_{9}=272}}}}{\textstyle ,L_{10}=336}\atop }$\tabularnewline
\hline 
$\:1$ & $255$ &  & - &  & - & -\tabularnewline
$\:3$ & $255$ &  & $\vartheta_{3}^{(3,255)}=L_{12}=756$ &  & - & -\tabularnewline
$\:7$ & $255$ &  & - &  & - & $L_{8}=240$\tabularnewline
$15$ & $255$ & \multicolumn{1}{l}{} & - &  & - & -\tabularnewline
$31$ & $255$ &  & - &  & - & $0,L_{4}=24,L_{12}=756$\tabularnewline
$63$ & $255$ &  & - &  & - & $0,L_{4}=24,L_{6}=72,L_{8}=240,{\textstyle \fbox{\ensuremath{{\textstyle \mathit{L}_{16}=4320}}}}$\tabularnewline
$127$ & $255$ &  & - &  & - & $0,L_{3}=12,L_{4}=24,L_{10}=336,L_{12}=756,L_{16}=4320$\tabularnewline
\hline 
$1,{\scriptstyle \ldots},15$ & $511$ &  & - &  & - & -\tabularnewline
$31$ & $511$ &  & $\vartheta_{9}^{(31,511)}=L_{8}=4320$ &  & - & $L_{10}=336$\tabularnewline
$63$ & $511$ &  & - &  & - & $0,L_{3}=12,L_{16}=4320$\tabularnewline
$127$ & $511$ &  & - &  & - & $0,L_{3}=12,L_{8}=240,L_{10}=336,L_{16}=4320$\tabularnewline
$255$ & $511$ &  & - &  & - & $0,L_{3}=12,L_{4}=24,L_{9}=272,{\textstyle L_{10}=336,L_{12}=756,L_{15}=2340}$\tabularnewline
\end{tabular}\lyxline{\normalsize}

$^{*})$ first appearance only
\end{table}

\noindent The Green squares $\alpha(p_{l}-\alpha+1)$ and $\beta(p_{u}-\beta+1)$
define points on two sets of parabolas -- with extremal ordinates
$2^{2l-2}$ and $2^{2u-2}$, respectively, and general ordinate differences
divided by two defining (generalized) spin. Vanishing ordinate differences
define a spin-0 subclass which is given by $\alpha=(2^{l}-1)2^{m},$
$\beta=2^{m},$ $p_{l}=2^{l+m+1}-1$, $p_{u}=2^{2l+m}-1$ $(l=2,3,4,\ldots;m=0,1,2,\ldots)$.
With no larger spin than of two, Green squares can be subdivided into
two classes -- a degenerate one, with remainder 1 under $\alpha(p_{l}-\alpha+1)\textrm{ mod }4$
for $p_{l}=\alpha=1$, and a regular one with remainders 0 and 3 under
both $\alpha(p_{l}-\alpha+1)\textrm{ mod }4\:(l>1)$ and $\beta(p_{u}-\beta+1)\textrm{ mod }4$.
It follows that every kissing number with factorization $L_{x}=2^{z_{1}}\cdot\pi_{r}^{z_{r}}\cdot\ldots\:(z_{1}>1,\pi_{r}>2)$
possesses at least one representation $\vartheta_{\alpha\beta}^{(p_{l},p_{u})}$,
while a kissing number with factorization $L_{y}=2\cdot\pi_{r}^{z_{r}}\cdot\ldots\:(\pi_{r}>2)$
by necessity changes the remainders $0,3$ under $[\beta(p_{u}-\beta+1)-L_{y}]\textrm{ mod }4$
to $1,2$ and requires the degenerate class $\alpha(p_{l}-\alpha+1)\textrm{ mod }4\:(p_{l}=\alpha=1)$
to counterbalance that change. Out of the first sixteen kissing numbers
containing base prime two non-exponentiated -- $L_{1}=2,L_{2}=6$,
$L_{7}=126$, $L_{11}=438(582)$, $L_{13}=918(1154)$, $L_{14}=1422(1606)$
\footnote{{\tt http://www.math.rwth-aachen.de/}$\!\thicksim${\tt Gabriele.Nebe/LATTICES/kiss.html};
the numbers in parentheses spring from nonlattice calculations%
} -- only three conform with this condition via the representation
$\vartheta_{1,1}^{(p_{l},p_{u})}$. In other words, they are equivalent
to the Mersenne number differences $L_{1}=3-1,L_{2}=7-1,L_{7}=127-1$.
Together with the kissing numbers $L_{x}=2^{z_{1}}\cdot\pi_{r}^{z_{r}}\cdot\ldots\:(z_{1}>1,\pi_{r}>2)$
and zeroing differences, they are shown in Table \ref{tab:The-first-sixteen}.
A surprising feature can be read from that table: If only a kissing
number's first appearance as $\vartheta_{\alpha\alpha}^{(p_{l},p_{u})}\:(\alpha\textrm{ odd})$
or $\vartheta_{\alpha\beta}^{(p_{l},p_{u})}\:(\alpha\neq\beta)$ is
considered (framed items), first appearances as $\vartheta_{\alpha\alpha}^{(p_{l},p_{u})}\:(\alpha\textrm{ odd})$
lead to an $\alpha$ sequence $1,1,3,3,5,3,7,1$. Taken pairwise,
the $\alpha$'s follow the alternation map $\Lambda=(\mathrm{mod}\:8)\circ(\times3)$
characteristic of $\mathrm{LL}(G_{\mu\nu}^{(p)})$ structure -- which
is the reason for our usage of the single index $\lambda$ in the
table. We arrive at the unlooked-for topologic\medskip{}

\begin{conj}
\label{con:A-hypersphere-configuration}$f^{(p)}$- (or $h^{(p)}$-)
parafermions are containers of hypersphere configurations of densest
packing, in Euclidean \emph{$D^{(p)}$-space} down to those in \emph{$D^{(1)}$}-space,
where \emph{$D^{(p)}$} is the largest dimension for which the kissing
number $L_{D^{(p)}}$ determines the row (column) structure of $\mathrm{LL}(G_{\mu\nu}^{(p)})$
(or $\mathrm{LL}(J_{\mu\nu}^{(p)})$) -- \emph{$D^{(15)}=7$}, for
instance.%
\footnote{As Conjecture \ref{con:In--(or} suggests, \emph{$D^{(31)}=31$.}%
} Dual to this inner structural connection is the spin-based connection
relating densenst-packing hypersphere configurations to Green square
differences of parafermions of Mersennian order such that either $L_{D}=\vartheta_{\alpha\beta}^{(p_{l},p_{u})}$
or $L_{D}=\vartheta_{\lambda_{1}}^{(p_{l_{1}},p_{u_{1}})}\pm\vartheta_{\lambda_{2}}^{(p_{l_{2}},p_{u_{2}})}\pm\ldots$
. This latter connection is exterior in the sense that $p_{u}$ may
become larger than $p$. 
\end{conj}
\noindent The exterior connection can be delineated as follows:
\begin{cor}
\label{cor:If-a-kissing}For each pair $p_{l},p_{u}$, interordinal
differences of Green squares $\beta(p_{u}-\beta+1)-\alpha(p_{l}-\alpha+1)$
form a distinct rectangular matrix $\vartheta_{\alpha\beta}^{(p_{l},p_{u})}$.
Among its entries, which include zeroing differences leading to spin
0

, of particularly interest are those coinciding with members of the
class $L_{x}=2^{z_{1}}\cdot\pi_{r}^{z_{r}}\cdot\ldots\:(z_{1}>1,\pi_{r}>2)$
of kissing numbers. If only $\vartheta_{\alpha\alpha}^{(p_{l},p_{u})}\:(\alpha\textrm{ odd})$
and $\vartheta_{\alpha\beta}^{(p_{l},p_{u})}\:(\alpha\neq\beta)$
are considered, the first appearances of these kissing number in representations
$\vartheta_{\alpha\alpha}^{(p_{l},p_{u})}\:(\alpha\textrm{ odd})\equiv\vartheta_{\lambda}^{(p_{l},p_{u})}$
are characterized by indices $\lambda$ that (pairwise, in ascending
order) are in one-to-one correspondence with the patterns\emph{ }$\mathrm{({\scriptstyle {1,1\atop 1,1}})},$
$\Lambda\left(({\scriptstyle {1,1\atop 1,1}})\right)$\emph{, }$({\scriptstyle {5,3\atop 3,5}})$,
$\Lambda\left(({\scriptstyle {5,3\atop 3,5}})\right)$ where\emph{
}$\Lambda=(\mathrm{mod}\:8)\circ(\times3)$. 
\end{cor}
\noindent One approach to narrowing the range of pairs $p_{l},p_{u}$
leans on paraorder sums which are endowed with the identity\[
\Sigma_{i=1}^{2n-1}p_{i}-\frac{1}{C_{p_{n}}\mathrm{B}(p_{n},p_{n}+1)}=\Sigma_{i=1}^{2n-1}p_{i}-p_{n}(p_{n}+1)=-n+\Sigma_{i=1}^{n-1}p_{i}\]
\noindent (where B(,) is the beta function) and allow taking three
different paraorders into account on each assignment:\[
\begin{array}{cc}
p_{l}:=p_{n-1}, & p_{u}:=p_{n},\\
p_{l}:=p_{n-1}, & p_{u}:=p_{2n-1},\\
p_{l}:=p_{n}, & p_{u}:=p_{2n-1}.\end{array}\]
\noindent  It is easily shown that these choices include $L_{x}=\vartheta_{\alpha\alpha}^{(p_{l},p_{u})}\:(\alpha\textrm{ even})$
and thus are a broader approach than the $\Lambda$-approach.

\section{\label{sec:Synopsis-of-root-}Synopsis of root-{\em f}- and root-{\em
h} related coefficient differences}

\noindent Making the review more complete by a further sideglance
to the root-$h$ sequence is overdue. In the introduction it was already
stated that this sequence bears a resemblance to the root-$f$ sequence.
The kinship tellingly expresses itself in the relations, starting
with $p=7$, $\: p'=15$,

\begin{equation}
\mathrm{\mathrm{UR(UL(LL\mathit{h^{\mathrm{(}p'\mathrm{)}}\mathrm{))}}}=LL(LL\mathit{h}^{(\mathit{p})})}-2\,\mathrm{UR(LL\mathit{h^{\mathrm{(}p\mathrm{)}}\mathrm{)}},}\label{eq:inter-h}\end{equation}
 \begin{equation}
\mathrm{\mathrm{UR(LL(LL\mathit{h^{\mathrm{(}p'\mathrm{)}}\mathrm{))}}}=\mathrm{LL(UL(LL\mathit{h^{\mathrm{(}p'\mathrm{)}}\mathrm{))}}}-2\, UR(UL(LL\mathit{h}^{(\mathit{p'})}))}+2\,\mathrm{UR(LL\mathit{h^{(p)}\mathrm{)}},}\label{eq:inter-intra-h}\end{equation}
 only the first of which is purely interordinal, while the second
is a mixture of intra- and interordinal relationship. Juxtaposing
these opposite (\ref{eq:inter})-(\ref{eq:intra}) -- identities that
we remember are pure interordinal and intraordinal respectively --,
one is not surprised to find that the partial sequences $\left(J_{\omega}^{(p)}\right)$
-- with $J_{\omega}^{(p)}$ as representatives of $J_{\mu\nu}^{(p)}$
-- cease being monomial already at paraorder seven.%
\footnote{As opposed to the partial sequences $\left(G_{\rho}^{(p)}\right)$
which do not move on from monomiality until paraorder fifteen%
} Starting out with that order, partial sequences with differences
$\Delta J_{\omega^{*}}^{(p)}$, $\partial J_{\kappa}^{(p)}$%
\footnote{where we again encounter a reduced set of $\,\frac{(q+1)(q+2)}{2}\cdot4\,$
differences $\,\partial J_{\kappa}^{(p)}\:(p=15,31,\ldots;q=1,3,\ldots)$,
based on subsubquadrantwise subtraction performed along a tilted path
that pairs distinct $J_{\mu,\nu}^{(p)}$ from upper right to lower
left; even though performed in the same way, the subtraction process
does not automatically lead to a monotonously increasing sequence
of differences such as $\left(\partial G_{\kappa}^{(p)}\right)$%
} and $\,\Delta J_{\theta}^{(p,p')}$ then are readily formable. Briefly
expounding to what extent a synopsis between them and their $G$ counterparts
on the one hand and $\vartheta_{\lambda}^{(p_{l},p_{u})}$ on the
other can be used to the kissing-number problem is the subject of
this section. 

\noindent We have already learned three ways of expressing kissing
numbers: 

\noindent a) in the Introduction and in Sect. \emph{\ref{sub:Row-(column)-structure}},
by additive partitions within rows, or groups of rows, of $\textrm{LL(}\mathit{G^{\mathrm{(}p\mathrm{)}}\mathrm{)}}$;

\noindent b) in Sect\emph{. \ref{sub:Kissing-numbers---}}, by interordinal
differences $\vartheta_{\alpha\beta}^{(p_{l},p_{u})}$ of Green squares;

\noindent c) same place, by higher-order parafermial differences
$\vartheta_{\lambda_{1}}^{(p_{l_{1}},p_{u_{1}})}\pm\vartheta_{\lambda_{2}}^{(p_{l_{2}},p_{u_{2}})}\pm\ldots$
. 

\noindent It was also remarked upon the connection of the naive partial
sequences of Eq. (\ref{eq:naive-G31})\[
\begin{array}{cc}
(\Delta G_{\rho'^{*}}^{(15)})=(2,6,6,24,72),\\
(\Delta G_{\rho''^{*}}^{(31)})=(24,72,40,274,846,320,\ldots)\end{array}\]
 as well as the skewed sequence of scheme (\ref{eq:oblique-G15})
\[
(\partial G_{\kappa'}^{(15)})=(2,6,12,24)\]
with kissing numbers. What remained to be checked is whether $G_{\mu\nu}^{(\tilde{p})}$-
and $J_{\mu\nu}^{(\tilde{p})}$-derived differences, $\tilde{p}\in\left\{ p,p'\right\} $,
have a way of jointly determining these numbers. We therefore computed
certain $J_{\mu\nu}^{(p)}$-derived partial sequences of differences
for the occasion: 

\noindent the naive\begin{equation}
\begin{array}{c}
\left(\Delta J_{\omega^{*}}^{(7)}\right)=(4),\qquad\left(\Delta J_{\omega'^{*}}^{(15)}\right)=(38,\!6,\!14,\!134),\\
\left(\Delta J_{\omega''^{*}}^{(31)}\right)=(688974,53888,4474,388,54,104,26,1176,24,204,14000,2722,176724,28580,2662662),\end{array}\end{equation}
 \\
 and the skewed\begin{equation}
\begin{array}{c}
\left(\partial J_{\kappa'}^{(15)}{}\right)=(-6,20,-58),\\
\left(\partial J_{\kappa''}^{(31)}{}\right)=(104,-388,1404,1226,1202,-4474,-4394,16722,16442,14228,-53968,205584).\end{array}\end{equation}
 \\
With the root-$f$ related sequence $\!\left(\partial G_{\kappa''}^{(31)}\right)$
of scheme (\ref{eq:oblique-G31-head})-(\ref{eq:oblique-G31}) computed
to

\begin{equation}
\left(\partial G_{\kappa''}^{(31)}\right)=(136,386,1160,1440,1478,4390,4476,13792,14022,16994,53886,174074),\end{equation}
 \\
 we've actually found a scheme construing the values $L_{D}$ $(D\leq8)$
as second-order synoptic differences: 

\noindent those linked to odd-dimensional Euclidean spaces in representations
that mix $\Delta$- and $\partial$ terms,

\[
\!\!\!\begin{array}{ccl}
L_{1} & = & 2=\Delta J_{3}^{(15)}-\partial G_{3}^{(15)}=14-12,\\
L_{3} & = & 12=\Delta G_{7}^{(31)}-\partial J_{7}^{(31)}=38-26,\\
L_{5} & = & 40=\partial J_{1}^{(15)}+\Delta G_{1}^{(15)}=38+2,\\
L_{7} & = & 126=\Delta G_{5}^{(15)}+\partial J_{5}^{(31)}=72+54,\\
{\scriptscriptstyle } & {\scriptscriptstyle } & {\scriptscriptstyle }\end{array}\]
\noindent and those linked to even-dimensional spaces in representations
that are homogeneous in either $\Delta$ or $\partial$:

\[
\begin{array}{ccl}
L_{2} & = & 6=\Delta J_{1}^{(7)}+\Delta G_{1}^{(15)}=4+2,\\
L_{4} & = & 24=\Delta G_{7}^{(31)}-\Delta J_{3}^{(15)}=38-14,\\
L_{6} & = & 72=\Delta G_{9}^{(31)}-\Delta J_{3}^{(15)}=86-14,\\
L_{8} & = & 240=\partial G_{1}^{(31)}+\partial J_{1}^{(31)}=136+104.\\
{\scriptscriptstyle } & {\scriptscriptstyle } & {\scriptscriptstyle }\end{array}\]
 \noindent Without knowledge of the inputs $\left(\Delta G_{\rho'''^{*}}^{(63)}\right),\left(\partial G_{\kappa'''}^{(63)}\right)$
on the one hand and $\left(\Delta J_{\omega'''^{*}}^{(63)}\right),\left(\partial J_{\kappa'''}^{(63)}\right)$
on the other, one cannot be sure of how to enlarge that picture for
$D>8$. Fortunately some key information is still within reach:\bigskip{}
\begin{table}[H]

\caption{{\small \label{tab:Kissing-numbers-from}Kissing numbers $L_{9}$
to $L_{16}$ as 2nd-order synoptic\,/\,returning interordinal (or
higher-order synoptic\,/\,otherwise interordinal) differences}}

\lyxline{\normalsize}

\medskip{}

\hfill{}\begin{tabular}{cccccccc}
$L_{9}$ & $=$ & $272$ & $=$ & $\Delta G_{4}^{(31)}-\Delta G_{1}^{(15)}$ & $=$ & $274-2$ & \tabularnewline
$L_{10}$ & $=$ & $336$ & $=$ & $\vartheta_{7,1}^{31,511}$ &  &  & \tabularnewline
$L_{11}$ & $=$ & $438$ & $=$ & $\textrm{\ensuremath{{\scriptstyle \textrm{h.o.i./h.o.s}.}}}$ &  &  & \tabularnewline
$L_{12}$ & $=$ & $756$ & $=$ & $\textrm{\ensuremath{{\scriptstyle \textrm{o.i./h.o.s}.}}}$ &  &  & \tabularnewline
$L_{13}$ & $=$ & $918$ & $=$ & $\Delta G_{5}^{(31)}+\Delta G_{5}^{(15)}$ & $=$ & $846+72$ & \tabularnewline
$L_{14}$ & $=$ & $1422$ & $=$ & $\textrm{\ensuremath{{\scriptstyle \textrm{h.o.i./h.o.s.}}}}$ &  &  & \tabularnewline
$L_{15}$ & $=$ & $2340$ & $=$ & $\textrm{\ensuremath{{\scriptstyle \textrm{o.i./h.o.s.}}}}$ &  &  & \tabularnewline
$L_{16}$ & $=$ & $4320$ & $=$ & $\vartheta_{9,9}^{(31,511)}$ &  &  & \tabularnewline
\end{tabular}\hfill{}\medskip{}
 \lyxline{\normalsize}

\noindent \medskip{}
 
\end{table}

\noindent Now the alternation map $\Lambda$ is closely related to
the odd-integer partitions of the number 8: the quadripartite 1+1+3+3=8
and the bipartite 5+3=8 and 7+1=8. Thus the action of $\Lambda$ can
be put in one-to-one correspondence with either the alternation of
the halves of the quadripartite or the alternation of the full bipartite
partition(s). The alternation of the halves of the quadripartite partition
fits in one eight-period of dimensions and in fact is in one-to-one
correspondence with the action of $\Lambda$ on characteristic increments
in that eight-period; conversely, the alternation of the full bipartite
partitions should fit in two such periods and also be in one-to-correspondence
with characteristic index increments in there. Writing\begin{equation}
\begin{array}{rcl}
i & = & 2^{m}+\lambda_{1},\\
j & = & 2^{m}+\lambda_{1}+\lambda_{2},\qquad\qquad(\textrm{for some }m\geq3)\\
k & = & 2^{m}+\lambda_{1}+\lambda_{2}+\lambda_{3},\\
l & = & 2^{m}+\lambda_{1}+\lambda_{2}+\lambda_{3}+\lambda_{4,}\end{array}\label{eq:allocation}\end{equation}
 we arrive at the following conjecture which supplements Corollary
\ref{cor:If-a-kissing}:\\

\begin{conj}
\label{con:The-kissing-partite}The kissing number associated with
a hypersphere configuration of densest packing in Euclidean $D$-space
is representable both by a 2nd-order synoptic difference and an interordinal
difference\emph{ $\vartheta_{\alpha\beta}^{(p_{l},p_{u})}$ of Green
squares for $D\leq8$. }As\emph{ $D>8$,} representations can be assigned
accordingly and either are pairwise 2nd-order synoptic based on paraorders
$2^{m}-1$, $2^{m+1}-1$ and associated with dimensions\emph{ $D_{1}=i$,
$D_{3}=k$, }or returning interordinal \emph{$\vartheta_{\alpha_{j},\beta_{j}}^{(p_{l},p_{u})}$,
$\vartheta_{\alpha_{l},\beta_{l}}^{(p_{l},p_{u})}$} based on paraorders
$2^{m+1}-1$, $2^{m+5}-1$ and associated with dimensions \emph{$D_{2}=j$,
$D_{4}=l$, }while consisting of higher-order synoptic\,/\,otherwise
interordinal differences at interstitial dimensions. \emph{$i,j,k$
}and\emph{ $l$} are determined by the above system of equations,
and the span of dimensions taken is one full\emph{ 8}-period for $(\lambda_{1},\lambda_{2})\cong\mathrm{({\scriptstyle {1,1\atop 1,1}})}$,
$(\lambda_{3},\lambda_{4})\cong\mathrm{\Lambda\left(({\scriptstyle {1,1\atop 1,1}})\right)}$
and two successive \emph{8}-periods for\emph{$(\lambda_{1},\lambda_{2})\cong({\scriptstyle {5,3\atop 3,5}})$,
$(\lambda_{3},\lambda_{4})\cong\Lambda\left(({\scriptstyle {5,3\atop 3,5}})\right)$.}
\end{conj}

\section{\label{sec:Kissing-number-representation}Kissing number representation
and dimensional periodicities\,}

\noindent A natural question to ask is if Conjecture \ref{con:The-kissing-partite}
allows instantiations of\emph{ $\lambda_{1}$ etc. }to repeat periodically.
If $L_{2^{s}+b_{0}}$ $((2^{s}+b_{0})\textrm{mod}\,8=\lambda_{1})$
is representable 2nd-order synoptic, so too could $L_{2^{s+ct}+b_{t}}$
$((2^{s+ct}+b_{t})\textrm{mod}\,8=\lambda_{1},c=\mathrm{const.})$.
While the input/output entries of Table \ref{tab:Continuation-of-interordinal}
and the $\lambda$-sequence of Table \ref{tab:The-first-sixteen},
which share with\emph{ } \emph{$\lambda_{1},\lambda_{2},$$\lambda_{3},\lambda_{4}$
}the $\Lambda$-mapping precept, signal nothing of the kind, the stenoscopy
of kissing numbers described in Sect.\emph{ \ref{sub:Row-(column)-structure}}
hints at such a possibility. The argument proceeds as follows. Via
the stenoscopic coupling\bigskip{}
\begin{equation}
D^{(p)}=\max(D)\mid L_{D}\leq\Sigma_{i=1}^{(q+1)/2}G_{q+1,i}^{((p-1)/2)}\label{eq:max-L_D}\end{equation}
 \bigskip{}

\noindent (dotted underlined in Table \ref{tab:Stenoscopy-of-kissing}),
a least dimension falling within the subsequent interordinal corridor
is being defined:\bigskip{}
\begin{equation}
D_{\textrm{lowest}}^{(p')}=D^{(p)}+1.\label{eq:lowest-next D}\end{equation}
\bigskip{}
 For reasons given below, we relate to the latter the dimension $N$
taking values\bigskip{}
\[
N={\scriptscriptstyle -}\left\lceil {\scriptscriptstyle {\textstyle (n+1)/{\scriptstyle 2}}}\right\rceil {\scriptscriptstyle +}\Sigma_{i=1}^{n-1}p_{i}=\left\lfloor \log_{2}C_{q'}\right\rfloor \]
\bigskip{}

\noindent (an identity that first sprang up in Eq. (\ref{eq:log C identity})).
We have assembled a selection of the numbers $D_{\textrm{lowest}}^{(p')}$
and $N$ in Table \ref{tab:The-lowest-dimensions}:

\bigskip{}
\begin{table}[H]
\caption{\label{tab:The-lowest-dimensions}least dimension $D_{\textrm{lowest}}^{(p')}$
covered by $f$- ($h$-) parafermion of order $p'\equiv p_{n+1}$,
and companion dimension $N$ }

\lyxline{\normalsize}\begin{tabular}{lccccccccccccccc}
$p'$ & $2^{4}\textrm{-1}$ & $2^{5}\textrm{-1}$ & $2^{6}\textrm{-1}$ & $2^{7}\textrm{-1}$ & $2^{8}\textrm{-1}$ & $2^{9}\textrm{-1}$ & $2^{10}\textrm{-1}$ & $\textrm{\ensuremath{2^{11}}-1}$ & $2^{12}\textrm{-1}$ & $2^{13}\textrm{-1}$ & $2^{14}\textrm{-1}$ & $2^{15}\textrm{-1}$ & $2^{16}\textrm{-1}$ & $2^{17}\textrm{-1}$ & $\cdots$\tabularnewline
\cline{1-15} 
$\!\!\! D_{\textrm{l'st}}^{(p')}$ & 2 & 8 & 32 & 112 & 416 & 1\,640 & 6\,586 & 25\,504 & 101\,132 & 407\,154 & 1\,642\,292 & 6\,618\,374 & 26\,638\,982 & 107\,107\,722 & $\cdots$\tabularnewline
$N$ & 2 & 8 & 23 & 53 & even & even & 497 & 1007 & even & even & 8171 & 16361 & even & even & $\cdots$\tabularnewline
\end{tabular}\lyxline{\normalsize}\bigskip{}

\end{table}
\bigskip{}
 \noindent As far as that Table \ref{tab:The-lowest-dimensions}
goes, $D_{\textrm{lowest}}^{(p')}$ runs from $2$ to 107\,107\,722,
and $N$ from 2 to the two even numbers following 16361. From Table
\ref{tab:Periodicity-of-binary} we see two possibilities for follow-up.
If only $D_{\textrm{lowest}}^{(p')}-N\;(N\:\textrm{odd-numbered})$
is being realized, Table \ref{tab:Periodicity-of-binary} suggests
the possibility of second-order synoptic (or returning interordinal)
kissing number representability for\begin{equation}
L_{2^{s+16t}+b_{t}}\quad((2^{s+16t}+b_{t})\textrm{mod}\,8=\lambda;\; t=0,1,2,\ldots),\label{eq:simple period}\end{equation}

\noindent with $c$ coinciding with a key figure of Table \ref{tab:The-first-sixteen},
$\frac{\max(p_{u}+1)}{\max(p_{l}+1)}=16$. If, on the other hand,
$D_{\textrm{lowest}}^{(p')}-N'\;(N'\:\textrm{too odd-numbered})$
is being realized as well, then representations with multiple periodicities
might come to light. Either way, if so, relating $D_{\textrm{lowest}}^{(p')}$
to $N$ (and $N'$) would lie at the heart of the effectiveness of
the interordinal map $\Lambda$ for the kissing number representations
in question.%
\begin{table}[H]
\caption{\label{tab:Periodicity-of-binary}Possibility of (double) periodicity
in second-order synoptic and/or returning interordinal kissing number
representation}

\lyxline{\normalsize}\medskip{}

$p'\equiv p_{n+1}$$\qquad\qquad\qquad\qquad\qquad\qquad\qquad\qquad\qquad\qquad\qquad\qquad$${{\displaystyle D_{\textrm{lowest}}^{(p')}-N}\atop {\displaystyle D_{\textrm{lowest}}^{(p')}-N'}}\qquad(N,N'\:\textrm{odd-numbered})$
\medskip{}

\lyxline{\normalsize}\bigskip{}
$\begin{array}{crc}
2^{7}-1\qquad\qquad\qquad\qquad\qquad & {{\textstyle }\atop {\textstyle }} & {{\textstyle 112-23=2^{5}+57,\quad(2^{5}+57)\textrm{mod}\,8=1}\atop {\textstyle 112-53=2^{5}+27,\quad(2^{5}+27)\textrm{mod}\,8=3}}\\
\\2^{11}-1\qquad\qquad\qquad\qquad\qquad & {{\textstyle }\atop {\textstyle }} & {{\textstyle 25\,504-497=2^{9}+24\,495,\quad(2^{9}+24\,495)\textrm{mod}\,8=7}\atop {\textstyle 25\,504-1\,007=2^{9}+23\,985,\quad(2^{9}+23\,985)\textrm{mod}\,8=1}}\\
\\2^{15}-1\qquad\qquad\qquad\qquad\qquad & {{\textstyle }\atop } & {{\textstyle 6\,618\,374-8\,171=2^{13}+6\,602\,011,\quad(2^{13}+6\,602\,011)\textrm{mod}\,8=3}\atop {\textstyle 6\,618\,374-16\,361=2^{13}+6\,593\,821,\quad(2^{13}+6\,593\,821)\textrm{mod}\,8=5}}\\
\\2^{19}-1\qquad\qquad\qquad\qquad\qquad & {{\textstyle }\atop {\textstyle }} & {{\textstyle 598\,753\,098-131\,045=2^{17}+598\,490\,981,\quad(2^{17}+598\,490\,981)\textrm{mod}\,8=5}\atop {\textstyle 598\,753\,098-262\,115=2^{17}+598\,359\,911,\quad(2^{17}+598\,359\,911)\,\textrm{mod}\,8=7}}\\
\\2^{23}-1\qquad\qquad\qquad\qquad\qquad & {{\textstyle }\atop {\textstyle }} & {{\textstyle 54\,868\,958\,480-2\,097\,119=2^{21}+54\,864\,764\,209,\quad(2^{21}+54\,864\,764\,209)\textrm{mod}\,8=1}\atop {\textstyle 54\,868\,958\,480-4\,194\,269=2^{21}+54\,862\,667\,059,\quad(2^{21}+54\,862\,667\,059)\textrm{mod}\,8=3}}\end{array}$\medskip{}
 \lyxline{\normalsize}\bigskip{}

\end{table}
\noindent To understand why, we have to recall the definition of
the suffix of consecutive prime factors, $\mathrm{\: SCPF}(G_{q+2,1}^{(p)})$
from Sect. \emph{\ref{sub:The-factorization-aspect}}. Obviously,
our $\lambda$'s \,3, 5 and 7 and Mersenne primes > 7 represent the
only base primes that are not SCPF primes in the factorization of
$C_{q}\:(q=7,15,31,\ldots)$. When all base primes defining the infix
(Mersenne primes > 7) are mapped to the remaining value $\lambda=1$,
the factorization has a prefix that can be defined by \begin{equation}
\mathrm{P}^{(p)}=\textrm{(\{base primes}>7\}\rightarrow1)\circ(\textrm{factorization of}\; C_{q}/\mathrm{\: SCPF}(G_{q+2,1}^{(p)}),\quad p=31,63,\ldots\;.\label{eq:prefix-def}\end{equation}
The resulting factorization beginnings are given in%
\begin{table}[H]
\caption{\label{tab:Prefixes-of-}Prefixes of $C_{q}$ factorization when infix
base primes $(\textrm{those}>7)$ are mapped to 1}
\lyxline{\normalsize}

\medskip{}
\begin{tabular}{lccccccccccccccc}
$p$ & $2^{4}\textrm{-1}$ & $2^{5}\textrm{-1}$ & $2^{6}\textrm{-1}$ & $2^{7}\textrm{-1}$ & $2^{8}\textrm{-1}$ & $2^{9}\textrm{-1}$ & $2^{10}\textrm{-1}$ & $\textrm{\ensuremath{2^{11}}-1}$ & $2^{12}\textrm{-1}$ & $2^{13}\textrm{-1}$ & $2^{14}\textrm{-1}$ & $2^{15}\textrm{-1}$ &  & $\cdots$ & \tabularnewline
\cline{1-15} 
$\mathrm{P}^{(p)}$ & $-$ & $3$ & $3^{2}\cdot5$ & 7 & $3\cdot5^{3}$ & $3^{3}\cdot7$ & $7$ & $3^{3}\cdot5$ & $3^{3}\cdot5^{3}\cdot7^{2}$ & $3^{3}\cdot5^{2}\cdot7^{3}$ & $3^{6}\cdot5^{3}\cdot7^{3}$ & $3^{2}\cdot5^{4}\cdot7$ &  & $\cdots$ & \tabularnewline
 &  &  &  &  &  &  &  &  &  &  &  &  &  &  & \tabularnewline
\end{tabular}\lyxline{\normalsize}
\end{table}
\noindent Now consider the paraorder products $\prod_{r=1}^{n+1}p_{r}$
(which as we shall see in the next section play a vital role in the
interordinal preon model to be presented there) and extract from them
factorization beginnings in analogous fashion, namely\begin{equation}
\coprod\,^{(p')}=\textrm{( \{base primes}>7\}\rightarrow1)\circ(\textrm{factorization of}\;\prod_{r=1}^{n+1}p_{r}),\label{eq:paraprod-def}\end{equation}
as summarized in %
\begin{table}[H]
\caption{Paraorder-product factorization beginnings when base primes > 7 are
mapped to 1}
\lyxline{\normalsize}

\medskip{}
\begin{tabular}{lcccccccccccc}
$p'$ & $2^{2}\textrm{-1}$ & $2^{3}\textrm{-1}$ & $2^{4}\textrm{-1}$ & $2^{5}\textrm{-1}$ & $2^{6}\textrm{-1}$ & $2^{7}\textrm{-1}$ & $2^{8}\textrm{-1}$ & $\textrm{\ensuremath{2^{9}}-1}$ & $2^{10}\textrm{-1}$ & $2^{11}\textrm{-1}$ & $2^{12}\textrm{-1}\quad\cdots$ & \tabularnewline
\cline{1-12} 
$\coprod\,^{(p')}$ & 3 & $3\cdot7$ & $3^{2}\cdot5\cdot7$ & $3^{2}\cdot5\cdot7$ & $3^{4}\cdot5\cdot7^{2}$ & $3^{4}\cdot5\cdot7^{2}$ & $3^{5}\cdot5^{2}\cdot7^{2}$ & $3^{5}\cdot5^{2}\cdot7^{3}$ & $3^{6}\cdot5^{2}\cdot7^{3}$ & $3^{6}\cdot5^{2}\cdot7^{3}$ & $3^{8}\cdot5^{3}\cdot7^{4}\:\cdots$ & \tabularnewline
 &  &  &  &  &  &  &  &  &  &  &  & \tabularnewline
\end{tabular}\lyxline{\normalsize}

\end{table}
\noindent Then we find that $\coprod\,^{(p')}$ can be expressed
as a multiplicative partition $\mathrm{P}^{(r)}\mathrm{P}^{(s)}\cdots$
$(r,s,\ldots\in\{2q+1,4q+3,8q+7,\ldots\})$ of length $2^{m}$ for
some $m$, requiring all instantions of $r,s,\ldots$ to be mutually
distinct. This starts working out as \[
\begin{array}{cccccc}
 & \!\!\!\!\!\!\!\!\coprod\,^{(3)}=\mathrm{P}^{(31)},\\
\coprod\,^{(7)}=\mathrm{P}^{(31)}\mathrm{P}^{(127)}, &  & \coprod\,^{(15)}=\coprod\,^{(31)}=\mathrm{P}^{(63)}\mathrm{P}^{(127)},\\
\coprod\,^{(63)}=\coprod\,^{(127)}=\mathrm{P}^{(31)}\mathrm{P}^{(127)}\mathrm{P}^{(1023)}\mathrm{P}^{(2047)}, &  & \coprod\,^{(255)}=\mathrm{P}^{(63)}\mathrm{P}^{(127)}\mathrm{P}^{(1023)}\mathrm{P}^{(2047)},\end{array}\]
only to stop thereafter:\[
\begin{array}{cccccc}
\coprod\,^{(511)}\neq\mathrm{P}^{(r)}\mathrm{P}^{(s)}\cdots, & \coprod\,^{(1023)}\neq\mathrm{P}^{(r)}\mathrm{P}^{(s)}\cdots, & \textrm{}\end{array}\]
which suggests that $N$ is a (first) bound to the multiplicatice
partition process for $n+1>N=8$. From Tables \ref{tab:The-lowest-dimensions}
and \ref{tab:Periodicity-of-binary} we expect $N'=23$ etc. to act
as (second etc.) rebound to any further partitionability of $\coprod\,^{(2^{n+1}-1)}$
arising in the sequel,%
\footnote{the prefixes $\mathrm{P}^{(65535)}$ and beyond are not easily assessable%
} thus indicating a link between the $\lambda$'s -- interpreted as
base primes $3,5,7$ including map of higher base primes to 1 -- and
$N,N'$.

\section{\label{sec:An-interordinal-preon}An interordinal preon model}

\noindent To our knowledge, Oscar Wallace Greenberg was the first
to recognize that quarks can be viewed as parafermions of order 3.
But with the advent of QCD, and the experimental findings to date
that quarks are pointlike down to $10^{-20}\,$m, preons, parafermionic
or otherwise, have not found much acclaim among physicists. This is
not the place to review the variously theorized preon types in the
literature, including Green's own proposal, nor is the following meant
to be a worked out physical model of the subatomic onion -- it rather
tries to point the way to the putative mathematical structure of hadronic
matter, which likely is interordinal in the Mersennian sense. A more
self-contained elaboration of the present ideas will appear in a separate
paper. Here, we choose the symbols $\mathrm{p}_{\textrm{up}}^{(p)}$
and $\mathrm{p}_{\textrm{down}}^{(p)}$ to denote up-type and down-type
preons of paraorder $p$ respectively.
\begin{conj}
\label{con:Preons-of-order}Preons of order $p_{n+1}$ are either
up-type or down-type, $\mathrm{p}_{\textrm{\emph{up}}}^{(p_{n+1})}$
or $\mathrm{p}_{\textrm{\emph{down}}}^{(p_{n+1})}$. The electric
charge (in \emph{e}) of up-type items is given by expressions $\mathrm{c}_{\mathrm{up}}^{(p_{n+1})}=\dfrac{p_{n+1}-\sum_{s=0}^{n}p_{s}}{\prod_{r=1}^{n+1}p_{r}}=\dfrac{n+1}{\prod_{r=1}^{n+1}p_{r}}$
and the charge of down-type items by $\mathrm{c}_{\mathrm{down}}^{(p_{n+1})}=\dfrac{-\sum_{s=0}^{n}p_{s}}{\prod_{r=1}^{n+1}p_{r}}$.
The charge of up-type items transforms as \textup{$\mathrm{c}_{\mathrm{up}}^{(p_{n})}=(p_{n+1}-1)\mathrm{c}_{\mathrm{up}}^{(p_{n+1})}+\mathrm{c}_{\mathrm{down}}^{(p_{n+1})}$
and the charge of down-type items as $\mathrm{c}_{\mathrm{down}}^{(p_{n})}=(p_{n}+1)\mathrm{c}_{\mathrm{down}}^{(p_{n+1})}+p_{n}\mathrm{c}_{\mathrm{up}}^{(p_{n+1})}$.}
\end{conj}
See the tabularized values:

\begin{table}[h]
\caption{\label{tab:Interordinal-preon-model}Interordinal preon model}
\lyxline{\normalsize}\medskip{}
\begin{tabular}{ccccccccc}
$p_{n+1}$ &  &  & up-type charge  &  &  & down-type charge  &  & \tabularnewline
\hline
$1$ & $\qquad\qquad\qquad\qquad$ & $\qquad\qquad\qquad\qquad$ & $1$ & $\qquad\qquad\qquad\qquad$ & $\qquad\qquad\qquad\qquad$ & $0$ &  & \tabularnewline
$3$ &  &  & $\frac{2}{3}$ &  &  & $-\frac{1}{3}$ &  & \tabularnewline
$7$ &  &  & $\frac{3}{21}$ &  &  & $-\frac{4}{21}$ &  & \tabularnewline
$15$ &  &  & $\frac{4}{315}$ &  &  & $-\frac{11}{315}$ &  & \tabularnewline
$31$ &  &  & $\frac{5}{9765}$ &  &  & $-\frac{26}{9765}$ &  & \tabularnewline
$63$ &  &  & $\frac{6}{615195}$ &  &  & $-\frac{57}{615195}$ &  & \tabularnewline
$\vdots$ &  &  & $\vdots$ &  &  & $\vdots$ &  & \tabularnewline
\end{tabular}\lyxline{\normalsize}
\end{table}
\noindent  Let us start with $p_{n}=7$. According to Conjecture
\ref{con:Preons-of-order}, level-7 preons will contain fourteen $\mathrm{p}_{\textrm{up}}^{(15)}$
and one $\mathrm{p}_{\textrm{down}}^{(15)}$ in the up-type case,
and eight $\mathrm{p}_{\textrm{down}}^{(15)}$ plus seven $\mathrm{p}_{\textrm{up}}^{(15)}$
in the down-type case. For sublevel-3 preons, the constituents are
$\mathrm{6p}_{\textrm{up}}^{(7)}+\mathrm{p}_{\textrm{down}}^{(7)}$
in the up-type and $\mathrm{4p}_{\textrm{down}}^{(7)}+3\mathrm{p}_{\textrm{up}}^{(7)}$
in the down-type case, and for sublevel-1 items, the constituents
are valence quarks, $\mathrm{2p}_{\textrm{up}}^{(3)}+\mathrm{p}_{\textrm{down}}^{(3)}$
for the proton and $\mathrm{2p}_{\textrm{down}}^{(3)}+\mathrm{p}_{\textrm{up}}^{(3)}$
for the neutron. Bearing the lesson of Sect. \emph{\ref{sub:Row-(column)-structure}}
in mind, a larger confining space -- only the $3D$ lower end of which
is familiar to us -- is required to implement the Russian-doll preon
structure. The operators $f^{(15)}$ and $h^{(15)}$ encode the necessary
topological information to successively guide configurations across
extra-dimensional (Euclidean) space. This is propounded in our second
\begin{conj}
\label{con:To-serve-as}To qualify as constituents of a superordinate
preon, hyperspheres must be used in \#'s that divide the kissing number
of the space they live in.
\end{conj}
Conjecture \ref{con:To-serve-as} implies that down-type constituents
are subspace inhabitants relative to the space that up-type constituents
live in. Since there are at least three generations of quarks, there
must be room for enlarged configurations. For instance, while the
down quark with electric charge $-\dfrac{1}{3}$ is a level-3 configuration
on its own, an {}``uptype-preon plus downtype-antipreon'' configuration
formed by level-15 constituents, plus an {}``uptype antipreon''
configuration formed by level-7 constituents, is required to yield
the same result for the strange quark:

\[
[(14\mathrm{p}_{\textrm{up}}^{(15)}+\mathrm{p}_{\textrm{down}}^{(15)})+(8\mathrm{\overline{p}}_{\textrm{down}}^{(15)}+\mathrm{7\overline{p}}_{\textrm{up}}^{(15)})]+(6\mathrm{\overline{p}}_{\textrm{up}}^{(7)}+\mathrm{\overline{p}}_{\textrm{down}}^{(7)}),\]

\noindent combines to yield

\[
\frac{14\cdot4+1\cdot(-11)+8\cdot11+7\cdot(-4)}{315}+\frac{6\cdot(-3)+1\cdot(4)}{21}=-\frac{1}{3}.\]

\noindent Correspondingly, level-31 constituents are part of configurations
that reproduce the charges of quarks of the third generation. Just
as level-15 constituents require $\mathrm{14}_{\textrm{up}}\mid L_{x}$,
$\mathrm{8}_{\textrm{down}}\mid L_{y}$ $(x>y)$, level-31 constituents
require $30_{\textrm{up}}\mid L_{w}$, $\mathrm{16}_{\textrm{down}}\mid L_{z}$
$(w>z)$. According to Table \ref{tab:Stenoscopy-of-kissing}, the
only certified kissing numbers covered by $f^{(31)}$ and $h^{(31)}$
that are divisible by 30 are also divisible by 16: $L_{8}=240,L_{16}=4320$
and $L_{24}=196560.$ Close to the confining $D^{(31)}=31$, there
is a further, if uncertified, one: %
\footnote{{\tt http://www.math.rwth-aachen.de/}$\!\thicksim${\tt Gabriele.Nebe/LATTICES/kiss.html}%
} $L_{29}=207930$. It's divisible by 30 but not by 16. There is the
faint hope that this kissing number will be certified one day.

\section{\label{sec:A-proposal-for}A proposal for a planar geometric model}

\subsection{\label{sub:The-cardioid-and}The cardioid and her arclength}

\begin{figure}[H]

\caption{\label{fig:card-turn-1}$f$- ($h$-) parafermion in cardiodic representation}
\medskip{}
\hfill{}\includegraphics[scale=0.33]{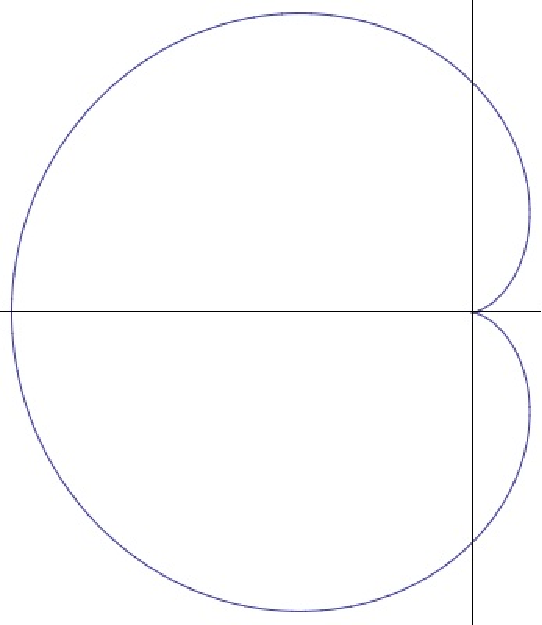}\quad{}\includegraphics[scale=0.25]{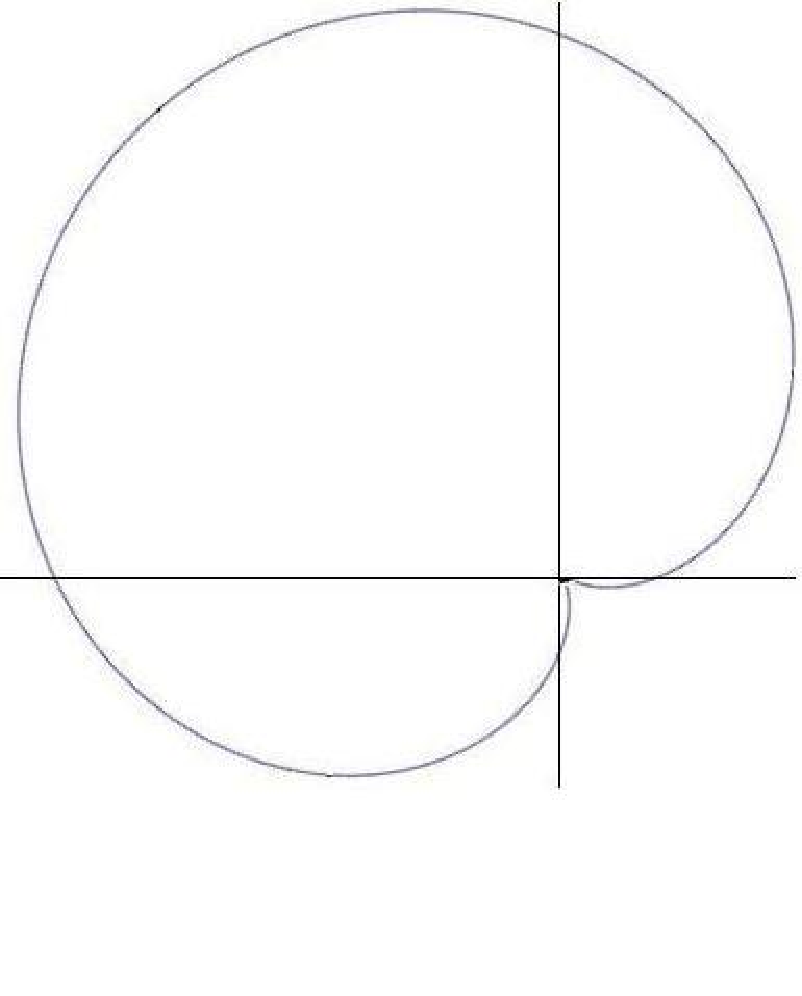}\quad{}\includegraphics[scale=0.25]{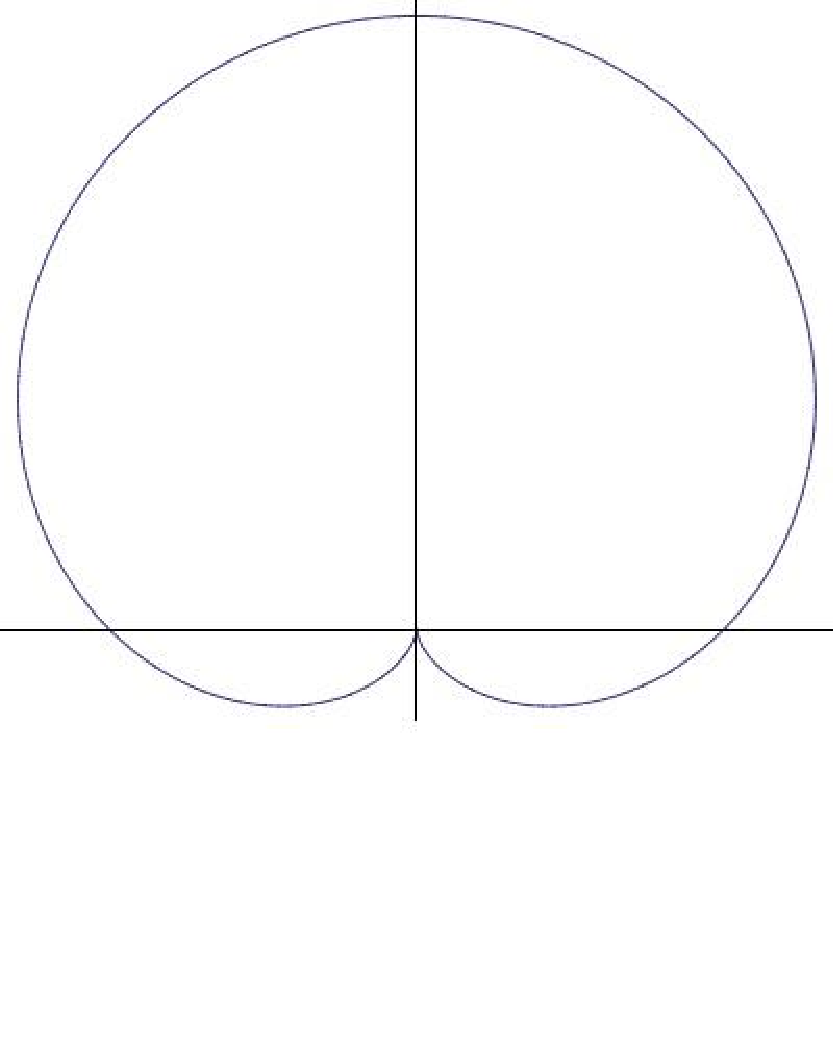}\quad{}\includegraphics[scale=0.33]{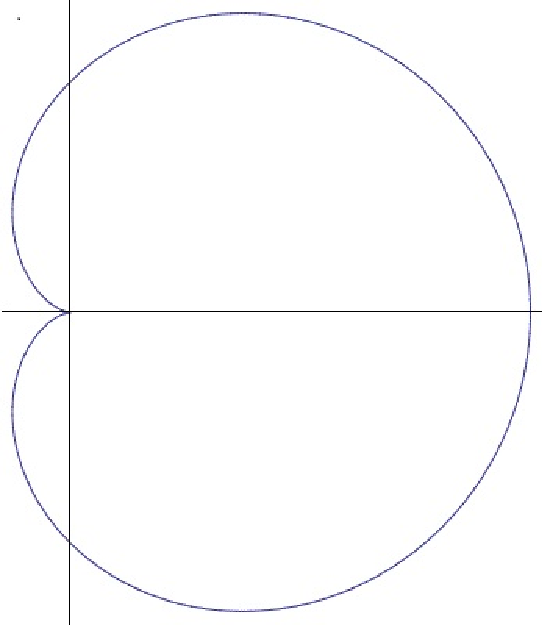}\hfill{}\bigskip{}

\hfill{}{\small$\infty\leftarrow p$\hfill$p=7$\hfill$p=3$\hfill$p=1$}
\hfill{}
\end{figure}

\noindent Apart from the conjectured connection with sphere packing
in Euclidean $D$-space, nilpotent operators such as $f^{(p)}$ and
$h^{(p)}$ have interesting representations in ordinary plane trig\-on\-ometry,
involving cardiods whose cusps are located at the origin. Consider
the rightmost cardioid in Fig. \ref{fig:card-turn-1} which has the
polar representation ($a$ a parameter) \begin{equation}
r=a(1+\cos\theta),\mbox{}\label{eq:Cardioid0}\end{equation}
 and compare it to the cardiod left to it with polar representation
\begin{equation}
r=a(1+\sin\theta).\label{eq:Cardioid-90}\end{equation}
\bigskip{}
\noindent Obviously the transformation implies a quarter-turn around
the origin, and the transformation into the leftmost cardioid a half-turn
or flip-over. Recalling that $c_{3}$ and $c_{2}$, the basic building
blocks of $f^{(p)}$ and $h^{(p)}$ for $p>1$, realize such transformations
in matrix form, we may conclude they lie at the basis of planar representations
of $f^{(p)}$ and $h^{(p)}$. Whichever of the two one uses, they
should be made an infinite process to mirror in a geometric spirit
the forming of the root-$f$- and root-$h$ sequence. Now reflection
is an operation indivisible within the framework of plane trigonometry,
so we are left with rotation as a vehicle to express the infinite
sequence. It would consist, first of no turn, followed by a quarter-turn,
followed by further turns of ever-halving angle measured in radians
as shown in Fig. \ref{fig:card-turn-1} (see the second cardioid from
the left as an example of an intermediate stage in the process corresponding
to $p=7$): \begin{equation}
r_{n}=a\left(1+\cos(\theta+\tfrac{p-1}{p+1}\cdot\pi)\right),\quad n=\log_{2}(p+1),\: p\in\left\{ 1,3,7,\ldots\right\} ,\label{eq:card-trans90-0}\end{equation}
 whence in the limit -- as though it was effected by $c_{2}$ -- the
rightmost cardioid would take the position of the leftmost:\begin{equation}
\lim_{p\rightarrow\infty}a\left(1+\cos(\theta+\tfrac{p-1}{p+1}\cdot\pi)\right)=a(1+\cos(\theta+\pi))=a(1-\cos\theta).\label{eq:card-translimit}\end{equation}
 Let us expound the details of the envisioned representation. The
arc\-length of the cardioid is determined by the integral \begin{equation}
\intop_{a_{1}}^{a_{2}}\mathrm{d}s,\label{eq:arclength}\end{equation}
 where $\mathrm{d}s=\sqrt{r^{2}\mathrm{d\text{\ensuremath{\theta}}^{2}+d\mathit{r^{\mathrm{2}}}}}$.
For the cardioid at rest, $\mathrm{d}s$ equals $2a\cos\frac{\theta}{2}\mathrm{d}\text{\ensuremath{\theta}}$.
This function remains valid if the arc traced between the limits doesn't
cross the cusp or its antipodal point. The maximal admissible interval
for the lower and upper limits thus is {[}$0,\pi]$, and the circumference
becomes \begin{equation}
C=2\intop_{0}^{\pi}\mathrm{d}s=8a;\label{eq:Circum}\end{equation}
 otherwise the circumference turns 0, a result which is in accord
with the nullity of an order-$p$ nilpotent operator for exponents
larger than $p$, and immediately makes clear that this operator exponentiated
has to be represented by a compound of arclengths which eventually
transgresses the 4$a$ boundary. In the case of the cardioid set into
motion, the total angle accruing from counterclockwise rotations according
to Eq. \ref{eq:card-trans90-0} toward its end position does not exceed
$\pi$. To properly map the nilpotence condition, however, we must
after each step use the $x$-axis as an equator%
\footnote{In case of a process $r_{i}=a(1+\sin(\theta+\frac{p-1}{p+1}\cdot\pi))$
the $y$-axis had to be used as equator%
} and 1) separate arc parts from the upper half-plane with lower and
upper azimuths $u_{n}$ and $v_{n}$ from those of the lower half-plane
with azimuths $\bar{u}_{n}$ and $\bar{v}_{n}$, and 2) keep track
of the gap left behind in the upper half-plane by the moving cardioid
with lower and upper azimuths $\textrm{co-}u_{n}$ and $\textrm{co-}v_{n}$.
As for actual arclength computations, an option has to be taken of
either using the cardioid-at-rest arclength function or its cardioid-in-motion
counterpart. Let us first examine option one according to which we
have to compare those arclength parts of the cardioid in motion and
the cardioid at rest that are in correspondence with each other. It
turns out that, to accommodate cardioidic motion, the lower and upper
integral limits have to be determined by the coordinate transformations
\begin{equation}
\begin{array}{l}
\textrm{upper}\,\textrm{half-plane:}\\
(u_{n},v_{n})=\left(\tfrac{p-1}{p+1}\cdot\pi,\pi\right)\rightarrow(w_{n},z_{n})=\left(0,\left(1-\tfrac{p-1}{p+1}\right)\pi\right),\\
(\textrm{co-}u_{n},\textrm{co-}v_{n})=(\textrm{co-}w_{n},\textrm{co-}z_{n})=\left(0,\tfrac{p-1}{p+1}\cdot\pi\right),\\
\textrm{lower}\,\textrm{half-plane:}\\
(\bar{u}_{n},\bar{v}_{n})=\left(\pi,\left(1+\tfrac{p-1}{p+1}\right)\pi\right)\rightarrow(\bar{w}_{n},\bar{z}_{n})=\left(\left(1-\tfrac{p-1}{p+1}\right)\pi,\pi\right),\end{array}\label{eq:rot-to-rest}\end{equation}
 where of course $n=\log_{2}(p+1)$, $p\in\left\{ 1,3,7,\ldots\right\} $
both times. Labeling the corresponding arclengths $A_{n}$ and $\bar{A}_{n}$
and computing them for the first four cardioid stops $r_{1}=a(1+\cos\theta)$,
$r_{2}=a(1+\cos(\theta+\frac{\pi}{2}))$, $r_{3}=a(1+\cos(\theta+\frac{3}{4}\pi))$,
$r_{4}=a(1+\cos(\theta+\frac{7}{8}\pi))$, one finds \[
\begin{array}{ll}
 & A_{1}=\intop_{0}^{\pi}2a\cos\frac{\theta}{2}\mathrm{d}\theta=4a;\\
 & A_{2}=\intop_{0}^{\frac{\pi}{2}}2a\cos\mathrm{\frac{\theta}{2}\mathrm{d}}\theta=4a\cdot\frac{\sqrt{2}}{2};\\
 & A_{3}=\intop_{0}^{\frac{\pi}{4}}2a\cos\frac{\theta}{2}\mathrm{d}\theta=4a\cdot\frac{\sqrt{2-\sqrt{2}}}{2};\\
 & A_{4}=\intop_{0}^{\frac{\pi}{8}}2a\cos\frac{\theta}{2}\mathrm{d}\theta=4a\cdot\frac{\sqrt{2-\sqrt{2+\sqrt{2}}}}{2};\\
\! & \cdots\end{array}\]
 and for the subequatorials, \[
\begin{array}{ccl}
 &  & \bar{A}_{1}=\intop_{\pi}^{\pi}2a\cos\frac{\theta}{2}\mathrm{d}\theta=0;\\
 &  & \bar{A}_{2}=\intop_{\frac{\pi}{2}}^{\pi}2a\cos\frac{\theta}{2}\mathrm{d}\theta=4a\left(1-\frac{\sqrt{2}}{2}\right);\\
 &  & \bar{A}_{3}=\intop_{\frac{\pi}{4}}^{\pi}2a\cos\frac{\theta}{2}\mathrm{d}\theta=4a\left(1-\frac{\sqrt{2-\sqrt{2}}}{2}\right);\\
 &  & \bar{A}_{4}=\intop_{\frac{\pi}{8}}^{\pi}2a\cos\frac{\theta}{2}\mathrm{d}\theta=4a\left(1-\frac{\sqrt{2-\sqrt{2+\sqrt{2}}}}{2}\right);\\
 &  & \qquad\qquad\qquad\qquad\cdots\end{array}\]
 A nilpotent operator of order $p\in\left\{ 1,3,7,\ldots\right\} $
then is representable by $A_{n}$ $\:(n=\log_{2}(p+1)$), and the
action on itself by the operation\begin{equation}
(A_{n},A_{n})\equiv2\,\frac{A_{n}\cdot\textrm{co-}A_{n}}{A_{n}+\bar{A}_{n}}=A_{n-1}.\label{eq:f-squared}\end{equation}
 The auxiliary expressions $\textrm{co-}A_{n}$ (for co-arclength
on the cardioid fixed at rest) are given by \[
\begin{array}{lcl}
 &  & \textrm{co-}A_{1}=\intop_{0}^{0}2a\cos\frac{\theta}{2}\mathrm{d}\theta=0;\\
 &  & \textrm{co-}A_{2}=\intop_{0}^{\frac{\pi}{2}}2a\cos\frac{\theta}{2}\mathrm{d}\theta=4a\cdot\frac{\sqrt{2}}{2};\\
 &  & \textrm{co-}A_{3}=\intop_{0}^{\frac{3}{4}\pi}2a\cos\frac{\theta}{2}\mathrm{d}\theta=4a\cdot\frac{\sqrt{2+\sqrt{2}}}{2};\\
 &  & \textrm{co-}A_{4}=\intop_{0}^{\frac{7}{8}\pi}2a\cos\frac{\theta}{2}\mathrm{d}\theta=4a\cdot\frac{\sqrt{2+\sqrt{2+\sqrt{2}}}}{2};\end{array}\]
 \[
\cdots\]
 and obey the Vieta condition, i.e., for $a=\frac{1}{4}$ form the
Euler product\begin{equation}
\lim_{n\rightarrow\infty}\prod_{i=2}^{n}\textrm{co-}A_{i}=\frac{2}{\pi}.\label{eq:vieta}\end{equation}
 The {}``unmoved-mover'' representation constructed this way seems
to be akin to the root-$f$ sequence since the integration limits
are derived by coordinate rotation (see transformations (\ref{eq:rot-to-rest})).
Following this reasoning, a {}``moving-mirror'' representation that
was akin to the root-$h$ sequence would be expected to ensue from
option two with integration limits derived by coordinate reflection.
To see if this is true we introduce the arclength function of the
cardioid in motion, $\mathrm{d}s=2a\cos(\frac{\theta}{2}+\frac{p-1}{p+1}\cdot\frac{\pi}{2})\mathrm{d}\text{\ensuremath{\theta}}$,
and look what else is needed to reproduce identical results in terms
of arclengths as this function is used. It turns out that in all azimuth-to-integration
limit transformations there is a flip over the equator ($x$-axis)
in this case,\begin{equation}
\begin{array}{l}
\textrm{upper}\,\textrm{half-plane:}\\
(u_{n},v_{n})=\left(\tfrac{p-1}{p+1}\cdot\pi,\pi\right)\rightarrow(w_{n},z_{n})=\left(\left(2-\tfrac{p-1}{p+1}\right)\pi,\pi\right),\\
(\textrm{co-}u_{n},\textrm{co-}v_{n})=\left(0,\tfrac{p-1}{p+1}\cdot\pi\right)\rightarrow(\textrm{co-}w_{n},\textrm{co-}z_{n})=\left(2\pi,\left(2-\tfrac{p-1}{p+1}\right)\pi\right),\\
\textrm{lower}\,\textrm{half-plane:}\\
(\bar{u}_{n},\bar{v}_{n})=\left(\pi,\left(1+\tfrac{p-1}{p+1}\right)\pi\right)\rightarrow(\bar{w}_{n},\bar{z}_{n})=\left(\pi,\left(1-\tfrac{p-1}{p+1}\right)\pi\right),\end{array}\label{eq:rot-reflect}\end{equation}

\vspace*{0.4cm}\noindent confirming the expectation that the corresponding
set of arclength formulae is akin to the root-$h$ sequence. This
conclusion is also supported by interordinality considerations. Recalling
that for the root-$f$ sequence the carry-bit neighborhood link is
characterized by identities that either are purely interordinal (Eq.
(\ref{eq:inter})) or intraordinal (Eq. (\ref{eq:intra})) whereas
for the root-$h$ sequence the corresponding link mixes interordinal
and intraordinal relationship (Eq. \ref{eq:inter-intra-h}), we can
observe a similar phenomenon in the present planar representations:
For the upper definite-integral limit $\textrm{co-}z_{n'}$, e.g.,
we find \begin{equation}
\begin{array}{l}
\textrm{in the ``unmoved-mover'' representation:}\\
\textrm{co-}z_{n'}=\frac{p}{p+1}\cdot\pi,\\
\textrm{and in the ``moving-mirror'' representation:}\\
\textrm{co-}z_{n'}=\left(\frac{p'+1}{p+1}-\frac{p}{p+1}\right)\pi,\end{array}\label{eq:inter-intra-cardioid}\end{equation}

\vspace*{0.4cm}\noindent where $n'=\log_{2}(p'+1)$, $p'=2p+1$ and
$p\in\left\{ 1,3,7,15,\ldots\right\} $, which completes the analogy.

\subsection{\label{sub:Continued-fraction-representation-}Cardioidic arclength
and its relation to Catalan structure}

\noindent When focusing on cardioidic (co-)arclength per se:\begin{equation}
\begin{array}{c}
A_{n}=4a\sin\frac{\pi}{p+1},\\
\textrm{co-}A_{n}=4a\cos\frac{\pi}{p+1},\end{array}\label{eq:arcs-per-se}\end{equation}

\noindent no such distinction can be made. This aside, simple features
of the Catalan structure such as the departure from $G_{\mu\nu}=1$
and onset of $(q,p)$-interordinality at $p=15$ still find a planar
analogy, namely in the loss of homogeneity for $n>3$ of the factoring
of $x^{2^{n}}+y^{2^{n}}$ -- whose limit contour is a square centered
at the origin of the $x$-$y$ plane -- into polynomials whose contours
are diagonally-oriented crossing ellipses%
\footnote{see \cite{ShuKo10} for further study of crossing ellipses%
} with eccentricities rooted in Eqs. (\ref{eq:arcs-per-se}). In what
follows we always assume $a=\frac{1}{4}$. Then the first approximation
to the square is an inscribed circle:\[
x^{2}+2\textrm{co-}A_{1}xy+y^{2}\]
($\textrm{co-}A_{1}=0$ marks the degeneracy of the case); the second
approximation is given by\[
\begin{array}{rcl}
x^{4}+y^{4} & = & (x^{2}+2\textrm{co-}A_{2}xy+y^{2})(x^{2}-2\textrm{co-}A_{2}xy+y^{2})=(x^{2}+2A_{2}xy+y^{2})(x^{2}-2A_{2}xy+y^{2});\end{array}\]

\noindent and the third one by\[
x^{8}+y^{8}=(x^{2}+2\textrm{co-}A_{3}xy+y^{2})(x^{2}-2\textrm{co-}A_{3}xy+y^{2})(x^{2}+2A_{3}xy+y^{2})(x^{2}-2A_{3}xy+y^{2}).\]

\noindent As early as in the next higher instance, however, interpolating
(co-)arclengths -- as a footprint of $(p,q)$-interordinality -- \[
\begin{array}{ccc}
\tilde{A}_{n-1} & \equiv & A_{n}\textrm{co-}A_{n_{q}}+\textrm{co-}A_{n}A_{n_{q}}=4a\cdot\sin(\frac{\pi}{p+1}+\frac{\pi}{q+1}),\\
\textrm{co-}\tilde{A}_{n-1} & \equiv & \textrm{co-}A_{n}\textrm{co-}A_{n_{q}}-A_{n}A_{n_{q}}=4a\cdot\cos(\frac{\pi}{p+1}+\frac{\pi}{q+1})\end{array}\quad(n_{q}=n-2)\]

\noindent  join in: \[
\begin{array}{rcl}
x^{16}+y^{16} & = & (x^{2}+2\textrm{co-}A_{4}xy+y^{2})(x^{2}-2\textrm{co-}A_{4}xy+y^{2})(x^{2}+2\textrm{co-}\tilde{A}_{3}xy+y^{2})(x^{2}-2\textrm{co-}\tilde{A}_{3}xy+y^{2})\cdot\\
 &  & \cdot\:(x^{2}+2A_{4}xy+y^{2})(x^{2}-2A_{4}xy+y^{2})(x^{2}+2\tilde{A}_{3}xy+y^{2})(x^{2}-2\tilde{A}_{3}xy+y^{2}).\end{array}\]

\noindent Further insight is gained by considering the continued
fraction expansions\begin{equation}
A_{n}=a_{0}^{(n)}+\frac{1}{a_{1}^{(n)}+{\displaystyle \frac{1}{\; a_{2}^{(n)}{\displaystyle +\frac{1}{a_{3}^{(n)}+\ddots}}}}}\;\equiv\;[a_{0}^{(n)};a_{1}^{(n)},a_{2}^{(n)},a_{3}^{(n)},\ldots]\quad(n>1;\; a={\textstyle \frac{1}{4}})\label{eq:arclength-contfrac}\end{equation}

\noindent in conjunction with the accordingly defined continued fractions\[
\begin{array}{cc}
\textrm{co-}A_{n}\equiv[\textrm{co-}a_{0}^{(n)};\textrm{co-}a_{1}^{(n)},\textrm{co-}a_{2}^{(n)},\ldots], & \bar{A}_{n}\equiv[\bar{a}_{0}^{(n)};\bar{a}_{1}^{(n)},\bar{a}_{2}^{(n)},\ldots],\\
\textrm{co-}A_{n}^{\,2}\equiv[\textrm{sqco-}a_{0}^{(n)};\textrm{sqco-}a_{1}^{(n)},\textrm{sqco-}a_{2}^{(n)},\ldots], & A_{n}^{\,2}=1-\textrm{co-}A_{n}^{\,2}\equiv[\overline{\textrm{sqco-}a}_{0}^{(n)};\overline{\textrm{sqco-}a}_{1}^{(n)},\overline{\textrm{sqco-}a}_{2}^{(n)},\ldots],\end{array}\]
 and their associated identities for $n>2$,\begin{equation}
\begin{array}{cc}
\textrm{co-}a_{1}^{(n)}=\bar{a}_{1}^{(n)}=1,a_{1}^{(n)}=1+\bar{a}_{2}^{(n)}, & \; a_{\alpha}^{(n)}=\bar{a}_{\alpha+1}^{(n)}\enskip\mathrm{for}\:\alpha>1,\end{array}\label{eq:cf-ident1}\end{equation}
\begin{equation}
\left\lfloor (1+\textrm{co-}a_{2}^{(n)})/2\right\rfloor =\textrm{sqco-}a_{2}^{(n)}\enskip\mathrm{for}\enskip n\equiv1\,\mathrm{mod\,4}\enskip\mathrm{else\,\left\lfloor (1+\textrm{co-}\mathit{a}_{2}^{(\mathit{n})})/2\right\rfloor =\overline{\textrm{sqco-}\mathit{a}}_{1}^{(\mathit{n})}},\label{eq:cf-ident2}\end{equation}
as well as the special cases\[
A_{2}=\textrm{co-}A_{2}=[0;1,\bar{2}],\quad\bar{A}_{2}=[0;3,\bar{2}].\]
 Then the leading $A_{n}$ CF expansion coefficient will be found
to mimic a carry-bit neighborhood $p'=2p+1$:\begin{equation}
a_{1}^{(n+1)}=2a_{1}^{(n)}+1+\delta_{1}^{(n)}\label{eq:generalized-neigh}\end{equation}

\noindent where $\delta_{1}^{(n)}\in\{-3,-1,0\}$; and the next-to-leading
$\textrm{co-}A_{n}$ CF expansion coefficient a second-closest-carry-bit
neighborhood $p''=4p+3$:

\begin{equation}
\textrm{co-}a_{2}^{(n+1)}=4\textrm{co-}a_{2}^{(n)}+3+\delta_{2}^{(n)}\label{eq:generalized-next-but-one}\end{equation}

\noindent where $\delta_{2}^{(n)}\in\{-1,0,1,2,3\}$. As we shall
see, these expansion coefficients can directly be tied to the intrinsic
Catalan structure of $(G_{\mu\nu}^{(p)})$ (or $(J_{\mu\nu}^{(p)})$).
This is because the carry-bit neighborhood $p'=2p+1$ and its extension,
the second-closest-carry-bit neighborhood $p''=4p+3=2p'+1=2(2p+1)+1$,
have structural analogs in the interordinal identities (\ref{eq:inter},\ref{eq:sine-like})
and intraordinal identities (\ref{eq:intra},\ref{eq:cosine-like}):
In the light of Eqs. (\ref{eq:generalized-neigh},\ref{eq:generalized-next-but-one}),
these identities can be dubbed sine-like and cosine-like respectively.%
\footnote{For $p''\!\!=31$, for instance, the assignments $\widetilde{p''}:=$
$\left(\begin{array}{cc}
4905c_{3} & 1633c_{3}\\
15297c_{3} & 4905c_{3}\end{array}\right)$ , $\widetilde{p}:=$ $\left(\begin{array}{cc}
c_{3} & c_{3}\\
c_{3} & c_{3}\end{array}\right)$, $1_{\widetilde{p}}:=$ $\left(\begin{array}{cc}
41c_{3} & 17c_{3}\\
113c_{3} & 41c_{3}\end{array}\right)$ and $1_{\widetilde{p'}}:=$$\left(\begin{array}{cc}
4819c_{3} & 1595c_{3}\\
15067c_{3} & 4819c_{3}\end{array}\right)$ furnish the structural analog of second-closest-carry-bit neighborhood:
$\widetilde{p''}=2(2\widetilde{p}+1_{\widetilde{p}})+1_{\widetilde{p'}}$.%
}We conveniently harmonize Eqs. (\ref{eq:generalized-neigh},\ref{eq:generalized-next-but-one})
by constraining Eq. (\ref{eq:generalized-neigh}) to its second-closest-neighbor
form:\begin{equation}
a_{1}^{(2n+r)}=4a_{1}^{(2n+r-2)}+3+2\delta_{1}^{(2n+r-2)}+\delta_{1}^{(2n+r-1)}\qquad(n>2;\: r\in\{0,1\}).\label{eq:4-form}\end{equation}
\newpage{}\noindent The last two terms in Eq. (\ref{eq:4-form})
suggest that, in Eq. (\ref{eq:generalized-next-but-one}), $\delta_{2}^{(n)}$
may be resolved in a cognate way, \emph{viz}.\begin{equation}
\delta_{2}^{(n)}=2\delta_{1,\textrm{co-}a_{3}^{(n)}}+\varepsilon_{\textrm{bool}_{1}^{(n)},\textrm{bool}_{2}^{(n)}}\label{eq:cog-form}\end{equation}
where $\delta_{ab}$ is the Kronecker symbol and $\varepsilon_{\textrm{bool}_{1}\textrm{bool}_{2}}$
a Levi-Civita symbol with $\varepsilon_{\textrm{FF}}=\varepsilon_{\textrm{TT}}=0,\:\varepsilon_{\textrm{FT}}=1,\varepsilon_{\textrm{TF}}=-1.$
It turns out that the boolean data in question ensue from the truth
values of distinct inequalities $\textrm{co-}a_{i}^{(n)}\gtrless\textrm{co-}a_{i+1}^{(n)}$,
$\textrm{co-}a_{i+2}^{(n)}\gtrless\textrm{co-}a_{i+3}^{(n)}$ in modulo-8
arithmetics. Their onset at $\textrm{co-}a_{i}^{(n)}$ is determined
by the index of the Fibonacci number $F_{i}$ as $n$ progresses and
is preserved $F_{i}$ times. Following the pattern $F_{3,},F_{4},F_{5,},\ldots$,
we for $F_{i}\neq0\,(\textrm{mod}\,8)$ find unidirectionality of
inequality pairs whose sense alternates with progressing $n$, whereas
for $F_{i}=0\,(\textrm{mod}\,8)$, we see same-sense bidirectionality
within subdivisions $F_{i-1}$, $F_{i-2}$ and sense reversal at the
boundary. See Table below where onset is emphasized by a vertical
bar: \medskip{}
\begin{table}[H]
\caption{\label{tab:min-alpha-coeff-1-1} A Fibonacci way of representation
of $\delta_{2}^{(n)}=2\delta_{1,\textrm{co-}a_{3}^{(n)}}+\varepsilon_{\textrm{bool}_{1}^{(n)},\textrm{bool}_{2}^{(n)}}$
using modulo-8 arithmetics}

\medskip{}
\lyxline{\normalsize}

\noindent $\enskip F_{i}\qquad\qquad n=\log_{2}(p+1)$$\,$$\qquad\qquad\qquad\delta_{2}^{(n)}$$\qquad\qquad\qquad\qquad\qquad\textrm{\qquad co-}A_{n}\qquad\qquad\qquad\qquad\qquad\qquad\quad\delta_{1,\textrm{co-}a_{3}^{(n)}}+\varepsilon_{\textrm{bool}_{1}^{(n)},\textrm{bool}_{2}^{(n)}}$

\lyxline{\normalsize}\vspace{1mm}

\begin{tabular}{ccccc}
$F_{3}=2\qquad$ & \qquad{}3\qquad{}\qquad{} & \qquad{}0\qquad{} & \qquad{}\qquad{}{[}0;1,$\underrightarrow{12,}$$\left|7,\right.\!$\negthinspace{}3,2,1,...{]} & \qquad{}\qquad{}$2\delta_{1,7}+\varepsilon_{7<3,2<1}=0\qquad\qquad$\qquad{}\tabularnewline
 & \qquad{}4\qquad{}\qquad{} & \qquad{}-1\qquad{} & \qquad{}\qquad{}{[}0;1,51,$\left|23,\right.\!$\negthinspace{}43,8,1,...{]} & \qquad{}\qquad{}$2\delta_{1,7}+\varepsilon_{7>3,0>1}=-1$\qquad{}\qquad{}\qquad{}\tabularnewline
\hline 
$F_{4}=3\qquad$ & \qquad{}5\qquad{}\qquad{} & \qquad{}2\qquad{} & \qquad{}\qquad{}{[}0;1,206,$\underrightarrow{1,}$$\left|2,\right.\!$\negthinspace{}18,1,1,...{]} & \qquad{}\qquad{}$2\delta_{1,1}+\varepsilon_{2<2,1<1}=2$\qquad{}\qquad{}\qquad{}\tabularnewline
 & \qquad{}6\qquad{}\qquad{} & \qquad{}0\qquad{} & \qquad{}\qquad{}{[}0;1,829,5,$\left|3,\right.\!$\negthinspace{}1,2,1,...{]} & \qquad{}\qquad{}$2\delta_{1,5}+\varepsilon_{3>1,2>1}=0$\qquad{}\qquad{}\qquad{}\tabularnewline
 & \qquad{}7\qquad{}\qquad{} & \qquad{}0\qquad{} & \qquad{}\qquad{}{[}0;1,3319,3,$\left|1,\right.\!$\negthinspace{}6,32,1,...{]} & \qquad{}\qquad{}$2\delta_{1,3}+\varepsilon_{1<6,0<1}=0$\qquad{}\qquad{}\qquad{}\tabularnewline
\hline 
$F_{5}=5\qquad$ & \qquad{}8\qquad{}\qquad{} & \qquad{}1\qquad{} & \qquad{}\qquad{}{[}0;1,13279,1,$\underrightarrow{1,}$$\left|6,\right.\!$\negthinspace{}3,1,1,...{]} & \qquad{}\qquad{}$2\delta_{1,1}+\varepsilon_{6>3,1>1}=1$\qquad{}\qquad{}\qquad{}\tabularnewline
 & \qquad{}9\qquad{}\qquad{} & \qquad{}2\qquad{} & \qquad{}\qquad{}{[}0;1,53120,1,1,$\left|1,\right.\!$\negthinspace{}5,10,13,...{]} & \qquad{}\qquad{}$2\delta_{1,1}+\varepsilon_{1<5,2<5}=2$\qquad{}\qquad{}\qquad{}\tabularnewline
 & \qquad{}10\qquad{}\qquad{} & \qquad{}-1\qquad{} & \qquad{}\qquad{}{[}0;1,212485,11,6,$\left|6,\right.\!$\negthinspace{}1,2,2,...{]} & \qquad{}\qquad{}$2\delta_{1,3}+\varepsilon_{6>1,2>2}=-1$\qquad{}\qquad{}\qquad{}\tabularnewline
 & \qquad{}11\qquad{}\qquad{} & \qquad{}2\qquad{} & \qquad{}\qquad{}{[}0;1,849942,1,6,$\left|16,\right.\!$\negthinspace{}1,13,7,...{]} & \qquad{}\qquad{}$2\delta_{1,1}+\varepsilon_{0<1,5<7}=2$\qquad{}\qquad{}\qquad{}\tabularnewline
 & \qquad{}12\qquad{}\qquad{} & \qquad{}3\qquad{} & \qquad{}\qquad{}{[}0;1,3399773,1,14,$\left|225,\right.\!$\negthinspace{}1,5,2,...{]} & \qquad{}\qquad{}$2\delta_{1,1}+\varepsilon_{1>1,5>2}=3$\qquad{}\qquad{}\qquad{}\tabularnewline
\hline 
$F_{6}=8\qquad$ & \qquad{}13\qquad{}\qquad{} & \qquad{}0\qquad{} & \qquad{}\qquad{}{[}0;1,13599098,4,3,$\underrightarrow{1,}$$\left|1,\right.\!$\negthinspace{}13,1,1,...{]} & \qquad{}\qquad{}$2\delta_{1,4}+\varepsilon_{1>5,1<1}=0$\qquad{}\qquad{}\qquad{}\tabularnewline
 & \qquad{}14\qquad{}\qquad{} & \qquad{}1\qquad{} & \qquad{}\qquad{}{[}0;1,54396395,2,3,3,$\left|1,\right.\!$\negthinspace{}2,1,3,...{]} & \qquad{}\qquad{}$2\delta_{1,2}+\varepsilon_{1>2,1<3}=1$\qquad{}\qquad{}\qquad{}\tabularnewline
 & \qquad{}15\qquad{}\qquad{} & \qquad{}0\qquad{} & \qquad{}\qquad{}{[}0;1,217585584,4,3,1,$\left|4,\right.\!$\negthinspace{}6,11,1,...{]} & \qquad{}\qquad{}$2\delta_{1,4}+\varepsilon_{4>6,3<1}=0$\qquad{}\qquad{}\qquad{}\tabularnewline
 & \qquad{}16\qquad{}\qquad{} & \qquad{}1\qquad{} & \qquad{}\qquad{}{[}0;1,870342339,2,3,1,$\left|1,\right.\!$\negthinspace{}3,1,14,...{]} & \qquad{}\qquad{}$2\delta_{1,2}+\varepsilon_{1>3,1<6}=1$\qquad{}\qquad{}\qquad{}\tabularnewline
 & \qquad{}17\qquad{}\qquad{} & \qquad{}0\qquad{} & \qquad{}\qquad{}{[}0;1,3481369360,3,1,17,$\left|4,\right.\!$\negthinspace{}7,2,17,...{]} & \qquad{}\qquad{}$2\delta_{1,3}+\varepsilon_{4>7,2<1}=0$\qquad{}\qquad{}\qquad{}\tabularnewline
\cline{3-5} 
 & \qquad{}18\qquad{}\qquad{} & \qquad{}1\qquad{} & \qquad{}\qquad{}{[}0;1,13925477443,1,1,17,$\left|2,\right.\!$\negthinspace{}30,2,4,...{]} & \qquad{}\qquad{}$2\delta_{1,1}+\varepsilon_{2<6,2>4}=1$\qquad{}\qquad{}\qquad{}\tabularnewline
 & \qquad{}19\qquad{}\qquad{} & \qquad{}1\qquad{} & \qquad{}\qquad{}{[}0;1,55701909776,1,1,3,$\left|1,\right.\!$\negthinspace{}494,8,1,...{]} & \qquad{}\qquad{}$2\delta_{1,1}+\varepsilon_{1<6,0>1}=1$\qquad{}\qquad{}\qquad{}\tabularnewline
 & \qquad{}20\qquad{}\qquad{} & \qquad{}2\qquad{} & \qquad{}\qquad{}{[}0;1,222807639108,1,2,1,$\left|1,\right.\!$\negthinspace{}1,1,30,...{]} & \qquad{}\qquad{}$2\delta_{1,1}+\varepsilon_{1<1,1>6}=2$\qquad{}\qquad{}\qquad{}\tabularnewline
\hline
$F_{7}=13\qquad$ & \qquad{}21\qquad{}\qquad{} & \qquad{}1\qquad{} & \qquad{}\qquad{}{[}0;1,891230556437,2,1,1,$\underrightarrow{3,}$$\left|7,\right.\!$\negthinspace{}2,5,23,...{]} & \qquad{}\qquad{}$2\delta_{1,2}+\varepsilon_{7<2,5<7}=1$\qquad{}\qquad{}\qquad{}\tabularnewline
 & \qquad{}. . .\qquad{}\qquad{} &  &  & \qquad{}\tabularnewline
\end{tabular}\vspace{3mm}
\lyxline{\normalsize}\medskip{}

\end{table}

\noindent On closer inspection, it appears that $\delta_{2}^{(n)}$
provides the clue to the envisioned continued-fraction Catalan-structure
link: 
\begin{conj}
Let $\mathrm{LL}(G_{\mu\nu}^{(p)})$ be the product of a $(\frac{p+1}{4})\times m_{c}^{(n)}$
matrix $(e_{rs})$ and a $m_{c}^{(n)}\times(\frac{p+1}{4})$ matrix
$(\chi_{sr}).$ Then $e_{rs}\in\{\delta_{2}^{(2)},\ldots,\delta_{2}^{(n)}\}$
and $\chi_{sr}\in\{C_{\min(1,q)},\ldots C_{2q}\},$ where $p=3,7,15,\ldots,\: q=\frac{p-3}{4},\: n=\log_{2}(p+1)$. 
\end{conj}
This is trivially true in the case $p=3$, where $\textrm{LL}(G_{\mu\nu}^{(3)})=C_{0}=1$
coincides with $\delta_{2}^{(2)}=1$ in the computation: $\textrm{co-}a_{2}^{(3)}=12=4\textrm{co-}a_{2}^{(2)}+3+\delta_{2}^{(2)}=4\cdot2+3+1$.
For paraorder seven, or $n=3$, every $G_{\mu\nu}^{(7)}$ from $\textrm{LL}(G_{\mu\nu}^{(7)})$
can be represented by a dot product of a vector containing two elements
$\in\{1\}\cup\{\delta_{2}^{(3)}\}$, where $\delta_{2}^{(3)}=0$,
and the vector $(C_{1},C_{1})$. Hence also trivially:\begin{equation}
\textrm{LL}(G_{\mu\nu}^{7)})=\left(\begin{array}{cc}
1 & 0\\
0 & 1\end{array}\right)\left(\begin{array}{cc}
C_{1} & C_{1}\\
C_{1} & C_{1}\end{array}\right);\label{eq:G7-epitom}\end{equation}
 and the $\delta_{2}^{(3)}$ value added here coincides with the one
used in the computation: $\textrm{co-}a_{2}^{(4)}=4\textrm{co-}a_{2}^{(3)}+3+\delta_{2}^{(3)}=4\cdot12+3+0=51$.
For paraorder fifteen, every $G_{\mu\nu}^{(15)}$ from $\textrm{LL}(G_{\mu\nu}^{(15)})$
can be represented by a dot product of a vector with elements $\in\{0,1\}\cup\{\delta_{2}^{(4)}\}$,
where $\delta_{2}^{(4)}=-1$, and vector $(C_{1,},C_{2,},\ldots,C_{6})$.
It turns out that we have to use $4\!\times\!11$ and $11\!\times\!4$
matrices to, for the first time nontrivially, epitomize $\textrm{LL}(G_{\mu\nu}^{(15)})$
in product form:\begin{equation}
\textrm{LL}(G_{\mu\nu}^{15)})=\left(\begin{array}{ccccccccccc}
1 & 0 & 0 & 0 & -1 & 0 & 0 & 0 & 0 & 0 & 0\\
0 & 1 & 0 & 0 & 0 & -1 & 1 & 0 & 0 & 0 & 0\\
0 & 0 & 1 & 0 & 0 & 0 & 0 & -1 & 1 & 0 & 0\\
0 & 0 & 0 & 1 & 0 & 0 & 0 & 0 & 0 & -1 & -1\end{array}\right)\left(\begin{array}{cccc}
C_{3} & C_{3} & C_{1} & C_{1}\\
C_{4} & C_{3} & C_{1} & C_{1}\\
C_{5} & C_{4} & C_{3} & C_{3}\\
C_{6} & C_{5} & C_{4} & C_{3}\\
0 & C_{2} & 0 & 0\\
C_{3} & 0 & 0 & 0\\
C_{2} & 0 & 0 & 0\\
C_{1} & C_{2} & 0 & C_{2}\\
0 & C_{3} & 0 & 0\\
C_{4} & C_{1} & C_{2} & 0\\
C_{3} & 0 & C_{1} & 0\end{array}\right);\label{eq:G15-epitom}\end{equation}
 here, the $\delta_{2}^{(4)}$ coincides with the one used in $\textrm{co-}a_{2}^{(5)}=4\textrm{co-}a_{2}^{(4)}+3+\delta_{2}^{(4)}=4\cdot51+3-1=206$.
Continuing for paraorder thirty-one, every $G_{\mu\nu}^{(31)}$ can
be represented by a dot product of a vector with elements $\in\{-1,0,1\}\cup\{\delta_{2}^{(5)}\}$,
where $\delta_{2}^{(5)}=2$, and vector $(C_{1,},C_{2,},\ldots,C_{14})$.
We find that it takes $8\!\times\!48$ and $48\!\times\!8$ matrices
to render $\textrm{LL}(G_{\mu\nu}^{31)})$ in product form; saving
space, we name only those $G_{\mu\nu}^{(31)}$ that require $\delta_{2}^{(5)}=2$,
namely the secondary-diagonal entries from $\textrm{LL(LL}(G_{\mu\nu}^{31)}))$;
but its presence benefits other entries and has the dimensions of
the factor matrices shrink from $8\!\times\!51$ and $51\!\times\!8$
to $8\!\times\!48$ and $48\!\times\!8$. And so on.

\noindent The number of columns in the left (rows in the right) matrix
factor, $m_{c}^{(n)}$, can now be linked to both leading/next-to-leading
continued-fraction coefficients, $a_{1}^{(n_{a})}$ and $\textrm{co-}a_{2}^{(n_{b})}$,
via the harmonized second-closest-carry-bit neighborhood equation\begin{equation}
m_{c}^{(n+2)}=\textrm{co-}a_{2}^{(n+1)}-a_{1}^{(2n-2)}+n-1\qquad(n>1).\label{eq:n_c-relation}\end{equation}
Thus $m_{c}^{(4)}=12-2+1=11$; $m_{c}^{(5)}=51-5+2=48$; also $m_{c}^{(6)}=206-20+3=189$;
and so on. We further note that the associated paraorder is always
contained among the infinity of coefficients $a_{\alpha}^{(n)}$ or
$\textrm{co-}a_{\alpha}^{(n)}$ for some $\alpha$, where the connection
between value and least place is provided by $m_{c}^{(n)}$. It appears
that relation (\ref{eq:cf-ident2}) decides which of the two is the
representative of $p$ with $\alpha$ least. Where $\alpha_{\textrm{least}}=m_{c}^{(n)}$,
$n>3$: if $\left\lfloor (1+\textrm{co-}a_{2}^{(n-1)})/2\right\rfloor =\textrm{sqco-}a_{2}^{(n-1)}$,
which is the case for $n-1\equiv1\,\mathrm{mod\,4}$, then $p=a_{\alpha_{\textrm{least}}}^{(n)}$
else $\left\lfloor (1+\textrm{co-}a_{2}^{(n-1)})/2\right\rfloor =\overline{\textrm{sqco-}a}_{1}^{(n-1)}$
and $p=\textrm{co-}{a}_{\alpha_{\textrm{least}}}^{(n)}$, as shown
in Table \ref{tab:min-alpha-coeff-1}:\medskip{}
\begin{table}[H]
\caption{\label{tab:min-alpha-coeff-1} Representation of $p$ by $a_{\alpha_{\mathrm{\textrm{least}}}}^{(n)}$
or $\textrm{co-}a_{\alpha_{\mathrm{\textrm{least}}}}^{(n)}$, based
on the hypothesis: $\alpha_{\mathrm{least}}=m_{c}^{(n)},\quad n>3$}

\medskip{}
\lyxline{\normalsize}

\noindent $n=\log_{2}(p+1)$$\,$$\qquad\qquad\qquad\qquad m_{c}^{(n)}$$\qquad\qquad\qquad\qquad\qquad\qquad p$

\lyxline{\normalsize}\vspace{1mm}
$\begin{array}{cccccc}
\quad3 & \qquad\qquad &  & \qquad\qquad\qquad\qquad- &  & \qquad\qquad\qquad\textrm{co-}a_{3}^{(3)}=7\\
\quad4 & \qquad &  & \qquad\qquad\qquad\qquad12-2+1 &  & \qquad\qquad\qquad\textrm{co-}a_{11}^{(4)}=15\\
\quad5 & \qquad &  & \qquad\qquad\qquad\qquad51-5+2 &  & \qquad\qquad\qquad\textrm{co-}a_{48}^{(5)}=31\\
\quad6 & \qquad &  & \qquad\qquad\qquad\qquad206-20+3 &  & \qquad\qquad\qquad a_{189}^{(6)}=63\\
\quad7 & \qquad &  & \qquad\qquad\qquad\qquad829-81+4 &  & \qquad\qquad\qquad\textrm{co-}a_{752}^{(7)}{{{\scriptstyle {\scriptscriptstyle }}\atop {\textstyle {\scriptstyle {\scriptscriptstyle {\scriptstyle ?}}}}}\atop {{\textstyle =}\atop }}127\\
\quad\vdots & \qquad &  & \qquad\qquad\qquad\qquad\vdots &  & \qquad\qquad\qquad\vdots\end{array}$\vspace{3mm}
\lyxline{\normalsize}\bigskip{}

\end{table}

\noindent

\subsection{\label{sub:The-positional-number}The positional number system aspect
and ``intensional'' Catalan numbers}

\noindent  We have seen the continued-fraction representation (CFR)
of cardiodic arclength is closely bound up with basic properties of
the Catalan structure of $\textrm{LL}(G_{\mu\nu}^{(p)})$. Since we
are also interested in the finer points of Catalan structure -- such
as were hinted at in the CFR $\Phi^{(p)}=(G_{\textrm{max}}^{(p)}/C_{q})^{1/q}$
of Sect. \emph{\ref{sub:Trace-structure-vs.}} --, we are searching
here for specially shaped numbers whose CFR properties would help
examine them side by side with the cardioidic arclength case. As a
starting point, we have a closer look at a conspicuous region where
$\delta_{1}^{(n)}$ first begins constant and then resumes deviating
shortly after: \begin{equation}
\begin{array}{cccccccccc}
{\scriptstyle \underrightarrow{\delta_{1}^{(3)}=0}} & {\atop {\displaystyle a_{1}^{(4)}=5}} & {\scriptstyle \underrightarrow{\delta_{1}^{(4)}=-1}} & {\atop {\displaystyle a_{1}^{(5)}=10}} & {\scriptstyle \underrightarrow{\delta_{1}^{(5)}=-1}} & {\atop {\displaystyle a_{1}^{(6)}=20}} & {\scriptstyle \underrightarrow{\delta_{1}^{(6)}=-1}} & {\atop {\displaystyle a_{1}^{(7)}=40}} & {\scriptstyle \underrightarrow{\delta_{1}^{(7)}=0}} & {\atop {\displaystyle a_{1}^{(8)}=81}}\end{array}{\atop {\displaystyle \cdots}}\label{eq:positional-number-sys}\end{equation}

\noindent The ensuing $a_{1}^{(n)}$ hints at a base-$5\cdot2^{n-4}$
positional number system origin, but the question is, if there are
deviations pending, is the relationship between $a_{1}^{(n)}$ and
$5\cdot2^{n-4}\;(n\geq4)$ stable enough so that $5\cdot2^{n-4}$
could be called their base? The answer is in the affirmative since
the ratio $a_{1}^{(n)}/(5\cdot2^{n-4})$ is fast approaching the constant
$1.01859{\scriptstyle \ldots}$. There is a neat interpretation at
hand for this phenomenon in terms of two rival angular measurement
systems -- the SI system, in which there are $2000\pi$ milliradians
in the circle, and the NATO system, with $6400$ angular mil in the
circle. Their conversion ratio coinciding with $1.01859{\scriptstyle \ldots}$,
we can set \[
\frac{a_{1}^{(n)}}{5\cdot2^{n-4}}\sim\frac{6400}{2000\pi}=\frac{2^{5}}{2\cdot5\cdot\pi}\]
to see that $a_{1}^{(n)}$ is just the integer approximation of $2^{n}/\pi$
and the change in $\delta_{1}^{(n)}$ due to decimal switching $<0.5\;\leftrightarrow\;>0.5$,
\emph{viz}. \[
\begin{array}{cccccc}
 & a_{1}^{(4)}=5 & a_{1}^{(5)}=10\quad & a_{1}^{(6)}=20\quad & a_{1}^{(7)}=40\quad & a_{1}^{(8)}=81\quad\\
\\ & {\scriptstyle {\scriptstyle {\displaystyle \frac{2^{4}}{\pi}}}}=5.09\quad & {\scriptstyle {\displaystyle \frac{2^{5}}{\pi}}}=10.18\quad & {\scriptstyle {\displaystyle \frac{2^{6}}{\pi}}}=20.37\quad & {\scriptstyle {\displaystyle \frac{2^{7}}{\pi}}}=40.74\quad & {\scriptstyle {\displaystyle \frac{2^{8}}{\pi}}}=81.48\quad\end{array}\cdots\]
\[
\]

\noindent  Candidate numbers that would allow taking such features
into account are the special Catalan's $C_{\textrm{-}1/k}\equiv\frac{2^{-2/k}\Gamma(1/2-1/k)}{\sqrt{\pi}\,\Gamma(2-1/k)}$,
\begin{equation}
C_{\textrm{-}1/k}=l_{0}^{(k)}+\frac{1}{l_{1}^{(k)}+{\displaystyle \frac{1}{\; l_{2}^{(k)}{\displaystyle +\frac{1}{l_{3}^{(k)}+{\displaystyle \ddots}}}}}}\;\equiv\;[l_{0}^{(k)};l_{1}^{(k)},l_{2}^{(k)},l_{3}^{(k)},\ldots],\label{eq:contfrac_C-1/k}\end{equation}

\medskip{}

\noindent  which satisfy the successor axiom\begin{equation}
l_{1}^{(k+1)}=l_{0}^{(k)}+l_{1}^{(k)}=1+l_{1}^{(k)}\qquad(k>3),\label{eq:analog-inter}\end{equation}
\medskip{}

\noindent and possess incidences of $l_{\varkappa}^{(k+1)}=1+l_{\varkappa}^{(k)}$
for $\varkappa>1$ where the successor relation breaks after a few
increases in $k$, \emph{viz}.\[
\]
\begin{equation}
l_{4}^{(6)}=1,l_{4}^{(7)}=2,l_{4}^{(8)}=3,l_{4}^{(9)}=4\quad(\textrm{but }l_{4}^{(10)}=7).\label{eq:sameness-breaking}\end{equation}
\medskip{}

\noindent Apparently there is successorhip in example (\ref{eq:sameness-breaking}),
expressed by $k-l_{4}^{(k)}=\textrm{const}.$, which lasts until carry
occurs in the quinary or the decimal system, as the case may be. This
suggests arranging $k$ such that the numbers $C_{\textrm{-}1/k}$
are indexed by pa$\mathbf{\underline{i}}$rs of alter$\mathbf{\underline{n}}$a$\underline{\mathbf{t}}$ing
Mers$\underline{\mathbf{en}}$ne number$\underline{\mathbf{s}}$ ,
ar$\underline{\mathbf{i}}$thmetically-averaged and $\underline{\mathbf{o}}$rganized
as $\mathbf{\underline{n}}$egative reciproc$\underline{\mathbf{al}}$s,
or ``intensional'', for short:\begin{equation}
C_{\textrm{-}1/k}\qquad\textrm{with}\quad k=\frac{p+q}{2}=4,9,19,39,\ldots\qquad\textrm{for}\quad p=7,15,31,63,\ldots,\quad q=(p-3)/4.\label{eq:intense}\end{equation}
\[
\]
\noindent In terms of positional number systems, the first member
of the above $k$ sequence, 4, is about to carry in the quinary system;
the second, 9, is about to carry in the decimal or the quinary system,
and mutatis mutandis for the further members with respect to the vigesimal,
quadragesimal, etc. systems. Thus, $k$ and $k'$, endowed with the
relation\[
k'=2k+1,\]
\[
\]
form carry-digit neighborhoods in all base-$5\cdot2^{n}$ systems,
and the original carry-bit neighborhood of $p$ and $p'$ is recovered
cutting by the rightmost digit, most easily recognizable for the decimal
system:%
\footnote{The naturalness of $k=\frac{p+q}{2}=9,19,\ldots$ can be recognized
by the behavior of the first nontrivial representative sequence $\left(G_{\rho}^{(15)}\right)=(3,5,11,17,41,113)$.
By constructing the tuples $(3,5),(3,11),\ldots,$ $(3,5,11),(3,5,17),\ldots,$
$(3,5,11,17,41,113)$ and checking their scalar products (in absolute
value) with appropriately-sized tuples $(1,1)^{T},(1,-1)^{T},\ldots,$
$(1,1,1)^{T},(1,1,-1)^{T},\ldots,$ $(1,-1,-1,-1,-1,-1)^{T}$, one
sees that the natural numbers in the range 1 to $190=\sum G_{\rho}^{(15)}$
are covered, leaving but an unrepresentable rest of twenty numbers:
$7,34,48,\ldots,189$. \emph{Nineteen} of these can be lifted by adding
the first member of the representative sequence $\left(G_{\rho}^{(31)}\right)=(19,\ldots)$,
adapting the tuples appropriately; but at the expense of -- now \emph{nine}
-- new exceptions $183,191,\ldots208$ within the enlarged (interordinal)
range 1 to 209.%
}\[
19\!\!\!/\quad39\!\!\!/\quad79\!\!\!/\quad159\!\!\!/\quad\cdots\]
\bigskip{}

\noindent As Catalan structure is characterized by the way the Catalan
numbers $C_{q}$ thru to $C_{2q}$ are partitioned in $\textrm{LL}(G_{\mu\nu}^{(p)})$
(or $\textrm{LL}(J_{\mu\nu}^{(p)})$ for that matter), we first have
to search for an algorithm that looks for an intensional-Catalan-number
CFR for $C_{q}=G_{q+2,1}^{(p)}$. Determining intensional-Catalan-number
CFR for the remaining entries $G_{\mu\nu}^{(p)}$ then consists in
further refinement steps. In other words, where $\varphi^{(k)}:\mathbb{N}_{0}\rightarrow\mathbb{N}$
is defined by $\varphi^{(k)}(\varkappa):=l_{\varkappa}^{(k)}\;(k>3)$,
we search for a partial inverse map $(\varphi^{-1})^{(k)}:\mathbb{N}\rightarrow\mathbb{N}$
defined by $(\varphi^{-1})^{(k)}(y)=\varkappa$, for select values
$y$ and begin with $y=G_{q+2,1}^{(p)}$. The first algorithm we propose
embodies an interordinal relationship:
\begin{alg}
\label{alg:intense Cq}Where $j$ and $m$ are natural numbers, pick
the paraorders $p$ and $p'=2p+1$, with $q=(p-3)/4$ and $q'=(p'-3)/4$,
and initialize $j$ with $q$ and $m$ with $\max(2^{q}C_{q},20)$.
Vary $m$ by successive increases or decreases, and if needed reinitialize
$m$ and increase $j$, until for some pair $(j,m)$ and for some
prime $\pi_{s}>2$ the condition (1) $4^{j-1}m-\pi_{s}=G_{3q+5}^{(p)}$
is fulfilled under the constraint (2) $j<q'$. Then $\varkappa=2^{j}m+2^{j}-1+C_{q+1}^{2}$.\end{alg}
\begin{lyxlist}{00.00.0000}
\item [{Case$\: p'=15$:}] This is an example where with any contfrac calculator
we can find $l_{47}^{((9)}=C_{3}=G_{5,1}^{(15)}=5$ and check Algorithm
\ref{alg:intense Cq} for this solution.\\
 Given are $q'=3\rightsquigarrow C_{q'}=5$; $q=1\rightsquigarrow C_{q}=1,\max(2^{q}C_{q},20)=20,C_{q+1}=2$;
and $p=7\rightsquigarrow G_{3q+5}^{(p)}=1$.\\
 With $j=1,m=20$, right from the start we have $\varkappa=20\cdot2+1+2^{2}=47$,
and the pair $(j,m)$ fulfils condition (1) $4^{j-1}m-\pi_{8}=20-19=1$
as well as constraint (2).\\

\item [{Case$\: p'=31$:}] Given are $q'=7\rightsquigarrow C_{q'}=429$;
$q=3\rightsquigarrow C_{q}=5,\max(2^{q}C_{q},20)=40,C_{q+1}=14$;
and $p=15\rightsquigarrow G_{3q+5}^{(p)}=41$.\\
Then $\varkappa=13\cdot2^{4}+2^{4}-1+14^{2}=419$ and the pair
$(j,m)$ fulfils condition (1) $4^{j-1}m-\pi_{5}=52-11=41$ as well
as constraint (2).\\

\item [{Case$\: p'=63$:}] Given are $q'=15\rightsquigarrow C_{15}=9694845$;
$q=7\rightsquigarrow C_{q}=429,\max(2^{q}C_{q},20)=54912,C_{q+1}=1430$;
and $p=31\rightsquigarrow G_{3q+5}^{(p)}=58781$.\\
Then $\varkappa=59764\cdot2^{7}+2^{7}-1+1430^{2}=9694819$ and
the pair $(j,m)$ fulfils condition (1) $4^{j-1}m-\pi_{166}=59764-983=58781$
as well as constraint (2).
\end{lyxlist}
\noindent We have found a second algorithm that delivers identical
results for $p=15,31$, but differs for $p=63$.
\begin{alg}
\label{alg:second intense Cq}Where $\bar{C}(p)$ is the largest even
Catalan number $C_{r}<p\;(p=15,31,\ldots;q=(p-3)/4)$, choose a prime
number $\pi_{6m}$ such that $m\geq q$ is least under the constraint
$\pi_{6m}>C_{q}$. Then $\varkappa=\pi_{6m}-\bar{C}(p)$.\end{alg}
\begin{lyxlist}{00.00.0000}
\item [{Case$\: p=15$:}] Given are $q=3,C_{3}=5,\bar{C}(15)=C_{4}=14$;
then $\varkappa=\pi_{6\cdot3}-14=61-14=47$.\\

\item [{Case$\: p=31$:}] Given are $q=7,C_{7}=429,\bar{C}(31)=C_{4}=14$;
then $\varkappa=\pi_{6\cdot14}-14=433-14=419$.\\

\item [{Case$\: p=63$:}] Given are $q=15,C_{15}=9694845,\bar{C}(63)=C_{5}=42$;
then $\varkappa=\pi_{6\cdot107624}-42=9694877-42=9694835$.
\end{lyxlist}
\noindent Proving one of these algorithms wrong lies beyond the scope
of present-day online computing capabilities yet. In what follows
we stick to $\varkappa\leq500$ to address finer points of intensional-Catalan-number
CFR, and also attempt disambiguating the result for $p=63$. An important
aid in this enterprise is supplied by the $\mathrm{\mathbf{\underline{p}}}$ara$\mathrm{\mathbf{\underline{o}}}$rder
$\underline{\mathrm{\mathbf{s}}}$ums $\Sigma_{i=1}^{n}p_{i}=2^{n+1}-n-2$,
and their $\mathbf{\underline{\mathrm{\mathbf{e}}}}$ntourage

\begin{equation}
\begin{array}{l}
\textrm{pose}_{1}(p)\equiv2^{\Theta_{n}}\Sigma_{i=1}^{n}p_{i}-(\varphi^{-1})^{(k)}(C_{q}),\\
\textrm{pose}_{2}(p)\equiv2^{\Theta_{n}-1}\Sigma_{i=1}^{n+1}p_{i}-(\varphi^{-1})^{(k)}(C_{q}),\end{array}\quad(n>2)\label{eq:pose numbers}\end{equation}
\medskip{}

\noindent where $\Theta_{n}=\mathcal{C}_{n}^{+}+\Sigma_{i=1}^{n}\Delta r_{i}$,
$\mathcal{C}_{n}^{+}$ being the $n$th member of the ordered sequence
$\mathcal{C}^{+}$(see Eq. (\ref{eq:Upsilonplus})) and $\Sigma_{i=2}^{n}\Delta r_{i}$
the sum of index increments $\Delta r_{i}=\left\lceil \left|r_{i}\right|\right\rceil -\left\lceil \left|r_{i-1}\right|\right\rceil $
for $C_{r_{i-1}}=\mathcal{C}_{i-1}^{+}$ and $C_{r_{i}}=\mathcal{C}_{i}^{+}\;(n>2)$. 

\begin{table}[H]
\caption{\label{tab:pose and pose}Paraorder sums and their entourage up to
$n=8$}

\medskip{}
\lyxline{\normalsize}

\noindent $\quad p_{n}$$\qquad\qquad\qquad k$$\qquad\qquad\qquad$$\Sigma_{i=1}^{n}p_{i}\qquad\qquad\qquad\qquad\qquad\textrm{pose}_{1}(p)$$\qquad\qquad\qquad\qquad\qquad\qquad\qquad\qquad\textrm{pose}_{2}(p)$

\lyxline{\normalsize}\vspace{-5mm}
$\begin{array}{cccccc}
 & \qquad\qquad &  & \qquad\qquad\qquad & \qquad\qquad & \qquad\qquad\\
\quad1 & \qquad & - & \qquad\qquad\qquad1 & \qquad\qquad- & \qquad\qquad-\\
\quad3 & \qquad & - & \qquad\qquad\qquad4 & \qquad\qquad- & \qquad\qquad-\\
\quad7 & \qquad & 4 & \qquad\qquad\qquad11 & \qquad\qquad2^{0+0}\cdot11-1=10 & \qquad\qquad2^{-1}\cdot26-1=12\\
\quad15 & \qquad & 9 & \qquad\qquad\qquad26 & \qquad\qquad\!\!2^{1+0}\cdot26-47=5 & \qquad\qquad2^{0}\cdot57-47=10\\
\quad31 & \qquad & 19 & \qquad\qquad\qquad57 & \qquad\qquad\!\!2^{2+1}\cdot57-419=37 & \qquad\qquad2^{2}\cdot120-419=61\\
\quad63 & \qquad & 39 & \qquad\qquad\qquad120 & \qquad\qquad\!\!2^{14+3}\cdot120-(\varphi^{-1})^{(39)}(C_{15}){{{\scriptstyle {\scriptscriptstyle }}\atop {\textstyle {\scriptstyle {\scriptscriptstyle {\scriptstyle ?}}}}}\atop {{\textstyle =}\atop }}6033821 & \qquad\qquad2^{16}\cdot247-(\varphi^{-1})^{(39)}(C_{15}){{{\scriptstyle {\scriptscriptstyle }}\atop {\textstyle {\scriptstyle {\scriptscriptstyle {\scriptstyle ?}}}}}\atop {{\textstyle =}\atop }}6492573\\
\quad127 & \qquad & 79 & \qquad\qquad\qquad247 & \qquad\qquad\!\!2^{42+4}\cdot247-(\varphi^{-1})^{(79)}(C_{31})=\:? & \qquad\qquad2^{45}\cdot502-(\varphi^{-1})^{(79)}(C_{31})=\:?\\
\quad255 & \qquad & 159 & \qquad\qquad\qquad502 & \qquad\qquad\!\!? & \qquad\qquad?\end{array}$ \lyxline{\normalsize} \bigskip{}

\end{table}
\medskip{}
\noindent We recall: $C_{7}=429$ is the constitutive Catalan representative
$G_{9,1}^{(31)}$ of $(G_{\mu\nu}^{(31)})$. It is one of the results
predicted by Algs. \ref{alg:intense Cq} and \ref{alg:second intense Cq}
that this number is matched by the 419th expansion coefficient of
$C_{\textrm{-}1/19}$\begin{equation}
l_{419}^{(19)}=429.\label{eq:first stub}\end{equation}
\[
\]
Out of the remaining $G_{\mu\nu}^{(31)}$, only those that belong
to the nonbracketed, nonparenthesized part of the corridor $G$-set,%
\footnote{for a definition, and the meaning of the parentheses and brackets,
see Sect. \emph{\ref{sub:Row-(column)-structure}}%
} $G_{\textrm{cor}}^{(31)}=\{[1,3,5,11,17,41],(19,43),115,155,429\}$,
are allocated in the vicinity of $\varkappa=419$, which means Eqs.
(\ref{eq:first stub}) and (\ref{eq:sec stub a}) constitute an intensional
Catalan-number CF description of that part of the corridor $G$-set.
Thus, $G_{9,2}^{(31)}=155$ is matched by the 408th expansion coefficient,
and $G_{10,3}^{(31)}=115$ by the 397th,\begin{equation}
l_{408}^{(19)}=155,\qquad l_{397}^{(19)}=115.\label{eq:sec stub a}\end{equation}
\emph{Viz}.\[
\begin{array}{cccccccccccccccccc}
C_{\textrm{-}1/19}= & \left[1;\right. & 16, & 2, & 4, & \ldots, & 115, & \;2, & 13, & \ldots, & 155, & 97, & \;1, & \ldots, & 429, & \;2, & \;4, & \ldots\,].\\
{\atop } & {0\atop } & {1\atop } & {2\atop } & {3\atop } & {\cdots\atop } & {397\atop } & {398\atop } & {399\atop } & {\cdots\atop } & {408\atop } & {409\atop } & {410\atop } & {\cdots\atop } & {419\atop } & {420\atop } & {421\atop } & {\cdots\atop }\end{array}\]
 \noindent It turns out that $G_{9,2}^{(31)}=155$ as a CF denominator
of $C_{\textrm{-}1/19}$ occurs at a distance $11=\Sigma_{i=1}^{3}p_{i}$
off $\varkappa=419$, and the same distance lies between $(\varphi^{-1})^{(19)}(G_{9,2}^{(31)})$
and $(\varphi^{-1})^{(19)}(G_{10,3}^{(31)})$, \begin{equation}
\begin{array}{l}
\\\begin{array}{cccccccccc}
{\atop {\displaystyle \left(l_{419}^{(19)}=429\right)}}{\scriptstyle \underrightarrow{\mathcal{D}^{(19)}=11}} & {\atop {\displaystyle l_{408}^{(19)}=155}} & {\scriptstyle \underrightarrow{\mathcal{D}^{(19)}=11}} & {\atop {\displaystyle l_{397}^{(19)}=115,}}\end{array}{\atop }\\
\\\end{array}\label{eq:regular traversal}\end{equation}

\noindent so that the distance (edge length) between the entries
(nodes), $\mathcal{D}^{(19)}$, in this case coincides with the paraorder
sum \begin{equation}
11=\Sigma_{i=1}^{n-2}p_{i}.\label{eq:distlaw}\end{equation}

\noindent It's interesting to compare this pattern with that corresponding
to the alternating-sign corridor $J$-set of $\textrm{LL}(J_{\mu\nu}^{(31)})$,
although a slightly different methodology is required to this end.
Let $C_{\textrm{-}1/k}$ alternatively be given by the expansion%
\footnote{for further details, see the contfrac-routine options provided by
wims.unice.fr%
} \begin{equation}
C_{\textrm{-}1/k}=\ell_{0}^{(k)}-\frac{1}{\ell_{1}^{(k)}+{\displaystyle \frac{1}{\;\ell_{2}^{(k)}{\displaystyle -\frac{1}{\ell_{3}^{(k)}+\ddots}}}}}\;\equiv\;[\ell_{0}^{(k)};\ell_{1}^{(k)},\ell_{2}^{(k)},\ell_{3}^{(k)},\ldots],\label{eq:alter-intens-Cat}\end{equation}

\noindent and the associated map $\psi^{(k)}:\mathbb{N}_{0}\rightarrow\mathbb{N}$
by $\psi^{(k)}(\varkappa):=\ell_{\varkappa}^{(k)}\;(k>3)$ with partial
inverse $(\psi^{-1})^{(k)}:\mathbb{N}\rightarrow\mathbb{N}$. Then,
for entry $J_{9,1}^{(19)}$ there exists a denominator whose place
$(\psi^{-1})^{(19)}(-J_{9,1}^{(19)})$ lies in the vicinity of the
place $(\varphi^{-1})^{(19)}(G_{9,1}^{(31)})$ predicted by Algs.
\ref{alg:intense Cq} and \ref{alg:second intense Cq}, namely:\[
-J_{9,1}^{(31)}=429=\ell_{438}^{(19)};\]

\noindent but that vicinity by necessity now leads to branchings
for the remaining two entries which do not obey a strict ordering,
$J_{9,2}^{(31)},J_{10,3}^{(31)}\ngtr J_{\max}^{(15)}\,$ instead of
$G_{9,2}^{(31)},G_{10,3}^{(31)}>G_{\max}^{(15)}\,$:

\begin{equation}
J_{9,2}^{(31)}=\psi^{(19)}(\varkappa_{a})\pm\psi^{(19)}(\varkappa_{b}),\quad J_{10,3}^{(31)}=\psi^{(19)}(\varkappa_{c})\pm\psi^{(19)}(\varkappa_{d}).\label{eq:branching}\end{equation}

\noindent Thus, \[
J_{9,2}^{(31)}=117=\ell_{411}^{(19)}+\ell_{409}^{(19)}=116+1,\]
and\[
J_{10,3}^{(31)}=143=\ell_{425}^{(19)}-\ell_{414}^{(19)}=156-13.\]
\[
\]
Including the sublevels created, in contradistinction to the edge
length $11=\Sigma_{i=1}^{3}p_{i}$ of example (\ref{eq:regular traversal}),
the average edge length now equals $((438-425)+(425-414)+(425-411)+(411-409))/4$
$=10=\textrm{pose}_{1}(7)$\emph{:} 

\begin{equation}
\begin{array}{l}
\begin{array}{llllllllll}
{\atop {\displaystyle \left(\ell_{438}^{(19)}=429\right)}} & {\scriptstyle \underrightarrow{\mathcal{D}_{0,c}^{(19)}=13}} & {\atop {\displaystyle \ell_{425}^{(19)}=156}} & {\scriptstyle \underrightarrow{\mathcal{D}_{c,a}^{(19)}=14}} & {\atop {\displaystyle \ell_{411}^{(19)}=116}}\end{array}\\
\quad\qquad\qquad\qquad\begin{array}{cccccccccc}
{\scriptstyle \mathcal{D}_{c,d}^{(19)}=11} & {{\displaystyle \mid}\atop {\displaystyle \downarrow}} & \qquad\qquad & {\scriptstyle \mathcal{D}_{a,b}^{(19)}=2} & {{\displaystyle \mid}\atop {\displaystyle \downarrow}}\end{array}{\atop }\\
\quad\qquad\qquad\qquad\qquad\begin{array}{cccccccccc}
 & {\atop {\displaystyle \enskip\ell_{414}^{(19)}=13}} & \qquad\qquad & {\atop {\displaystyle \ell_{409}^{(19)}=1.}}\end{array}{\atop }\\
{\atop }\end{array}\label{eq:rerversed traversal}\end{equation}

\noindent Now the restriction of $C_{\textrm{-}1/19}$ CF denominators
to those qualifying as representatives of nonbracketed, nonparenthesized
corridor $G$-set entries recalls a similar one of kissing number
representatives to those qualifying as simple interordinal $(\vartheta_{\lambda}^{(p_{l},p_{u})})$
or second-order synoptic differences at paraorder 31: the entries
$G_{9,1}^{(31)}=429$ and $G_{9,2}^{(31)}=155$ are $C_{\textrm{-}1/19}$
CF represented by $l_{419}^{(19)}$ and $l_{408}^{(19)}$, respectively;
their pendants (in the kissing-number representative sense) from Table
\ref{tab:The-first-sixteen} are given by $L_{8}=240=\vartheta_{1}^{(15,255)}=\vartheta_{2}^{(7,127)}$
and $L_{10}=336=\vartheta_{7}^{(15,63)}=\vartheta_{3}^{(15,127)}$.
Also, $G_{10,3}^{(31)}=115$, $C_{\textrm{-}1/19}$ CF represented
by $l_{397}^{(19)}$, has a pendant from Table \ref{tab:Kissing-numbers-from},
$L_{9}=272=\partial G_{10}^{(31)}-\partial J_{8}^{(31)}$. The simple-interordinal/second-order
synoptic difference representability desert following $L_{10}$, first
ending at $L_{13}$ for Table \ref{tab:Kissing-numbers-from}, and
at $L_{16}$ for Table \ref{tab:The-first-sixteen}, should be accompanied
by a similar desert in $C_{\textrm{-}1/19}$ CF representability,
whose discovery yet awaits improved CFR computability conditions.
\medskip{}

\noindent Before looping back to the case $k=9$, let us make it
clear that the ensuing places discussed this far are least, that is,
the $l$- or $\ell$-values they map may reappear at higher places.
Thus the results of Algs. \ref{alg:intense Cq} and \ref{alg:second intense Cq}
for case $p=63$, if meaningful, need not be conflicting: $\varkappa=9694835$
could be a place of recurrence of $C_{15}$, as suggested by $C_{15}-(\varphi^{-1})^{(39)}(C_{15})=\textrm{pose}_{2}(15)$,
while $\varkappa=9694819$ would be the least and supported by $C_{15}-(\varphi^{-1})^{(39)}(C_{15})=\Sigma_{i=1}^{4}p_{i}$.
Keeping this in mind, we can now turn to entries which, unavailable
though they seem for $k=19$, are reachable for $k=9$. Thus, the
place of $G_{8,1}^{(15)}=113$ can be computed as $\textrm{pose}_{1}(31)\cdot\Sigma_{i=1}^{3}p_{i}=37\cdot11$:
\begin{equation}
l_{37\cdot11}^{(9)}=113,\label{eq:Gmax stub}\end{equation}
\emph{viz.} \[
\begin{array}{ccccccccccccccccc}
C_{\textrm{-}1/9}= & \left[1;\right. & 6, & 1, & 1, & 4, & \ldots, & \;5, & \;3, & \,113, & \;1, & \;1, & \ldots & \left.\right]\\
{\atop } & {0\atop } & {1\atop } & {2\atop } & {3\atop } & {4\atop } & {\cdots\atop } & {405\atop } & {406\atop } & {407\atop } & {408\atop } & {409\atop } & {\cdots\atop } & {\atop } & {\atop } & {\atop } & {\atop }\end{array}\]
and there is no lower place than 407 with this property. Yet, there
is another occurrence of $113$, close to the first, \[
l_{414}^{(9)}=113,\]
\emph{viz}. \[
\begin{array}{ccccccccccccccccc}
C_{\textrm{-}1/9}= & \left[1;\right. & 6, & 1, & 1, & 4, & \ldots, & \;6, & \;1, & \,113, & \;1, & \;15, & \ldots & \left.\right]\\
{\atop } & {0\atop } & {1\atop } & {2\atop } & {3\atop } & {4\atop } & {\cdots\atop } & {412\atop } & {413\atop } & {414\atop } & {415\atop } & {416\atop } & {\cdots\atop } & {\atop } & {\atop } & {\atop } & {\atop }\end{array}\]
which falls into place in that $G_{8,1}^{(15)}=113$ recurs interordinally
as (non-corridor-$G$ set entry) $G_{12,5}^{(31)}$=113. Subtracting
$\textrm{pose}_{1}(15)=5$ from place $(\varphi^{-1})_{\textrm{least}}^{(19)}(429)=419$,
we get

\begin{equation}
(\varphi^{-1})_{\textrm{2nd least}}^{(9)}(113)=(\varphi^{-1})_{\textrm{least}}^{(19)}(429)-\textrm{pose}_{1}(15).\label{eq:Gmax alt a}\end{equation}
\[
\]
And there is a third occurrence of 113, doubled in value, and computable
using $(\varphi^{-1})_{\textrm{least}}^{(9)}(5)$, but $\textrm{pose}_{2}(15)=10$:

\begin{equation}
(\varphi^{-1})_{\textrm{least}}^{(19)}(226)=(\varphi^{-1})_{\textrm{least}}^{(9)}(5)-\textrm{pose}_{2}(15),\label{eq:Gmax alt b}\end{equation}
\emph{viz}. \[
\begin{array}{ccccccccccccccccc}
C_{\textrm{-}1/19}= & \left[1;\right. & 16, & 2, & 4, & 1, & \ldots, & \;2, & \;1, & \,226, & \;3, & \;1, & \ldots & \left.\right] & .\\
{\atop } & {0\atop } & {1\atop } & {2\atop } & {3\atop } & {4\atop } & {\cdots\atop } & {35\atop } & {36\atop } & {37\atop } & {38\atop } & {39\atop } & {\cdots\atop } & {\atop } & {\atop } & {\atop } & {\atop }\end{array}\]
So the preliminary interpretation of these observations would read:
if the map $\varphi^{-1}(y_{1})\mapsto\varphi^{-1}(y_{2})-\textrm{pose}(p)$
for key values $y_{1},y_{2}$ is associated with a context change
$k'\mapsto k$ and the subtrahend is $\textrm{pose}_{1}(p)$, the
result is non-minimal $\varkappa$, and, conversely, if this map is
associated with a context change $k\mapsto k'$ and $\textrm{pose}_{2}(p)$
is subtracted, the result is minimal $\varkappa$, but with a doubled
reference outcome. Further study is required to corroborate this point.\medskip{}

\subsection{\label{sub:The-kissing-number}The kissing number aspect revisited}

\noindent Figures linked to Catalan numbers in a fundamental way
like the kissing numbers can be expected to be present in more overt
form in the current framework. They lay hidden in cardiodic-arclength
CFR, where $C_{7}-m_{c}^{(6)}$ yields the eighth kissing number,
$429-189=L_{8}$, and $C_{6}-m_{c}^{(5)}$ equals the third plus the
sixth, $132-48=12+72=L_{3}+L_{6}$. They're also an implicit part
of the workings of our algorithms, where $\left|(\varphi^{-1})_{\textrm{least}}^{(9)}(5)-5\right|=L_{5}+L_{1}$,
$\left|(\varphi^{-1})_{\textrm{least}}^{(19)}(429)-429\right|=L_{3}-L_{1}$,
and $\left|(\varphi^{-1})_{\textrm{least}}^{(39)}(9694845)-9694845\right|=L_{4}+L_{1}$
according to Algorithm \ref{alg:intense Cq}, and $L_{3}-L_{1}$ according
to Algorithm \ref{alg:second intense Cq}. Plus, they led to a salient
interordinal corridor aspect in the previous section. Changing that
perspective of inner regulative to its dual -- exterior connection
of densest-packing hypersphere configurations with Green's parafermions
of Mersennian order such that $L_{D}=\vartheta_{\lambda_{1}}^{(p_{l_{1}},p_{u_{1}})}\pm\vartheta_{\lambda_{2}}^{(p_{l_{2}},p_{u_{2}})}\pm\ldots$
as outlined in Conjecture \ref{con:A-hypersphere-configuration} --,
there is no a priori reason why kissing numbers should not occur overtly
as expansion coefficients of suitably chosen irrationals in imitation
of this connection. This may be tested using the assessable case when
simple interordinal differences $\vartheta_{\lambda}^{(p_{l},p_{u})}$of
Green's squares suffice to represent $L_{D}.$ \newpage{}

\subsubsection{Detuning intensional Catalan numbers}

\noindent From Eq. (\ref{eq:intense}) it follows that $k\sim\frac{5}{8}p$
for large $p$, we might therefore take the integer approximation
of the mean value $\frac{5}{16}(p_{l}+p_{u})$ as a target index $k$
in $C_{\textrm{-}1/k}$ and look for occurrences of $l_{\varkappa_{\textrm{least}}}^{(k)}=\vartheta_{\lambda}^{(p_{l},p_{u})}=L_{D}$,
keeping the pairing that used in Table \ref{tab:The-first-sixteen}:
\bigskip{}
\begin{table}[H]
\caption{\label{tab:kiss-as-exp-1}The first eight kissing numbers $L_{D}$
represented by contfrac expansion coefficients $l_{\varkappa_{\textrm{least}}}^{(k)}$
from $C_{\textrm{-}1/k}$ for $k\sim\frac{5}{16}(p_{l}+p_{u})$}

\medskip{}
\lyxline{\normalsize}

\noindent $\quad D$$\qquad\qquad\qquad\qquad$$\frac{5}{16}(p_{l}+p_{u})\qquad\qquad\qquad\qquad\vartheta_{\lambda}^{(p_{l},p_{u})}$$\enskip\qquad\qquad\qquad\qquad\qquad k$$\quad\qquad\qquad\qquad\qquad\qquad\qquad C_{\textrm{-}1/k}$

\lyxline{\normalsize}\vspace{-5mm}
$\begin{array}{cccccc}
 & \qquad\qquad & \qquad\qquad & \qquad\qquad\qquad\qquad & \qquad\qquad\qquad\qquad & \quad\qquad\qquad\qquad\\
\quad1 & \qquad & \qquad\qquad1.25 & \qquad\qquad\qquad\qquad\vartheta_{1}^{(1,3)}=2 & \qquad\qquad\qquad\qquad4^{*} & \quad\qquad\qquad\qquad l_{3}^{(4)}=2\\
\quad2 & \qquad & \qquad\qquad2.5 & \qquad\qquad\qquad\qquad\vartheta_{1}^{(1,7)}=6 & \qquad\qquad\qquad\qquad4^{*} & \quad\qquad\qquad\qquad l_{45}^{(4)}=6\\
\quad3 & \qquad & \qquad\qquad3.125 & \qquad\qquad\qquad\qquad\vartheta_{3}^{(3,7)}=12 & \qquad\qquad\qquad\qquad4^{*} & \quad\qquad\qquad\qquad l_{42}^{(4)}=12\\
\quad4 & \qquad & \qquad\qquad6.87 & \qquad\qquad\qquad\qquad\vartheta_{3}^{(7,15)}=24 & \qquad\qquad\qquad\qquad\!\!7 & \quad\qquad\qquad\qquad l_{178}^{(7)}=24\\
\quad5 & \qquad & \qquad\qquad6.87 & \qquad\qquad\qquad\qquad\vartheta_{5}^{(7,15)}=40 & \qquad\qquad\qquad\qquad\!\!8 & \quad\qquad\qquad\qquad l_{118}^{(8)}=40\\
\quad6 & \qquad & \qquad\qquad11.87 & \qquad\qquad\qquad\qquad\vartheta_{3}^{(7,31)}=72 & \qquad\qquad\qquad\qquad\!\!11 & \quad\qquad\qquad\qquad l_{151}^{(11)}=72\\
\quad7 & \qquad & \qquad\qquad40 & \qquad\qquad\qquad\qquad\vartheta_{1}^{(1,127)}=126 & \qquad\qquad\qquad\qquad\!\!40 & \qquad\qquad\qquad\qquad\qquad l_{\textrm{n/a}}^{(40)}\text{}^{\dagger}\:(\textrm{but}\, l_{4}^{(42)}=126)\\
\quad8 & \qquad & \qquad\qquad84.37 & \qquad\qquad\qquad\qquad\vartheta_{1}^{(15,255)}=240 & \qquad\qquad\qquad\qquad\!\!84 & \qquad\qquad\qquad\qquad\qquad l_{\textrm{n/a}}^{(84)}\text{}^{\dagger}\:(\textrm{but}\, l_{401}^{(91)}=240)\end{array}$ \lyxline{\normalsize}$^{*}$) the case $k<4$ is outside the domain
of successor relation (\ref{eq:analog-inter})\\
$\text{}^{\dagger}$) not available due to limitation to < 500
contfrac steps \bigskip{}

\end{table}

\bigskip{}
\noindent The meaningfulness of detuning $k$ to $\frac{5}{16}(p_{l}+p_{u})$
is apparently limited: a) for the first three dimensions, the Catalan
numbers $C_{\textrm{-}1/k}$ fall out of the range $\left]1,C_{\textrm{-}1/4}\right]=\left]1,1.57..\right]$
obeyed for finite $k\geq4$: $C_{\textrm{-}1}=-0.5$, $C_{\textrm{-}1/2}=0$,
$C_{\textrm{-}1/3}=0.11..$; and b) for dimensions seven and eight,
the CFRs of $C_{\textrm{-}1/40}$ and $C_{\textrm{-}1/84}$ respectively
fail to include 126 or 240 among their (first 500) denominators. The
basic idea of incorporating $p_{l}$ and $p_{u}$ in the irrationals'
gradation yet seems sound and just calling for a different implementation.\subsubsection {A qphyletic approach}

\noindent What looks more promising is finding ways to exploit the
identity ${\scriptscriptstyle -}\left\lceil {\scriptscriptstyle {\textstyle n/2}}\right\rceil +\Sigma_{i=1}^{n-2}p_{i}=\left\lfloor \log_{2}C_{q}\right\rfloor \;(n>3)$.
Setting $q_{u}=(p_{u}-3)/4$, we may construct irrationals from $\log_{2}C_{q_{u}}$,
$2^{i}/\pi$ and $a_{1}^{(i)}$ for $i=1,2,\ldots,\left\lfloor \log_{2}C_{q_{u}}\right\rfloor $,
and mould them into graded sequences, obvious candidates being\[
\begin{array}{cc}
\left((2^{i}/\pi)^{-1}\left\lfloor \log_{2}C_{q_{u}}\right\rfloor \right),\\
\left((2^{i}/\pi)^{-1}\log_{2}C_{q_{u}}\right), & \qquad\quad\qquad\qquad\qquad\qquad{\scriptstyle {\textstyle {{\displaystyle q_{u}=3,7,15,31,63,127;}\atop {\displaystyle i=1,2,\ldots,\left\lfloor \log_{2}C_{q_{u}}\right\rfloor }.}}}\\
\left((a_{1}^{(i)})^{-1}\log_{2}C_{q_{u}}\right),\end{array}\]

\noindent We further introduce a regularized range $\left]0,1\right[$
for candidates to be admissible, constraining the gradings to\begin{equation}
\begin{array}{cc}
\left((2^{n_{l}}/\pi)^{-1}\left\lfloor \log_{2}C_{q_{u}}\right\rfloor \:\textrm{with regular CFR }{\scriptstyle \daleth}_{\varkappa}^{(n_{l},q_{u})}\right),\\
\left((2^{n_{l}}/\pi)^{-1}\log_{2}C_{q_{u}}\:\textrm{with regular CFR }\mathsf{{\scriptstyle \beth}}_{\lambda}^{(n_{l},q_{u})}\right), & \qquad{\scriptstyle {\textstyle {{\displaystyle q_{u}=7,15,31,63,127;}\atop {\displaystyle n_{l}=\log_{2}(p_{u}{\scriptstyle +{\textstyle 1}}),\log_{2}(p_{u}{\scriptstyle +{\textstyle 1}})+1,\ldots,\left\lfloor \log_{2}C_{q_{u}}\right\rfloor }.}}}\\
\left((a_{1}^{(n_{l})})^{-1}\log_{2}C_{q_{u}}\:\textrm{with regular CFR }\mathsf{{\scriptstyle \left(\daleth_{\mathit{A}}\right)}}_{\mu}^{(n_{l},q_{u})}\right),\end{array}\label{eq:regularized grades}\end{equation}
\noindent The results of this program are summarized in the table
below and compared to those of Table \ref{tab:The-first-sixteen}
(column labeled $\Lambda$ approach):

\begin{table}[H]
\caption{\label{tab:regularized irrationals}The first sixteen kissing numbers
$L_{D}$ represented by contfrac expansion coefficients ${\scriptstyle \daleth}_{\varkappa_{\textrm{least}}}^{(n_{l},q_{u})},{\scriptstyle \,\beth}_{\lambda_{\textrm{least}}}^{(n_{l},q_{u})}\,\textrm{and}\,{\scriptstyle \left(\daleth_{A}\right)}_{\mu_{\textrm{least}}}^{(n_{l},q_{u})}$ }

\noindent \medskip{}
\lyxline{\normalsize}

$\; D$$\quad\qquad\qquad\qquad$$L_{D}\textrm{ }(\Lambda\textrm{ approach)}\qquad\qquad$
$L_{D}\,\textrm{as }\,{\scriptstyle \daleth}_{\varkappa_{\textrm{least}}}^{(n_{l},q_{u})}\qquad\qquad$
$L_{D}\,\textrm{as }\,\mathsf{{\scriptstyle \beth}}_{\lambda_{\textrm{least}}}^{(n_{l},q_{u})}\qquad\qquad$
$L_{D}\,\textrm{as}\:\mathsf{{\scriptstyle \left(\daleth_{\mathrm{\mathit{A}}}\right)}}_{\mu_{\textrm{least}}}^{(n_{l},q_{u})}$

\lyxline{\normalsize}\medskip{}
\begin{tabular}{cccccccccccccccc}
1 & \hspace*{1cm} &  & $\vartheta_{1}^{(1,3)}=2$ & \hspace*{1cm} &  & ${\scriptstyle \daleth}_{7}^{(5,7)}$ & \hspace*{1cm} &  & ${\scriptstyle \beth}_{9}^{(5,7)}$ & \hspace*{1cm} &  & ${\scriptstyle \left(\daleth_{A}\right)}_{5}^{(5,7)}$ &  &  & \tabularnewline
2 &  &  & $\vartheta_{1}^{(1,7)}=6$ &  &  & ${\scriptstyle \daleth}_{66}^{(5,7)}$ &  &  & ${\scriptstyle \beth}_{2}^{(5,7)}$ &  &  & ${\scriptstyle \left(\daleth_{A}\right)}_{2}^{(5,7)}$ &  &  & \tabularnewline
3 &  &  & $\vartheta_{3}^{(3,7)}=12$ &  &  & ${\scriptstyle \daleth}_{412}^{(5,7)}$ &  &  & ${\scriptstyle \beth}_{103}^{(5,7)}$ &  &  & ${\scriptstyle \left(\daleth_{A}\right)}_{276}^{(5,7)}$ ${}^\dagger$ &  &  & \tabularnewline
4 &  &  & $\vartheta_{3}^{(7,15)}=24$ &  &  & ${\scriptstyle \daleth}_{39}^{(5,7)}$ &  &  & ${\scriptstyle \beth}_{239}^{(5,7)}$ ${}^\dagger$ &  &  & ${\scriptstyle \left(\daleth_{A}\right)}_{31}^{(5,7)}$ &  &  & \tabularnewline
5 &  &  & $\vartheta_{5}^{(7,15)}=40$ &  &  & ${\scriptstyle \daleth}_{77}^{(5,7)}$ &  &  & ${\scriptstyle \beth}_{152}^{(5,7)}$ &  &  & ${\scriptstyle \left(\daleth_{A}\right)}_{76}^{(5,7)}$ &  &  & \tabularnewline
6 &  &  & $\vartheta_{3}^{(7,31)}=72$ &  &  & ${\scriptstyle \daleth}_{8}^{(5,7)}$ &  &  & ${\scriptstyle \beth}_{55}^{(9,7)}$ &  &  & ${\scriptstyle \left(\daleth_{A}\right)}_{279}^{(11,7)}$ ${}^\dagger$ &  &  & \tabularnewline
7 &  &  & $\vartheta_{1}^{(1,127)}=126$ &  &  & ${\scriptstyle \daleth}_{164}^{(8,7)}$ &  &  & ${\scriptstyle \beth}_{63}^{(14,15)}$ &  &  & ${\scriptstyle \left(\daleth_{A}\right)}_{246}^{(16,15)}$ &  &  & \tabularnewline
8 &  &  & $\vartheta_{1}^{(15,255)}=240$ &  &  & ${\scriptstyle \daleth}_{201}^{(15,15)}$ &  &  & ${\scriptstyle \beth}_{398}^{(10,31)}$ ${}^*$ &  &  & n/a &  &  & \tabularnewline
9 &  &  & $\mathrm{{\scriptstyle 2nd\, o.p.e.}}=272$ &  &  & n/a &  &  & n/a &  &  & ${\scriptstyle \left(\daleth_{A}\right)}_{156}^{(8,7)}$ &  &  & \tabularnewline
10 &  &  & $\vartheta_{7}^{(15,63)}=336$ &  &  & n/a &  &  & ${\scriptstyle \beth}_{360}^{(18,15)}$ &  &  & n/a &  &  & \tabularnewline
11 &  &  & $\mathrm{{\scriptstyle 2nd\, o.p.e.}}=438$ &  &  & n/a &  &  & ${\scriptstyle \beth}_{27}^{(17,31)}$ &  &  & n/a &  &  & \tabularnewline
12 &  &  & $\mathrm{{\scriptstyle 2nd\, o.p.e.}}=648{}^{\diamond}$ &  &  & ${\scriptstyle \daleth}_{407}^{(8,7)}$ ${}^*$ &  &  & ${\scriptstyle \beth}_{236}^{(32,127)}$ ${}^\dagger$ &  &  & n/a &  &  & \tabularnewline
13 &  &  & $\mathrm{{\scriptstyle 2nd\, o.p.e.}}=906{}^{\diamond}$ &  &  & ${\scriptstyle \daleth}_{1}^{(16,15)}$ &  &  & n/a &  &  & n/a &  &  & \tabularnewline
14 &  &  & $\mathrm{{\scriptstyle 2nd\, o.p.e.}}=1422$ &  &  & n/a &  &  & n/a &  &  & n/a &  &  & \tabularnewline
15 &  &  & $\mathrm{{\scriptstyle 2nd\, o.p.e.}}=2340$ &  &  & n/a &  &  & n/a &  &  & n/a &  &  & \tabularnewline
16 &  &  & $\vartheta_{9}^{(31,511)}=4320$ &  &  & n/a &  &  & n/a &  &  & n/a &  &  & \tabularnewline
\end{tabular}

\noindent \lyxline{\normalsize}

\noindent \noindent${}^*)$ these denominators are related to $l_{408}^{(19)},l_{397}^{(19)}$ (see comparison of Diagrams (\ref{eq:regular traversal}) and (\ref{eq:just coincidence?}) in main text)\par \noindent${}^\dagger)$ these denominators form the $L_D$ divisor hierarchy up to $L_{12}$ (see discussion of Diagram (\ref{eq:or this a coincidence?}) in  main text)\par \noindent${}^{\diamond})$ these are nonlattice kissing numbers, their lattice counterparts being $L_{12}=756$ and $L_{13}=918$ respectively\bigskip{}

\end{table}
\bigskip{}
\noindent For dimensions one to eight, in the chosen range, all three
types of irrationals have denominators that cover (nearly) all of
the corresponding kissing numbers. While for dimensions nine and thirteen,
the irrationals $(a_{1}^{(8)})^{-1}\log_{2}C_{7}$ and $(2^{16}/\pi)^{-1}\left\lfloor \log_{2}C_{15}\right\rfloor $
supplement the picture with denominators respectively representing
$L_{9}=272$ and $L_{13}=906$, the two of them last seen in the synoptic
Table \ref{tab:Kissing-numbers-from}. The representations of $L_{11}=438$
and $L_{12}=648$ are entirely new additions. The phenomenon of apparent
non-representability of $L_{D}$ for dimensions fourteen to sixteen
in the chosen range could prove to  be real, but the artificial limitation
to 500 contfrac steps, accounted for by a ``n/a'' in the table, casts
doubts on that conclusion -- at least $L_{16}=4320$ should show up
somewhere, judging from parafermial and synoptic representability
considerations.\[
\]

\noindent The strength of some of the figures in Table \ref{tab:regularized irrationals}
can be acknowledged on the basis of the kinship to base-$5\cdot2^{n-4}$
positional number systems they reveal in much the same way as intensional
Catalan numbers do. We can recognize a similarity between structure
diagram (\ref{eq:regular traversal}) and \begin{equation}
\begin{array}{l}
{\scriptscriptstyle }\\
\begin{array}{cccccccccc}
{\atop {\displaystyle \left(L_{11}=438\right)}}{\scriptstyle \underrightarrow{\mathcal{D}=31}} & {\atop {\displaystyle {\scriptstyle \daleth}_{407}^{(8,7)}=648}} & {\scriptstyle \underrightarrow{\mathcal{D}=9}} & {\atop {\displaystyle {\scriptstyle \beth}_{398}^{(10,31)}=240,}}\end{array}{\atop }\\
{\scriptscriptstyle }\end{array}\label{eq:just coincidence?}\end{equation}
which has its origin in the branching\begin{equation}
l_{1}^{(19)}={\scriptstyle \daleth}_{1}^{(8,7)}+{\scriptstyle \beth}_{1}^{(10,31)}.\label{eq:ao-decomp}\end{equation}
\[
\]
\noindent While the average edge length in Dgs. (\ref{eq:regular traversal})
and (\ref{eq:rerversed traversal}) equals 11 and 10, with a margin
reflecting the swings $408\leftrightarrow407$, $397\leftrightarrow398$,
in Dg. (\ref{eq:just coincidence?}) it equals $(31+9)/2=20$, and
in the following hierarchy diagram of $L_{D}$ divisors up to $L_{12}=648$
it amounts to $(43+40+40+37)/4=40$:

\begin{equation}
\begin{array}{l}
\begin{array}{llllllllll}
 & \quad & {\atop {\displaystyle {\scriptstyle \beth}_{236}^{(32,127)}=648}} & {\scriptstyle \underleftarrow{\mathcal{D}=43}} & {\atop {\displaystyle {\scriptstyle \left(\daleth_{A}\right)}_{279}^{(11,7)}=72}}\end{array}\\
\begin{array}{cccccccccc}
 & {\scriptstyle \mathcal{D}=40} & {{\displaystyle \,{}_{{\displaystyle \uparrow\,}}}\atop ^{{\displaystyle {\scriptstyle ^{{\displaystyle {\textstyle {\displaystyle \mid}}}}}}}} &  & \:\qquad\qquad & {\scriptstyle \mathcal{D}=40} & {{{\atop {\displaystyle \vphantom{}}}\atop {\displaystyle \mid}}\atop {{\displaystyle \downarrow}\atop {\scriptscriptstyle }}}\end{array}\\
\begin{array}{cccccccccc}
 &  & {{\scriptscriptstyle }\atop {\displaystyle {\scriptstyle \quad\left(\daleth_{A}\right)}_{276}^{(5,7)}=12}} & \;{\scriptstyle {\scriptstyle \underrightarrow{\mathcal{D}=37}}} & {{\scriptscriptstyle }\atop {\displaystyle {\displaystyle {\scriptstyle \enskip\beth}_{239}^{(5,7)}=24.}}} &  & {\displaystyle {\scriptstyle }}\end{array}\\
{\atop }\end{array}\label{eq:or this a coincidence?}\end{equation}
\noindent What is novel about these diagrams is that they seem to
impart a terminologic shift on the labels intra- and interordinality.
Borrowing from systematic biology, we may quasi identify structures
arising in $l^{(k)}$ (or in $\ell^{(k)}$) with ``phyla'' , and Dgs.
(\ref{eq:regular traversal}) and (\ref{eq:rerversed traversal})
phyleticly intraordinal, accordingly. By contrast, Dgs. (\ref{eq:just coincidence?})
and (\ref{eq:or this a coincidence?}), by displaying parallelism
between ${\scriptstyle \,\beth}_{\lambda}^{(n_{l},q_{u})}\,\textrm{and}\,{\scriptstyle \left(\daleth_{A}\right)}_{\mu}^{(n_{l},q_{u})}$,
are phyleticly interordinal. As the latter radiate along distinct
$n_{l}$ and/or $q_{u}$, we could, in contradistinction, call them
qphyletic, qphyla becoming the name for superordinate structures fed
from both ${\scriptstyle \,\beth}_{\lambda}^{(n_{l},q_{u})}\,\textrm{and}\,{\scriptstyle \left(\daleth_{A}\right)}_{\mu}^{(n_{l},q_{u})}$.\\
 \bigskip{}

\noindent We have found that Eq. (\ref{eq:ao-decomp}) is only one
of three possible solutions to the ``ancestral branching'' aspect

\begin{equation}
l_{1}^{(p+q)/2}={\scriptstyle \daleth}_{1}^{(n_{1},q)}+{\scriptstyle \beth}_{1}^{(n_{2},p)},\label{eq:aspect-oriented}\end{equation}
\[
\]

\bigskip{}

\noindent the other two being\[
l_{1}^{(9)}={\scriptstyle \daleth}_{1}^{(5,3)}+{\scriptstyle \beth}_{1}^{(7,15)}\]

\noindent and

\[
l_{1}^{(39)}={\scriptstyle \daleth}_{1}^{(10,15)}+{\scriptstyle \beth}_{1}^{(13,63)}.\]
\medskip{}

\noindent The latter is strong -- and independent -- evidence that
in the neighborhood of $(\varphi^{-1})_{\textrm{least}}^{(39)}(9694845)$
analogs of Dgs. (\ref{eq:regular traversal}), (\ref{eq:rerversed traversal})
and (\ref{eq:just coincidence?}) could be hiding, which requires
for its verification knowing those irrationals to a precision that
would allow for a computation scope of roughly $10^{7}$ contfrac
denominators.


\vspace{1.8cm}

\begin{thebibliography}{ShuKo10}
\bibitem[Green98]{Green98}H.S.\,Green, Int.J.of Theor.Phys. \textbf{37}
No.11, 2735 (1998)

\bibitem[Green53]{Green53}H.S.\,Green, Phys.Rev. \textbf{90}, 270
(1953)

\bibitem[Merk89]{Merk89}U.\,Merkel, Nuov.Cim.\textbf{ 103B} No.6,
599 (1989) 

\bibitem[ALee76]{ALee76}A.\,Lee, Period.Math.Hungar. \textbf{7},
63 (1976)

\bibitem[Nebe11]{Nebe11}G.\,Nebe, Dichte Kugelpackungen, in: \emph{Facettenreiche
Mathematik}, ed. A.\,Werner, K.\,Wendland, Vieweg (2011)

\bibitem[Green11]{Greenberg11}O.W.\,Greenberg, arXiv:1102.0764v1
{[}physics.hist-ph{]} 

\bibitem[ShuKo10]{ShuKo10}J.A.\,Shuster,\,J.\,Köplinger, Appl.Math.Comp.\textbf{
216} No.12, 3497 (2010)
\end{thebibliography}
\end{document}